\documentclass[reqno]{amsart}
\usepackage{amssymb,wasysym,mathscinet,enumitem,mathtools}
\usepackage{mathrsfs}   

\usepackage[numbers,sort&compress]{natbib}
\usepackage{hyperref} 
\usepackage{esint} 

\numberwithin{equation}{section}

\theoremstyle{plain}
\newtheorem{theorem}{Theorem}[section]
\newtheorem{lemma}[theorem]{Lemma}
\newtheorem{corollary}[theorem]{Corollary}
\newtheorem{proposition}[theorem]{Proposition}

\theoremstyle{definition}
\newtheorem{definition}[theorem]{Definition}

\newtheorem{assumption}[theorem]{Assumption}

\theoremstyle{remark}
\newtheorem{remark}[theorem]{Remark}
\newtheorem*{remark*}{Remark}

 \makeatletter
 \newcommand{\norm}{\@ifstar{\@normb}{\@normi}}
 \newcommand{\@normb}[2]{\left\Vert{#1}\right\Vert_{#2}}
 \newcommand{\@normi}[2]{\Vert{#1}\Vert_{#2}}
 \makeatother

 \global\long\def\Sob#1#2{{W}^{#1}_{#2}} 
 \global\long\def\dSob#1#2{{\mathcal{H}}^{#1}_{#2}} 
  
 \global\long\def\nSob#1#2{\mathbb{H}^{#1}_{#2}}

 \global\long\def\Sobloc#1#2{W^{#1}_{#2,\mathrm{loc}}} 
 \global\long\def\dSobloc#1#2{\mathcal{H}^{#1}_{#2,\mathrm{loc}}} 
 \global\long\def\oSob#1#2{\mathring{W}^{#1}_{#2}}
 \global\long\def\Leb#1{L_{#1}} 
 \global\long\def\Lebloc#1{L_{#1,\mathrm{loc}}}

 \newcommand{\boldF}{\mathbf{F}} 
 \newcommand{\boldG}{\mathbf{G}} 
 \newcommand{\boldH}{\mathbf{H}}

   \DeclareMathOperator{\Div}{div}
 \newcommand{\relphantom}[1]{\mathrel{\phantom{#1}}}
 
 \newcommand{\myd}[1]{\,d{#1}}

 \DeclareMathOperator{\supp}{supp}

 \global\long\def\Sobd#1#2{{W}^{#1}_{#2,0}}
 \global\long\def\odSob#1#2{\mathring{\mathcal{H}}^{#1}_{#2}}

\usepackage{tcolorbox}

\allowdisplaybreaks

\makeatletter
\def\@tocline#1#2#3#4#5#6#7{\relax
  \ifnum #1>\c@tocdepth % then omit
  \else
    \par \addpenalty\@secpenalty\addvspace{#2}%
    \begingroup \hyphenpenalty\@M
    \@ifempty{#4}{%
      \@tempdima\csname r@tocindent\number#1\endcsname\relax
    }{%
      \@tempdima#4\relax
    }%
    \parindent\z@ \leftskip#3\relax \advance\leftskip\@tempdima\relax
    \rightskip\@pnumwidth plus4em \parfillskip-\@pnumwidth
    #5\leavevmode\hskip-\@tempdima
      \ifcase #1
       \or\or \hskip 1em \or \hskip 2em \else \hskip 3em \fi%
      #6\nobreak\relax
    \hfill\hbox to\@pnumwidth{\@tocpagenum{#7}}\par% <---- \dotfill -> \hfill
    \nobreak
    \endgroup
  \fi}
\makeatother

\begin{document}

\title[Mixed norm estimates for nonstationary Stokes equations]{Interior and boundary mixed norm derivative estimates for nonstationary Stokes equations}

\author[H. Dong]{Hongjie Dong}
\address{Division of Applied Mathematics, Brown University, 182 George Street, Providence, RI 02912, USA}
\email{hongjie\_dong@brown.edu }

\author[H. Kwon]{Hyunwoo Kwon}
\address{Division of Applied Mathematics, Brown University, 182 George Street, Providence, RI 02912, USA}
\email{hyunwoo\_kwon@brown.edu }
\thanks{H. Dong was partially supported by Simons Fellows Award 007638 and the NSF under agreement DMS-2055244. H. Kwon was partially supported by the NSF under agreement DMS-2055244.}

\begin{abstract}
We obtain weighted mixed norm Sobolev estimates in the whole space for nonstationary Stokes equations in divergence and nondivergence form with variable viscosity coefficients that are merely measurable in time variable and have small mean oscillation in spatial variables in small cylinders. As an application, we prove interior mixed norm derivative estimates for solutions to both equations. We also discuss boundary mixed norm Hessian estimates for solutions to equations in nondivergence form under the Lions boundary conditions.
\end{abstract}

\maketitle

%\tableofcontents

\section{Introduction}

This paper is devoted to studying weighted mixed-norm estimates, and interior and boundary mixed-norm derivative estimates for solutions to nonstationary Stokes equations with variable viscosity coefficients. We consider Stokes equations in nondivergence form:
\begin{equation}\label{eq:Stokes-nondiv}
\left\{
\begin{alignedat}{1}
\partial_t u-a^{ij}(t,x)D_{ij} u+\nabla p &=f\\
\Div u &=g
\end{alignedat}
\right.\,\,\text{in } U.
\end{equation}
Here $U$ is a cylindrical domain in $\mathbb{R}^{d+1}$, $d\geq 2$,  $u:U \rightarrow \mathbb{R}^d$ denotes the velocity field, {$p:U \rightarrow \mathbb{R}$} denotes the associated pressure, and $f$ and $g$ are a given vector field and a function defined on $U$. The variable viscosity coefficients $a^{ij}$ satisfy the following uniform ellipticity conditions: there exists $\nu \in (0,1)$ such that
\begin{equation}\label{eq:elliptic}
\nu|\xi|^2 \leq a^{ij}(t,x) \xi_i\xi_j,\quad |a^{ij}(t,x)|\leq \nu^{-1},\quad \text{for all } i,j\in \{1,\dots,d\}
\end{equation}
for all $\xi=(\xi_1,\dots,\xi_d)\in \mathbb{R}^d$ and $(t,x) \in \mathbb{R}^{d+1}$. We follow the Einstein summation convention for repeated indices.

We also  consider Stokes equations in divergence form:
\begin{equation}\label{eq:Stokes-div}
\left\{
\begin{alignedat}{1}
\partial_t u-D_i(a^{ij}(t,x)D_{j} u)+\nabla p &=\Div \boldF\\
\Div u &=g\\
\end{alignedat}
\right.\,\,\text{in } U,
\end{equation}
where $\boldF=(F^{ij})_{1\leq i,j\leq d}:U\rightarrow \mathbb{R}^{d\times d}$ is a 2-tensor and $\Div \boldF$ is a vector field defined by
\[   \Div \boldF =\left( D_j F^{1j},\dots, D_j F^{dj}\right). \]

Besides mathematical interests, Equations \eqref{eq:Stokes-nondiv} and \eqref{eq:Stokes-div} can be used to model non-Newtonian fluids that have thixotropy, i.e., time-dependent shear thinning property (see e.g. \cite{BGL16}). These equations are also naturally introduced when we consider Stokes equations on manifolds (see e.g. \cite{MT01,DM04}).

{Since the pioneering work of Calder\'on-Zygmund \cite{CZ52}, many researchers have established $\Leb{q}$-estimates for elliptic and parabolic equations with variable coefficients, possibly discontinuous in spatial or time variables. See e.g. \cite{D20} and their references therein. In this paper, we mainly review results on $\Leb{q}$-estimates for Stokes equations.}

For stationary Stokes equations, there is plenty of literature on {Sobolev type} estimates. When the viscosity coefficient is constant, Cattabriga \cite{C61} first obtained $\Sob{1}{q}$-estimates when $d=3$ and $1<q<\infty$ in a smooth domain. Later, it was extended by Amrouche-Girault \cite{AG91} to a bounded $C^{1,1}$-domain, $d\geq 2$ and $1<q<\infty$. This was further extended to bounded Lipschitz domains with small Lipschitz constants by Galdi-Simader-Sohr \cite{GSS94}. A complete solvability result was obtained by Dindo\v{s}-Mitrea \cite{DM04} on arbitrary bounded Lipschitz domain in $\mathbb{R}^d$, $d\geq 2$.  We refer to \cite{G11} for exterior problems of stationary Stokes equations. When viscosity coefficients are variable coefficients,  Dong-Kim \cite{DK18,DK19} obtained $\Sob{1}{q}$-estimates and weighted $\Sob{1}{q}$-estimates on Reifenberg flat domains even if the viscosity coefficient is merely measurable in one direction and has a small BMO seminorm in orthogonal directions.

Many authors have studied mixed-norm Sobolev estimates for nonstationary Stokes equations in various settings. When $a^{ij}=\delta^{ij}$, Solonnikov \cite{S73} obtained $\Leb{q}$-estimates and solvability result for \eqref{eq:Stokes-nondiv} under the Dirichlet boundary conditions on the half-space and bounded $C^2$-domains.  Later,  it was extended by Giga-Sohr \cite{GS91} to mixed-norm Sobolev estimates including exterior domains.  An elementary proof was given by Maremonti-Solonnikov \cite{MS95}, and later Geissert et. al. \cite{GMHS10}  gave a different proof via $H^\infty$-calculus. For the problem \eqref{eq:Stokes-div} under the Dirichlet boundary conditions, Giga-Giga-Sohr \cite{GGS93} obtained $\Leb{q}$-estimates on half-spaces without estimating pressure. Later, Koch-Solonnikov \cite{KS02} gave more precise $\Leb{q}$-estimates for the problem \eqref{eq:Stokes-div} on {the half-space} including estimates for the pressure. These results were later extended by Chang-Kang \cite{CK18} to anisotropic Sobolev spaces on the half-space under the Dirichlet boundary conditions.  For weighted estimates, Fr\"ohlich \cite{F07} obtained weighted mixed-norm estimates by employing $H^\infty$-calculus approach based on the Stokes resolvent estimates due to  Farwig-Sohr \cite{FS97} and Fr\"ohlich \cite{F03}.

For variable coefficients cases, there are {relatively} few results on mixed-norm Sobolev estimates. Solonnikov \cite{S01} first obtained $\Leb{q}$-estimates and solvability results for such problem when $a^{ij}$ is continuous in $t$ and belongs to $\Sob{1}{r}$ in $x$ for some $r$ when the domain is bounded. Later, Abels-Terasawa \cite{AT10} and Abels \cite{A10} extended this result to mixed-norm estimates on several unbounded domains. There are also results on $\Leb{q}$-estimates under different assumptions on the viscosity part. See Bothe-Pr\"u{\ss} \cite{BP07}, Pr\"u{\ss} \cite{P18},   Pr\"u{\ss}-Simonett \cite{PS16}, and the references therein.  We also note that the variable density case was considered by Ladyzhenskaya-Solonnikov \cite{LS75} and Danchin \cite{D06}.

{ We now consider interior estimates for Stokes equation.} In the case of the heat equation  $\partial_t v-\Delta v=0$,
it is {well known that}
\[    \norm{D^2v}{\Leb{2}(Q_{1/2})}\leq N\norm{v}{\Leb{2}(Q_1)}\]
for some constant $N=N(d)>0$. Here $Q_r(t_0,x_0)$ denotes the parabolic cylinder centered at $(t_0,x_0)\in\mathbb{R}^{d+1}$ with radius $r>0$:
\[
      Q_r(t_0,x_0) = (t_0-r^2,t_0)\times B_r(x_0),
\]
where $B_r(x_0)$ is the ball in $\mathbb{R}^d$ of radius $r$ centered at $x_0\in\mathbb{R}^d$. When $(t_0,x_0)=(0,0)$, we drop $(t_0,x_0)$ in the notation. However, it is nontrivial to show the validity of such estimates for nonstationary Stokes equations {because of the pressure and the zero divergence constraint.}

When $a^{ij}=\delta^{ij}$, Chen-Strain-Yau-Tsai \cite{CSTY08} proved that if $1<s,q<\infty$, $f\in \Leb{s,q}(Q_1)^d$, $g=0$, and $u\in \Leb{s,1}(Q_1)^d$ is a very weak solution to \eqref{eq:Stokes-nondiv} in $Q_1$, then $D^2 u \in \Leb{s,q}(Q_{1/2})$ and
\[
    \norm{D^2 u}{\Leb{s,q}(Q_{1/2})}\leq N \left(\norm{u}{\Leb{s,1}(Q_1)}+\norm{f}{\Leb{s,q}(Q_1)} \right)
\]
for some constant $N=N(d,s,q)>0$.  This inequality was independently proved by Jin \cite{J13} and Wolf \cite{W15} when $f=g=0$ and $s=q=2$. We also note that Hu-Li-Wang \cite{HLW14} obtained interior $\Leb{q}$-estimates via a different approach without using the representation formula for Stokes equations.  Recently, Dong-Phan \cite{DP21} obtained such estimates even if $a^{ij}$ is not constant and $\Div u =g$. More precisely, if $1<s,q<\infty$,  $f\in\Leb{s,q}(Q_1)^d$, $g\in\Sob{0,1}{s,q}(Q_1)$, and $(u,p) \in \tilde{W}^{1,2}_{s,q}(Q_1)^d\times \Sob{0,1}{1}(Q_1)$ is a strong solution to \eqref{eq:Stokes-nondiv} in $Q_1$, then under the assumption that $a^{ij}$ has small mean oscillation in spatial variables in small cylinders (see Assumption \ref{assump:VMO}), they proved that there exists a constant $N=N(d,s,q,\nu,R_0)>0$ such that
\begin{equation}\label{eq:caccioppoli}
    \norm{D^2 u}{\Leb{s,q}(Q_{1/2})}\leq N \left( \norm{u}{\Leb{s,1}(Q_1)}+\norm{f}{\Leb{s,q}(Q_1)}+\norm{Dg}{\Leb{s,q}(Q_1)}\right).
\end{equation}
Here $\tilde{W}^{1,2}_{s,q}(Q_1)$ is the space of all functions $u$ belonging to $D^k u \in \Leb{s,q}(Q_1)$, $k=0,1,2$, and $u_t \in \Leb{1}(Q_1)$ (see the lines above Theorem \ref{thm:B} for the definition of $\tilde{W}^{1,2}_{s,q}(Q_1)$).
Similarly, gradient estimates were obtained for the problem \eqref{eq:Stokes-div} even if $a^{ij}$ is unbounded. Note that these are only {\em a priori} estimates{, not a local regularity estimate.} In the same paper, they applied interior regularity results for \eqref{eq:Stokes-div} to the incompressible Navier-Stokes equations to improve known regularity criteria results. Very recently, via level set argument as in \cite{HLW14}, Dong-Li \cite{DL22} obtained interior $\Leb{q}$-regularity for Stokes equations in both divergence form and nondivergence form under the stronger assumption that the viscosity coefficients are H\"older continuous in spatial variables.

{  For boundary estimates, Seregin \cite{S00,S09} proved the local spatial smoothing property of strong solutions to nonstationary Stokes equations under the Dirichlet boundary conditions (or no-slip boundary conditions) and $\partial_t u,D^2 u, \nabla p \in \Leb{s,q}(Q_1^+)$, where $Q_r^+=Q_r \cap \mathbb{R}^d_+$. Later, several counterexamples were constructed to show that it is not possible to have spatial smoothing of such solutions under the Dirichlet boundary conditions if we do not impose regularity conditions on the pressure (see Kang \cite{K05} and Seregin-\v{S}ver\'ak \cite{SS10}). Related to our paper, Chang-Kang \cite{CK20} proved that boundary gradient estimates may fail for solutions to nonstationary Stokes equations under the Dirichlet boundary conditions. It is natural to ask what type of boundary conditions may yield the boundary derivative estimates of solutions to nonstationary Stokes equations. See also the review of Seregin-Shilkin \cite{SS14}.

One answer was given by Dong-Kim-Phan \cite{DKP22} who proved that boundary mixed-norm Hessian estimates for solutions to \eqref{eq:Stokes-nondiv} on $Q_1^+$ are possible if we consider the Lions boundary conditions (see \eqref{eq:Lions-boundary}) which were introduced by J.-L. Lions in \cite[pp. 87--98]{L69} (see also P.-L. Lions in \cite[pp. 129--131]{L98}). Such boundary conditions are a special case of the Navier boundary conditions which were introduced by Navier in 1827:
\begin{equation}\label{eq:Navier-boundary}
   u\cdot n=0\quad\text{and}\quad (2\mathbb{D}(u)n)_\tau+\alpha u_\tau =0\quad \text{on } \partial\Omega,
\end{equation}
where $\alpha \geq 0$ is the friction coefficient, $n$ is the outer unit normal vector to the boundary $\partial\Omega$, and $v_\tau = v - (v\cdot n)n$ is the tangential component of $v$ to the boundary $\partial\Omega$, and $\mathbb{D}(u)$ is the deformation tensor of $u$ defined by $[\mathbb{D}(u)]^{ij}=(D_i u^j + D_j u^i)/2$. Many researchers studied Stokes and Navier-Stokes equations under such boundary conditions for mathematical reasons and physical applications. See, for instance, \cite{AAE17,CQ10,CMS20,K06,GK12,LS72,GS80,BB07} and references therein. Very recently, Chen-Liang-Tsai \cite{CLT23} proved that gradient estimates for very weak solutions to nonstationary Stokes equations on $Q_1^+$ are possible under the Navier boundary conditions \eqref{eq:Navier-boundary} when $\Div u=0$ and $a^{ij}=\delta^{ij}$.}

The purpose of this paper is two-fold. We prove weighted mixed-norm Sobolev estimates and solvability of the Cauchy problems for \eqref{eq:Stokes-nondiv} and \eqref{eq:Stokes-div} in $(0,T)\times\mathbb{R}^d$ when the viscosity coefficients satisfy the $\mathrm{VMO}_x$ assumption (see Assumption \ref{assump:VMO}). As an application of these weighted mixed-norm estimates, we prove that if $(u,p)\in \tilde{W}^{1,2}_{q_0}(Q_1)^d\times \Sob{0,1}{1}(Q_1)$ is a strong solution to \eqref{eq:Stokes-nondiv} for some $1<q_0<\infty$ and $f\in \Leb{s,q}(Q_1)^d$, and $g\in \Sob{0,1}{s,q}(Q_1)$, then $D^2 u \in \Leb{s,q}(Q_{1/2})$ and \eqref{eq:caccioppoli} holds. For Stokes equations in nondivergence form, we also prove {boundary Hessian estimates} under the Lions boundary conditions. In contrast to Dong-Phan \cite{DP21} and Dong-Kim-Phan \cite{DKP22}, we do not a priori assume that our strong solution $u$ to \eqref{eq:Stokes-nondiv} belongs to $\tilde{W}^{1,2}_{s,q}$. A similar result holds for weak solutions $u$ to \eqref{eq:Stokes-div} for the interior case.

Let us briefly {outline the} proofs of main theorems. To prove the weighted mixed-norm Sobolev estimates (Theorem \ref{thm:A}) in $(0,T)\times\mathbb{R}^d$, we employ the perturbation technique utilizing the Fefferman-Stein theorem, which was first introduced by Krylov \cite{K07} (see also \cite{K08}). To do so, we need weighted mixed-norm Sobolev estimates for Stokes equations with measurable coefficients depending only on $t$ (Theorem \ref{thm:mixed-simple}), which {are} not available in the literature. Such coefficients are referred to as \emph{simple coefficients} in this article. To obtain the solvability, we consider the associated vorticity equation to remove the pressure term, and then we recover a solution using the divergence equation and the Newtonian potential. A proof is given in \ref{app:B}. Using this solvability result, we prove mean oscillation estimate of the gradient of the vorticity of a solution to \eqref{eq:Stokes-nondiv} to derive a priori estimates for solutions to \eqref{eq:Stokes-nondiv} by using generalized Fefferman-Stein theorem {established in} \cite{DK18-weight} (see Lemma \ref{lem:FS}). Then the desired result follows from the method of continuity {together} with the solvability results for Stokes equations with simple coefficients. A similar argument is also {applied to} Stokes equations in divergence form (Theorem \ref{thm:C}) with some modification.

To prove the interior {mixed-norm} Hessian estimates (Theorem \ref{thm:B}) of solutions to equations in nondivergence form, we mollify Equation \eqref{eq:Stokes-nondiv} in space and time to obtain
\begin{equation*}
\left\{
\begin{aligned}
    \partial_t u^{(\varepsilon)}-a^{ij}D_{ij} u^{(\varepsilon)}+\nabla p^{(\varepsilon)}&=f^{(\varepsilon)}+[a^{ij}D_{ij} u]^{(\varepsilon)}-a^{ij}D_{ij} u^{(\varepsilon)},\\
    \Div u^{(\varepsilon)}&=g^{(\varepsilon)}
\end{aligned}
\right.
\end{equation*}
and then decompose $u^{(\varepsilon)}=u_1^{\varepsilon}+u_2^{\varepsilon}$ and $p^{(\varepsilon)}=p_1^\varepsilon+p_2^\varepsilon$, where $(u_1^{\varepsilon},p_1^\varepsilon)$ satisfies the initial value problem for \eqref{eq:Stokes-nondiv} with {$u_1^\varepsilon(-1,\cdot)=0$} on $\mathbb{R}^d$ by replacing $f$ with $h^\varepsilon:=([a^{ij}D_{ij} u]^{(\varepsilon)}-a^{ij}D_{ij} u^{(\varepsilon)})1_{Q_{3/4}}$ and $g$ with zero, respectively. Using the aforementioned weighted solvability results, we will show that $u_1^\varepsilon \in \Sob{1,2}{s_1,q_1}((-1,0)\times\mathbb{R}^d)$ (see Lemma \ref{lem:regularity-Stokes}) for any $1<s_1,q_1<\infty$. Moreover, it follows from parabolic Sobolev embedding theorem  $\Sob{1,2}{q_0}(Q_1)\hookrightarrow \Leb{s,1}(Q_1)$ and $\Leb{q_0}$-estimates for $u_1^\varepsilon$  that $u_1^\varepsilon \rightarrow 0$ in $\Leb{s,1}(Q_1)$. Then we can apply the result of Dong-Phan \cite{DP21} mentioned above to $u_2^\varepsilon$ to get \eqref{eq:caccioppoli} by replacing $(u,f,g)$ with $(u_2^\varepsilon,f^{(\varepsilon)},g^{(\varepsilon)})$. Then using weak compactness result in $\Leb{s,q}(Q_{1/2})$, we can pass the limit to show that up to subsequence, $D^2 u^{\varepsilon_j}_2\rightarrow D^2 u$ weakly in $\Leb{s,q}(Q_{1/2})$. This implies the desired result in Theorem \ref{thm:B}.

To prove the interior gradient estimates (Theorem \ref{thm:D}) of solutions to equations in divergence form, we perform a similar strategy as in the case of equations in nondivergence form. However, the previous strategy cannot be directly applied since unlike the space $\Sob{1,2}{q_0}(Q_1)$, the space $\dSob{1}{q_0}(Q_1)$  is not always embedded into $\Leb{s,1}(Q_1)$ (see Section \ref{sec:2} for definitions of $\dSob{1}{q_0}(Q_1)$). To overcome this issue, if $s>q_0$, then since $p_1^\varepsilon$, $Du_1^\varepsilon \in \Leb{q_0}((-1,0)\times\mathbb{R}^d)$, $(u_1^\varepsilon)_t$ can be written as a divergence of some matrix field $\mathbf{G}^\varepsilon \in \Leb{q_0}(Q_{3/4})^{d\times d}$. Then by using the recent embedding result due to Kim-Ryu-Woo \cite{KRW22} (see Lemma \ref{lem:KRW}) and $\Leb{q_0}$-estimates for $u_1^\varepsilon$, there exists $q_0<s_1\leq s$ such that $u_1^\varepsilon \in \Leb{s_1,q_0}(Q_{3/4})$ and $u_1^\varepsilon \rightarrow 0$ in $\Leb{s_1,q_0}(Q_{3/4})$ as $\varepsilon\rightarrow 0$. Hence by using a similar argument that we used in the case of nondivergence form, we can show that $Du\in \Leb{s_1,q}(Q_{3/4})$. Then by applying the above argument again, we can prove that $Du \in \Leb{s,q}(Q_{1/2})$ and the corresponding estimate for $Du$ to \eqref{eq:caccioppoli}. The case $s\leq q_0$ is easy to prove. This outlines the proof of Theorems \ref{thm:B} and \ref{thm:D}.

Lastly, this approach also enables us to show the boundary mixed-norm Hessian estimates of strong solutions to \eqref{eq:Stokes-nondiv} in $Q_1^+:=(-1,0)\times \{ y : |y|<1, y_d >0\}$ under the Lions boundary conditions. See Section \ref{sec:8}. However, we mainly focus on the interior derivative estimates for simplicity.

{This paper proceeds in eight sections and three appendix sections.}  In Section \ref{sec:2}, we introduce some notation and state the main results of this paper. In Section \ref{sec:prelim}, we summarize known results on function spaces with and without weights, potential estimates, and solvability results on the divergence equation and parabolic equations with simple coefficients. In Section \ref{sec:4}, we derive solvability results in weighted mixed-norm Sobolev estimates and H\"older estimates for solutions to \eqref{eq:Stokes-nondiv} and \eqref{eq:Stokes-div} with simple coefficients in $(0,T)\times\mathbb{R}^d$. Then we prove weighted mixed-norm solvability results for \eqref{eq:Stokes-nondiv} and \eqref{eq:Stokes-div} in $(0,T)\times\mathbb{R}^d$ with variable viscosity coefficients in Sections \ref{sec:5} and \ref{sec:6}, respectively. In Section \ref{sec:7}, we prove the interior mixed-norm derivative estimates (Theorems \ref{thm:B} and \ref{thm:D}) for solutions to \eqref{eq:Stokes-nondiv} and \eqref{eq:Stokes-div}, respectively. In Section \ref{sec:8}, we give a brief description of proving boundary mixed-norm Hessian estimates for solutions to \eqref{eq:Stokes-nondiv} under the Lions boundary conditions.  Finally, we give {the} proofs of the solvability of Stokes equations with simple coefficients in mixed-norm weighted Sobolev spaces in \ref{app:B} and \ref{app:C}, respectively.

\section{Notation and Main results}\label{sec:2}

\subsection{Notation and assumptions}
By $N=N(p_1,\dots,p_k)$, we denote a generic positive constant depending only on the parameters $p_1,\dots,p_k$. For two Banach spaces $X$ and $Y$, we write $X\hookrightarrow Y$ if $X\subset Y$ and there exists a constant $N$ such that $\norm{u}{Y}\leq N\norm{u}{X}$ for all $u\in X$.

Let  $\Omega$ be any domain in $\mathbb{R}^d$, where $\mathbb{R}^d$ is the standard $d$-dimensional Euclidean space of points $x=(x_1,\dots,x_d)$, $d\geq 2$. For $0<T< \infty$, we write  $\mathbb{R}^d_T:=(0,T) \times \mathbb{R}^d$. We denote the point in $\mathbb{R}^d_T$ by $(t,x)=(t,x',x_d)$, where $x'\in \mathbb{R}^{d-1}$ and $x_d \in \mathbb{R}$. We also define $\mathbb{R}^d_+ := \{ (y',y_d) : y' \in \mathbb{R}^{d-1}, y_d>0\}$.

 For $r>0$ and $(t,x)\in \mathbb{R}^{d+1}$, we write
\[
Q_r(t,x):=(t-r^2,t)\times B_r(x),\quad Q_r=Q_r(0,0)
\]
where
\[   B_r(x):=\{ y \in \mathbb{R}^d : |x-y|<r\}.\]
{For $(t,x) \in \overline{\mathbb{R}^d_+}$,} we define $Q_r^+(t,x)=Q_r(t,x) \cap \mathbb{R}^d_+$ and we write $B_r'(x')$ the $(d-1)$-dimensional ball in $\mathbb{R}^{d-1}$ with the radius $r$ centered at $x' \in \mathbb{R}^{d-1}$.

Let $\mathbb{N}_0=\{0,1,2,\dots\}$ be the set of nonnegative integers. For multi-indices $\gamma=(\gamma_1,\dots,\gamma_d)\in \mathbb{N}^d$ and a function $u$, we define
\[   u_{x_i}=\frac{\partial u}{\partial x_i} = D_i u,\quad (1\leq i\leq d),\quad D^\gamma u =D_1^{\gamma_1}\cdots D_d^{\gamma_d} u,\quad x^\gamma=(x_1)^{\gamma_1}\cdots (x_d)^{\gamma_d}.\]
For $m \in \mathbb{N}$, we use $D^m$ to denote a partial derivative of order $m$ with respect to $x$. For a function $u$, we define
\[
\nabla u :=(D_1 u,\dots,D_d u)\quad\text{and}\quad \nabla^2 u:=[D_{ij}u]_{i,j=1}^d.
\]
Given a weakly differentiable vector field $u=(u^1,\dots,u^d)$, define its \emph{gradient} $\nabla u$ and \emph{vorticity} $\nabla \times u$ by
\[  (\nabla u)^{ij}:= D_j u^i,\quad \text{and}\quad {[\nabla \times u]_{ij}:= D_j u^i - D_i u^j},\quad 1\leq i,j\leq d,\]
respectively.

We use bold-roman to denote $2$-tensors, e.g., $\boldF:(0,T)\times\mathbb{R}^d\rightarrow\mathbb{R}^{d\times d}$. For two vectors $u=(u^1,\dots,u^d)$ and $v=(v^1,\dots,v^d)$, their inner product is defined by
\[   u\cdot v :=\sum_{i=1}^d u^i v^i.\]
For two $2$-tensors $\boldF=[F^{ij}]^d_{i,j=1}$ and $\boldG=[G^{ij}]_{i,j=1}^d$, their inner product is defined by
\[    \boldF : \boldG := \sum_{i,j=1}^d F^{ij} G^{ij}.\]

For a measurable set $A$ of $\mathbb{R}^d$, we use $|A|$ to denote the Lebesgue measure of $A$ and $1_A$ the indicator of $A$. If $0<|A|<\infty$, we write
\[   \fint_A f dx =(f)_{A}:=\frac{1}{|A|} \int_A f \myd{x}.\]

A function $w$ is a \emph{weight} on $\mathbb{R}^d$ if $w$ is nonnegative and $w>0$ a.e. on $\mathbb{R}^d$. For $1<q<\infty$, we write $w\in A_q(\mathbb{R}^d,dx)$ if
\[ [w]_{A_q(\mathbb{R}^d,dx)} :=\sup_{x_0\in \mathbb{R}^d, r>0} \left(\fint_{B_r(x_0)} w \myd{x} \right)\left(\fint_{B_r(x_0)} w^{-1/(q-1)}\myd{x}\right)^{q-1} <\infty.
\]
See basic properties of $A_q$-weights in Subsection \ref{subsec:weight}. We can also define $A_1$ weights, see e.g. \cite[Chapter 7]{G08}. For $k=1,2,\ldots$, $1\leq q<\infty$, and $w\in A_q(\mathbb{R}^d,dx)$, we define
\[
    \Sob{k}{q,w}(\Omega)=\{ u : u, Du, \dots, D^k u \in\Leb{q,w}(\Omega)\}.
\]
By $C_0^\infty(U)$, we denote the set of infinitely differentiable functions with compact support in $U$. For $-\infty<S<T<\infty$, we write
\[  C_0^\infty([S,T)\times\Omega):=\{ u|_{(S,T)\times\Omega} : u\in C_0^\infty((-\infty,T)\times\Omega) \}.\]
 We denote $\Sobd{1}{q,w}(\Omega)$ the closure of $C_0^\infty(\Omega)$ under $\norm{\cdot}{\Sob{1}{q,w}(\Omega)}$.

{For $1< s,q<\infty$, $-\infty\leq S<\infty$, $-\infty<T\leq \infty$, and a weight $w(t,x)=w_1(x)w_2(t)$ on $\mathbb{R}^{d+1}$, where $w_1 \in A_q(\mathbb{R}^d,dx)$ and $w_2\in A_s(\mathbb{R},dt)$, we define}
\begin{align*}
\norm{f}{\Leb{s,q,w}((S,T)\times\Omega)}&:=\left(\int_S^T\left(\int_\Omega |f|^q w_1 \myd{x}\right)^{s/q} w_2 \myd{t}\right)^{1/s}
\end{align*}
and
\[
  \Leb{s,q,w}((S,T)\times\Omega):=\{ f : \norm{f}{\Leb{s,q,w}((S,T)\times\Omega)}<\infty \}.
\]
Similarly, for  $1\leq s,q<\infty$, and $w(t,x)=w_1(x)w_2(t)$, where $w_1\in A_q(\mathbb{R}^d,dx)$, $w_2\in A_s(\mathbb{R},dt)$, we define weighted parabolic Sobolev spaces
\begin{align*}
{\Sob{0,1}{s,q,w}((S,T)\times\Omega)}&{:=\{ u : u, Du \in \Leb{s,q,w}((S,T)\times \Omega)\},}\\
    \Sob{1,2}{s,q,w}((S,T)\times\Omega)&:=\{ u : u, Du, D^2 u, u_t \in\Leb{s,q,w}((S,T)\times\Omega)\}
\end{align*}
with the norm
\begin{align*}
\norm{u}{\Sob{0,1}{s,q,w}((S,T)\times\Omega)}&:=\norm{u}{\Leb{s,q,w}((S,T)\times\Omega)}+\norm{Du}{\Leb{s,q,w}((S,T)\times\Omega)},\\
   \norm{u}{\Sob{1,2}{s,q,w}((S,T)\times\Omega)}&:=\norm{u_t}{\Leb{s,q,w}((S, T)\times\Omega)}+\sum_{k=0}^2\norm{D^k u}{\Leb{s,q,w}((S,T)\times\Omega)}.
\end{align*}
When $s=q$ and $w=1$, we write $\Leb{q}((S,T)\times\Omega)=\Leb{q,q,w}((S,T)\times\Omega)$ and $\Sob{1,2}{q}((S,T)\times\Omega)=\Sob{1,2}{q,q,w}((S,T)\times\Omega)$.  For a measurable function $u$ defined on $(S,T)\times \Omega$, we write $u\in \Lebloc{s,q}((S,T)\times\Omega)$ if $u\in \Leb{s,q}(K)$ for any compact subset $K$ of $(S,T)\times\Omega$. Similarly, we can define {$\Sobloc{0,1}{s,q}((S,T)\times\Omega)$} and $\Sobloc{1,2}{s,q}((S,T)\times\Omega)$.

For equations in divergence form, we introduce additional function spaces  $\nSob{-1}{s,q,w}$ and $\dSob{1}{s,q,w}$. We say that $f \in \nSob{-1}{s,q,w}((S,T)\times\Omega)$ if there exist $g_0,g=(g_1,\dots,g_d) \in \Leb{s,q,w}((S,T)\times\Omega)$ such that
\[   f=g_0+D_i g_i \quad \text{in } (S,T)\times\Omega \]
in the sense of distribution and the norm
\[
     \norm{f}{\nSob{-1}{s,q,w}((S,T)\times\Omega)}:=\inf\left\{ \sum_{i=0}^d \norm{g_i}{\Leb{s,q,w}((S,T)\times\Omega)} : f=g_0+D_i g_i \right\}
\]
is finite. We define
\[
\dSob{1}{s,q,w}((S,T)\times\Omega):=\{ u : u_t \in \nSob{-1}{s,q,w}((S,T)\times\Omega), u \in\Sob{0,1}{s,q,w}((S,T)\times\Omega) \}
\]
with the norm
\[
   \norm{u}{\dSob{1}{s,q,w}((S,T)\times\Omega)}:=\norm{u_t}{\nSob{-1}{s,q,w}((S,T)\times\Omega)}+\norm{u}{\Sob{0,1}{s,q,w}((S,T)\times\Omega)}.
\]
When $s=q$ and $w=1$, we write  $\dSob{1}{q}((S,T)\times\Omega)=\dSob{1}{q,q,w}((S,T)\times\Omega)$.

Now we define strong solutions of Stokes equations in nondivergence form \eqref{eq:Stokes-nondiv}.

\begin{definition}
Let $f \in \Lebloc{1}((S,T)\times\Omega)^d$ and $g \in \Lebloc{1}((S,T)\times\Omega)$. A pair $(u,p)$ is said to be a \emph{strong solution} to \eqref{eq:Stokes-nondiv} in $(S,T)\times\Omega$ if $u \in \Sobloc{1,2}{1}((S,T)\times \Omega)^d$ and $p \in \Sobloc{0,1}{1}((S,T)\times\Omega)$ satisfy
\[
      \partial_t u-a^{ij}D_{ij}u+\nabla p = f\quad\text{and}\quad\Div u=g\quad\text{a.e. in } (S,T)\times\Omega.
\]
\end{definition}

{Similarly,} we define weak solutions of Stokes equations in divergence form \eqref{eq:Stokes-div}.

\begin{definition}
Given $\boldF \in \Lebloc{1}((S,T)\times\Omega)^{d\times d}$  and $g \in \Lebloc{1}((S,T)\times\Omega)$,  $u$  is a \emph{weak solution} to \eqref{eq:Stokes-div} in $(S,T)\times\Omega$ if {$u\in \Sobloc{0,1}{1}((S,T)\times\Omega)^d$ satisfies}
\begin{equation}\label{eq:div-g}
    \int_\Omega u(t,x)\cdot\nabla \varphi(x) \myd{x}=-\int_\Omega g(t,x)\varphi(x)dx
\end{equation}
for a.e.  $t \in (S,T)$, for all  $\varphi \in C_0^\infty(\Omega)$, and
\begin{equation*}%\label{eq:weak-sol-definition}
-\int_S^T\int_\Omega u\cdot(\partial_t \phi)- D_i \phi \cdot a^{ij}D_j u\myd{x}dt=-\int_{S}^T \int_\Omega \boldF : \nabla \phi \myd{x}dt
\end{equation*}
for all $\phi \in C_0^\infty((S,T)\times\Omega)^d$ with $\Div \phi(t)=0$ for all $t\in (S,T)$.
\end{definition}

To discuss the solvability of the initial value problem for Stokes equations in divergence form and nondivergence form in $(S,T)\times \Omega$, we write $u\in \oSob{1,2}{s,q,w}((S,T)\times\Omega)$ if there exists $\tilde{u} \in \Sob{1,2}{s,q,w}((-\infty,T)\times\Omega)$ such that $\tilde{u}=u$ in $(S,T)\times\Omega$ and $\tilde{u}=0$ in $(-\infty,S)\times\Omega$. Similarly, we can define $\odSob{1}{s,q,w}((S,T)\times\Omega)$.

\begin{definition}\label{defn:weak-sol-initial-data} Let $1<s,q<\infty$ and $w\in A_{s,q}$.
\begin{enumerate}
\item[(i)] Given $\boldF \in \Lebloc{1}((S,T)\times\Omega)^{d\times d}$ and $g \in \Lebloc{1}((S,T)\times\Omega)$, we say that $(u,p) \in \dSob{1}{s,q,w}((S,T)\times\Omega)^d\times\Leb{s,q,w}((S,T)\times\Omega)$ is a \emph{weak solution {pair} to \eqref{eq:Stokes-div} in $(S,T)\times\Omega$ with $u(S,\cdot)=0$ on $\Omega$} if $u\in \odSob{1}{s,q,w}((S,T)\times\Omega)^d$ and $(u,p)$ satisfies
\[
-\int_S^T\int_\Omega u\cdot(\partial_t \phi)- D_i \phi \cdot a^{ij}D_j u+p \Div \phi\myd{x}dt=-\int_{S}^T \int_\Omega \boldF : \nabla \phi \myd{x}dt
\]
for all $\phi \in C_0^\infty([S,T)\times\Omega)^d$ {and $u$ satisfies \eqref{eq:div-g}}.
\item[(ii)] Given $f \in \Lebloc{1}((S,T)\times\Omega)^{d}$ and $g \in \Sobloc{0,1}{1}((S,T)\times\Omega)$, we say that $(u,p) \in \Sob{1,2}{s,q,w}((S,T)\times\Omega)^d\times \Sobloc{0,1}{1}((S,T)\times\Omega)$ is a \emph{strong solution to \eqref{eq:Stokes-nondiv} in $(S,T)\times\Omega$  with $u(S,\cdot)=0$ on $\Omega$} if $u\in\oSob{1,2}{s,q,w}((S,T)\times\Omega)^d$ and $u$ is a strong solution to \eqref{eq:Stokes-nondiv} in $(S,T)\times\Omega$.
\end{enumerate}
\end{definition}

 \subsection{Main results}

 Now we present the main results of this paper. The following is our assumption on the viscosity coefficient of Equations \eqref{eq:Stokes-nondiv} and \eqref{eq:Stokes-div}.

{
\begin{assumption}[$\delta$]\label{assump:VMO}
There exists $R_0 \in (0,1/4)$ such that for any $(t_0,x_0) \in \mathbb{R}^{d+1}$ and $0<r<R_0$, there exists $\overline{a}^{ij}(t)$ satisfying \eqref{eq:elliptic} and
\[	\fint_{Q_r(t_0,x_0)} |a^{ij}(t,x)-\overline{a}^{ij}(t)|\myd{x}dt \leq \delta \]
for all $i,j=1,2,\dots, d$.
\end{assumption}
}
\begin{remark}\leavevmode
\begin{enumerate}
\item[(i)] The condition is weaker than the usual full $\mathrm{VMO}$ condition in both $t$ and $x$ since it does not require any regularity condition in $t$. A typical example is $a^{ij}(t,x)=b(t)c^{ij}(x)$, where $b(t)$ and $c(t)$ satisfy
\[
    \nu\leq |b(t)|,|c^{ij}(x)|\leq \nu^{-1},\quad \text{for all } (t,x) \in \mathbb{R}^{d+1},\quad c^{ij} \in \mathrm{VMO} \quad \text{for all } i,j
\]
for some $\nu \in (0,1)$. Here $c^{ij} \in \mathrm{VMO}$ means
\[    \lim_{r\rightarrow 0+} \fint_{B_r(x)} |c^{ij}(y)-(c^{ij})_{B_r(x)}| \myd{y}=0.\]
\item[(ii)]  {By Assumption \ref{assump:VMO} $(\delta)$, there exists $R_0>0$ such that for $0<r<R_0$ and $(t,x_0) \in \mathbb{R}\times \mathbb{R}^d$, we have
\[
    \fint_{Q_r^+(t_0,x_0)} |a^{ij}-\overline{a}^{ij}(t)|\myd{x}dt\leq 2\delta
\]
and
\[
    \fint_{Q_r^+(t_0,x_0)} |a^{ij}-(a^{ij})_{B_r^+(x_0)}|\myd{x}dt\leq 4\delta.
\]}
\end{enumerate} \end{remark}

Our first result concerns the solvability of the initial-value problem for Stokes equations in nondivergence form on weighted mixed-norm Sobolev spaces on $\mathbb{R}^d_T$.

 \begin{theorem}\label{thm:A}
 Let $1<s,q<\infty$, $0<T<\infty$, and let $K_0\geq 1$ be constant, {$w=w_1(x)w_2(t)$, where $[w_1]_{A_q(\mathbb{R}^d,dx)}\leq K_0$ and $[w_2]_{A_s(\mathbb{R},dt)}\leq K_0$.} There exists $0<\delta<1$ depending only on $d$, $\nu$, $s$, $q$,  and $K_0$ such that under Assumption \ref{assump:VMO} $(\delta)$, for every $f\in \Leb{s,q,w}(\mathbb{R}^d_T)^d$ and $g\in\odSob{1}{s,q,w}(\mathbb{R}^d_T)$ and $g_t=\Div G$ for some vector field $G =(G_1,\dots,G_d)\in \Leb{s,q,w}(\mathbb{R}^d_T)^d$ in the sense that
\begin{equation}\label{eq:compatibility-nondiv}
\int_{\mathbb{R}^d_T} g\varphi_t \myd{x}dt=\int_{\mathbb{R}^d_T} G \cdot \nabla \varphi \myd{x}dt
\end{equation}
for any $\varphi \in C_0^\infty([0,T)\times\mathbb{R}^d)$, there exists a unique { $u$ and a unique $p$ upto additive constants such that $(u,p)$ is a strong solution to \eqref{eq:Stokes-nondiv} in $\mathbb{R}^d_T$ with $u(0,\cdot)=0$ on $\mathbb{R}^d$} satisfying
\[
u\in\oSob{1,2}{s,q,w}(\mathbb{R}^d_T)^d,\quad \nabla p \in\Leb{s,q,w}(\mathbb{R}^d_T)^d.
\]
 Moreover, we have
\[
  \norm{u}{\Sob{1,2}{s,q,w}(\mathbb{R}^d_T)}+\norm{\nabla p}{\Leb{s,q,w}(\mathbb{R}^d_T)}\leq N \left(\norm{f}{\Leb{s,q}(\mathbb{R}^d_T)}+\norm{Dg}{\Leb{s,q,w}(\mathbb{R}^d_T)}+\norm{G}{\Leb{s,q,w}(\mathbb{R}^d_T)}\right),
\]
where $N=N(d,s,q,K_0,\nu,R_0,T)>0$.
  \end{theorem}

  The second result describes the solvability of the initial-value problem for Stokes equations in divergence form on weighted mixed-norm Sobolev spaces on $\mathbb{R}^d_T$.

   \begin{theorem}\label{thm:C}
    Let $1<s,q<\infty$, $0<T<\infty$, and let $K_0\geq 1$ be constant, {$w=w_1(x)w_2(t)$, where $[w_1]_{A_q(\mathbb{R}^d,dx)}\leq K_0$ and $[w_2]_{A_s(\mathbb{R},dt)}\leq K_0$.}  There exists $0<\delta<1$ depending only on $d$, $\nu$, $s$, $q$,  and $K_0$ such that under Assumption \ref{assump:VMO} $(\delta)$, for every $\boldF\in \Leb{s,q,w}(\mathbb{R}^d_T)^{d\times d}$ and $g\in \Leb{s,q,w}(\mathbb{R}^d_T)$ satisfying $g_t=\Div\Div \boldG$ for some 2-tensor $\mathbf{G} \in \Leb{s,q,w}(\mathbb{R}^d_T)^{d\times d}$ in the sense that
\begin{equation}\label{eq:compatibility-div}
 \int_{\mathbb{R}^d_T} g \varphi_t \myd{x}dt=-\int_{\mathbb{R}^d_T} \boldG:\nabla^2 \varphi \myd{x}dt
\end{equation}
for all $\varphi \in C_0^\infty([0,T)\times\mathbb{R}^d)$, there exists a unique weak solution $(u,p)$   to \eqref{eq:Stokes-div} in $\mathbb{R}^d_T$ with $u(0,\cdot)=0$ on $\mathbb{R}^d$ satisfying
\[
u\in\odSob{1}{s,q,w}(\mathbb{R}^d_T)^d,\quad  p \in\Leb{s,q,w}(\mathbb{R}^d_T).
\]
Moreover, we have
\[
  \norm{u}{\dSob{1}{s,q,w}(\mathbb{R}^d_T)}+\norm{p}{\Leb{s,q,w}(\mathbb{R}^d_T)}\leq N \left(\norm{\boldF}{\Leb{s,q,w}(\mathbb{R}^d_T)}+\norm{g}{\Leb{s,q,w}(\mathbb{R}^d_T)}+\norm{\boldG}{\Leb{s,q,w}(\mathbb{R}^d_T)}\right),
\]
where $N=N(d,s,q,K_0,\nu,R_0,T)>0$.
   \end{theorem}

  As an application of Theorems \ref{thm:A} and \ref{thm:C}, we prove the interior mixed-norm derivative estimates for strong solutions and weak solutions of \eqref{eq:Stokes-nondiv} and \eqref{eq:Stokes-div}, respectively. To state results in a more compact way, we introduce additional function space
  \begin{align*}
     \tilde{W}_{s,q}^{1,2}(U)&=\{ u : u, D u, D^2 u \in \Leb{s,q}(U), u_t \in \Leb{1}(U)\},
  \end{align*}
where $U$ is an open subset of $\mathbb{R}^{d+1}$.

   \begin{theorem}\label{thm:B}
   Let $1<q_0,s,q<\infty$. Then there exists $\delta=\delta(d,s,q,q_0,\nu)>0$ such that under Assumption \ref{assump:VMO} $(\delta)$, if $(u,p)\in \tilde{W}^{1,2}_{q_0}(Q_{1})^d\times \Sob{0,1}{1}(Q_1)$ is a strong solution to \eqref{eq:Stokes-nondiv} in $Q_{1}$ for some $f\in \Leb{s,q}(Q_{1})^d$ and $g\in \Sob{0,1}{s,q}(Q_{1})$, then $D^2 u  \in \Leb{s,q}(Q_{1/2})$. Moreover, there exists a constant $N=N(d,s,q,q_0,\nu,R_0)>0$ such that
   \[
      \norm{D^2 u}{\Leb{s,q}(Q_{1/2})}\leq N \left(\norm{u}{\Leb{s,1}(Q_{1})}+\norm{f}{\Leb{s,q}(Q_{1})}+\norm{Dg}{\Leb{s,q}(Q_{1})} \right).
   \]
   \end{theorem}

\begin{theorem}\label{thm:D}
  Let $1<q_0,s,q<\infty$. Then there exists $\delta=\delta(d,s,q,q_0,\nu)>0$ such that under Assumption \ref{assump:VMO} $(\delta)$ that if $u \in \Sob{0,1}{q_0}(Q_{1})^d$ is a weak solution to \eqref{eq:Stokes-div} in $Q_{1}$ for some $\boldF\in \Leb{s,q}(Q_{1})^{d\times d}$ and $g\in \Leb{s,q}(Q_{1})$, then $D u  \in \Leb{s,q}(Q_{1/2})$. Moreover, there exists a constant $N=N(d,s,q,q_0,\nu,R_0)>0$ such that
   \[
      \norm{D u}{\Leb{s,q}(Q_{1/2})}\leq N \left(\norm{u}{\Leb{s,1}(Q_{1})}+\norm{\boldF}{\Leb{s,q}(Q_{1})}+\norm{g}{\Leb{s,q}(Q_{1})} \right).
   \]
\end{theorem}
\begin{remark}\leavevmode
\begin{enumerate}[label=\textnormal{(\roman*)}]
\item {Unlike \cite{DP21}, we do not impose $u\in \tilde{W}^{1,2}_{s,q}(Q_1)$ or $u\in W^{0,1}_{s,q}(Q_1)$. Hence, our result is not an a priori estimate but a regularity estimate.}
\item If $u\in \tilde{W}^{1,2}_{q_0}(Q_1)^d$, then by the parabolic Sobolev embedding theorem, $u\in \Leb{s,1}(Q_1)$. However, if $u\in W^{0,1}_{q_0}(Q_1)^d$, then the norm $\norm{u}{\Leb{s,1}(Q_1)}$ is not always finite.
\item Due to Serrin's counterexample in \cite{S62},  weak and strong solutions may not possess good regularity in the time variable, i.e., it is not expected that $u_t \in \Leb{s,q}(Q_{1/2})$ for the case of equations in nondivergence form. Similarly, it is not expected that   $u_t \in \nSob{-1}{s,q}(Q_{1/2})$ for the case of equations in divergence form.
\item When $a^{ij}$ is merely measurable in $t$, then  {Theorems \ref{thm:B} and \ref{thm:D}} hold even for very weak solutions $u\in \Leb{s,1}(Q_1)^d$ (see Remark \ref{rem:t-dependent}). However, when $a^{ij}$ depends on $x$, it is unclear to us whether we could obtain interior mixed-norm derivative estimates for very weak solutions to \eqref{eq:Stokes-nondiv} and \eqref{eq:Stokes-div} since it is ambiguous to define the notion of very weak solutions.
\item In contrast to Theorems \ref{thm:A} and \ref{thm:C}, we do not need compatibility conditions on $g$.
\item In fact, {from its proof}{,} Theorem \ref{thm:B} still holds if we assume that $\partial_t u+\nabla p\in\Leb{1}(Q_1)^d$ instead of {assuming that} $\partial_t u, \nabla p \in \Leb{1}(Q_1)^d$. In this case, due to the lack of regularity in the time variable, it is not always guaranteed that $u\in \Leb{s,1}(Q_1)^d$.
\end{enumerate}
\end{remark}

One can also obtain a boundary version of Theorem \ref{thm:B} if we consider the Lions boundary conditions. We assume the following condition on viscosity coefficients:{
\begin{assumption}[$\delta$]\label{assump:VMO-boundary}
There exists $R_0 \in (0,1/4)$ such that for any $(t_0,x_0) \in\overline{Q_2^+}$ and $0<r<R_0$, there exists $\hat{a}^{ij}(t)$ satisfying uniform ellipticity \eqref{eq:elliptic} and
\[   \fint_{Q_r^+(t_0,x_0)} |a^{ij}(t,x)-\hat{a}^{ij}(t)|\myd{xdt}\leq \delta,\quad \text{for } i,j=1,\dots,d. \]
\end{assumption}}

\begin{theorem}\label{thm:boundary-estimate}
Let $1<s,q,q_0<\infty$. Then there exists $\delta>0$ such that under Assumption \ref{assump:VMO-boundary} $(\delta)$, if $(u,p)\in \tilde{W}^{1,2}_{q_0}(Q_1^+)^d \times \Sob{0,1}{1}(Q_1^+)$ is a strong solution to \eqref{eq:Stokes-nondiv} in $Q_1^+$ satisfying the Lions boundary conditions
\begin{equation}\label{eq:Lions-boundary}
D_d u^k= u^d=0\quad \text{on } (-1,0]\times B'_1 \times \{0\},\quad k=1,\dots,d-1
\end{equation}
for some $f \in \Leb{s,q}(Q_1^+)^d$ and $g \in \Sob{0,1}{s,q}(Q_1^+)$, then $D^2 u \in \Leb{s,q}(Q_{1/2}^+)$. Moreover, there exists a constant $N=N(d,s,q,q_0,\nu,R_0)>0$ such that
\[  \norm{D^2 u}{\Leb{s,q}(Q_{1/2}^+)}\leq N \left(\norm{u}{\Leb{s,1}(Q_1^+)}+\norm{f}{\Leb{s,q}(Q_1^+)}+\norm{Dg}{\Leb{s,q}(Q_1^+)} \right).\]
\end{theorem}
\begin{remark}\leavevmode
\begin{enumerate}
\item[(i)] {Unlike \cite{DKP22} which they assumed that $u\in\tilde{W}^{1,2}_{s,q}(Q_1^+)$, our result assumes $u\in\tilde{W}^{1,2}_{q_0}(Q_1^+)$ for some $q_0>1$. Hence, our result is a regularity result.}
\item[(ii)] Suppose that $(u,p) \in \tilde{W}^{1,2}_{q_0}(Q_1^+)^d\times \Sob{0,1}{1}(Q_1^+)$ is a strong solution to \eqref{eq:Stokes-nondiv} in $Q_1^+$ satisfying the Navier boundary conditions:
\[  \qquad D_d u^k-\alpha u^k=u^d=0\quad \text{on } (-1,0]\times B_1'\times \{0\},\quad k=1,\dots,d-1 \]
for some $\alpha>0$. If in addition $u,Du,p \in \Leb{s,q}(Q_1^+)$, then we can apply Theorem \ref{thm:boundary-estimate} to $(v,\pi)$ defined by $v(t,x)=e^{-\alpha x_d} u(t,x)$ and $\pi(t,x)=e^{-\alpha x_d} p(t,x)$ to get $D^2 u \in \Leb{s,q}(Q_{1/2}^+)$ and
\begin{align*}
    &\norm{D^2 u}{\Leb{s,q}(Q_{1/2}^+)}\\
    &\leq N \left(\norm{u}{\Sob{0,1}{s,q}(Q_1^+)}+\norm{p}{\Leb{s,q}(Q_1^+)}+\norm{Dg}{\Leb{s,q}(Q_1^+)}+\norm{f}{\Leb{s,q}(Q_1^+)} \right)
\end{align*}
for some constant $N=N(d,s,q,q_0,\nu,R_0,\alpha)>0$. 
\end{enumerate}
\end{remark}

\section{Preliminaries}\label{sec:prelim}

This section consists of four parts. In Subsection \ref{subsec:weight}, we list embedding theorems of function space $\dSob{1}{s,q}((0,T)\times\Omega)$, properties of $A_p$-weights, and Poincar\'e inequality on weighted spaces. In Subsection \ref{subsec:HL-FS}, we introduce Hardy-Littlewood maximal operator and Fefferman-Stein sharp maximal operator that will be used in this paper. Next, in Subsection \ref{subsec:div}, we state the solvability of the divergence equation in weighted Sobolev spaces. Finally, we state estimates of potentials on weighted spaces and list weighted solvability results for parabolic equations with simple coefficients in Subsection \ref{subsec:ep-weights}. These results will be used to construct a solution from vorticity in the remaining sections \ref{sec:4}, \ref{sec:5}, and \ref{sec:6}.

\subsection{Function spaces with and without weights}\label{subsec:weight}
In this subsection, we summarize several properties of function spaces with and without weights.

The following embedding result is a special case of Kim-Ryu-Woo \cite[Theorem 5.2]{KRW22}.
\begin{lemma}\label{lem:KRW}
Let $0<T<\infty$ and let $\Omega$ be a smooth bounded domain in $\mathbb{R}^d$, $d\geq 2$. Suppose that $1<s_0,q_0,s,q<\infty$ satisfy $s\leq s_0\leq \infty$, $q\leq q_0\leq \infty$ and either
\begin{enumerate}[label=\textnormal{(\roman*)}]
\item $s_0=s$ and $d/q\leq 1+d/q_0$,  $q\neq d$ or $q_0\neq \infty$; or
\item $s_0>s$ and $d/q+2/s\leq 1+d/q_0+2/s_0$.
\end{enumerate}
Then there exists a constant $N=N(d,s,q,T,\mathrm{diam } \,\Omega)>0$ such that
\[
     \norm{u}{\Leb{s_0,q_0}((0,T)\times\Omega)}\leq N \left(\norm{u}{\Leb{s,q}((0,T)\times\Omega)}+\norm{G}{\Leb{s,q}((0,T)\times\Omega)}\right)
\]
for all $u \in \Sob{0,1}{s,q}((0,T)\times\Omega)$ satisfying $u_t=\Div G$ for some $G \in \Leb{s,q}((0,T)\times\Omega)^d$.
\end{lemma}

Next we summarize some properties of $A_p$ weights and results on weighted Sobolev spaces, see e.g. Farwig-Sohr \cite[Lemmas 2.2 and 2.3]{FS97} and Grafakos \cite[Chapter 7]{G08}.

\begin{proposition}\label{prop:weight-property}
Let $1<p<\infty$ and $w\in A_p(\mathbb{R}^d,dx)$.
\begin{enumerate}[label=\textnormal{(\roman*)}]
\item $w^{-1/(p-1)} \in A_{p'}$ and $[w^{-1/(p-1)}]_{A_{p'}}=[w]_{A_p}^{1/(p-1)}$;
\item If $1<p<q<\infty$, then $w\in A_q$ and $[w]_{A_q}\leq [w]_{A_p}$;
\item There exists $1<q=q(d,p,[w]_{A_p})<p$ such that $w \in A_q$;
\item The functions defined by
\[
   |x|^\alpha \quad\text{and}\quad (1+|x|)^\alpha
\]
 are $A_p$-weights for all $-d<\alpha<d(p-1)$;
 \item There exist $\delta \in (0,1)$ and $N>0$ depending only on $d$, $p$, and $[w]_{A_p}$ such that
 \[
      \frac{w(S)}{w(B)}\leq N \left(\frac{|S|}{|B|}\right)^\delta
 \]
 for any ball $B$ {in $\mathbb{R}^d$} and any measurable subset $S$ of $B$;
 \item $w(B_R)\rightarrow \infty$ as $R\rightarrow\infty$.
\end{enumerate}
\end{proposition}
\begin{proof}
(i) This follows directly from the definition.

(ii) This follows directly from the definition and H\"older's inequality.

(iii) See e.g. \cite[Theorem 7.2.2]{G08}.

(iv) See Farwig-Sohr \cite[Lemmas 2.2 and 2.3]{FS97}.

(v) See e.g. \cite[Proposition 7.2.8]{G08}.

(vi)  For $R>1$, choose $S=B_1$ and $B=B_R$ in (vi). Then
\[
         \frac{w(B_1)}{w(B_R)} \leq N(d,p,[w]_{A_p})\left(\frac{|B_1|}{|B_R|}\right)^\delta.
\]
Since $|B_R|\rightarrow\infty$ as $R\rightarrow\infty$, it follows that $w(B_R)\rightarrow \infty$ as $R\rightarrow\infty$. This completes the proof of Proposition \ref{prop:weight-property}.
\end{proof}

The following weighted Poincar\'e inequality was first proved by Fabes-Kenig-Serapioni \cite{FKS82} and later simplified by Chiarenza-Frasca \cite{CF85}.

{
\begin{lemma}\label{lem:weighted-Poincare} 
Let $1<p<\infty$ and $w\in A_p$. There exists a $\delta=\delta(d,p)>0$ such that for $1\leq k\leq d/(d-1)+\delta$, there exists a constant $N=N(d,p,[w]_{A_p})>0$ such that
\[
   \left(\frac{1}{w(B_R)} \int_{B_R} |u|^{kp}w \myd{x}\right)^{1/(kp)} \leq NR \left(\frac{1}{w(B_R)}\int_{B_R} |\nabla u|^p w \myd{x} \right)^{1/p}
\]
for all $u\in C_0^\infty(B_R)$, $R>0$, and
\[
   \left(\frac{1}{w(B_R)} \int_{B_R} |u-(u)_{B_R,w}|^{kp}w \myd{x}\right)^{1/(kp)} \leq NR \left(\frac{1}{w(B_R)}\int_{B_R} |\nabla u|^p w \myd{x} \right)^{1/p}
\]
for all $u\in C^\infty(\overline{B_R})$, where $(u)_{B_R,w}=w(B_R)^{-1}\int_{B_R} u w \myd{x}$.
\end{lemma}
}

\subsection{Hardy-Littlewood maximal function and Fefferman-Stein theorem on weighted spaces}\label{subsec:HL-FS}

{For $T \in (-\infty,\infty]$ and a locally integrable function $f:\mathbb{R}^{d}_T\rightarrow \mathbb{R}$, we define its \emph{Hardy-Littlewood maximal function} by
\[   M_Tf(t,x):=\sup_{Q_r(s,y)\ni (t,x)} \fint_{Q_r(s,y)} |f(r,z)|\myd{r}dz,\quad (t,x)\in \overline{\mathbb{R}^{d}_T}.\]
If $T=\infty$, we write $M_T f:=M f$.}
Muckenhoupt \cite{M72} first proved the boundedness of the Hardy-Littlewood maximal operator on weighted spaces $\Leb{q,w}(\mathbb{R}^d)$, $1<q<\infty$, and $w\in A_q(\mathbb{R}^d,dx)$.  By applying a version of the Rubio de Francia extrapolation theorem (see e.g.  \cite[Theorem 2.5]{DK18-weight}), we can also prove the mixed-norm version of the theorem of Muckenhoupt.
\begin{lemma}\label{lem:HL}
Let {$T \in (-\infty,\infty]$}, $1<s,q<\infty$, $K_0\geq 1$, $w(t,x)=w_1(x)w_2(t)$, $[w_1]_{A_q(\mathbb{R}^d,dx)}\leq K_0$, $[w_2]_{A_s(\mathbb{R},dt)}\leq K_0$. Then there exists a constant $N=N(d,s,q,K_0)>0$ such that
\[
\norm{M_Tf}{\Leb{s,q,w}(\mathbb{R}^{d}_T)}\leq N \norm{f}{\Leb{s,q,w}(\mathbb{R}^{d}_T)}
\]
for all $f\in \Leb{s,q,w}(\mathbb{R}^{d}_T)$.
\end{lemma}

To introduce another type of maximal operator that we need, let
\[
    \mathbb{C}_n:=\{ Q^n_{\vec{i}} = Q^n_{(i_0,i_1,\dots,i_d)} : \vec{i}=(i_0,i_1,\dots,i_d)\in \mathbb{Z}^d\},
\]
where $n\in \mathbb{Z}$ and
\[
    Q^n_{\vec{i}} =\left[\frac{i_0}{2^{2n}},\frac{i_0+1}{2^{2n}}\right)\times \left[\frac{i_1}{2^n},\frac{i_1+1}{2^n}\right)\times \cdots\times  \left[\frac{i_d}{2^n},\frac{i_d+1}{2^n}\right).
\]
Then the collection $\mathbb{C}_n$ is a filtration of partitions (see e.g. \cite[Chapter 3]{K08} or \cite[Theorem 2.1]{DK18-weight}). Define the \emph{dyadic sharp function} of $g$ by
\[
g_{\mathrm{dy}}^{\#}(t,x):=\sup_{n<\infty} \fint_{Q_{\vec{i}}^n\ni (t,x)} |g(s,y)- g_{|_n} (t,x)|\,dyds,
\]
where
\[  g_{|_n}(t,x)=\fint_{Q_{\vec{i}}^n} g(s,y)\,dyds,\quad (t,x) \in Q_{\vec{i}}^n.\]

The following version of the Fefferman-Stein theorem was proved in Dong-Kim \cite[Corollary 2.7]{DK18-weight}.

\begin{lemma}\label{lem:FS}
Let {$T\in (-\infty,\infty]$,} $p,q,\tilde{p},\tilde{q}\in (1,\infty)$, $K_0\geq 1$, $w_1 \in A_{\tilde{p}}(\mathbb{R}^d,dx)$, $w_2\in A_{\tilde{q}}(\mathbb{R},dt)$ { whose seminorm are less than $K_0$}. Then there exists a constant $N=N(K,p,q,\tilde{p},\tilde{q},K_0)>0$ such that
\[
 \norm{f}{\Leb{p,q,w}(\mathbb{R}^{d}_T)}\leq N \norm{f_{\mathrm{dy}}^{\#}}{\Leb{p,q,w}(\mathbb{R}^{d}_T)}
\]
for all $f\in \Leb{p,q,w}(\mathbb{R}^{d}_T)$.
\end{lemma}

\subsection{The equation \texorpdfstring{$\Div u =g$}{}} \label{subsec:div}
Let $\Omega$ be a bounded Lipschitz domain in $\mathbb{R}^d$, $d\geq 2$. We consider the following Dirichlet problem for the divergence equation:
\begin{equation}\label{eq:div-equation}
\left\{
\begin{aligned}
\Div u &=g &&\quad \text{in } \Omega,\\
u&=0 &&\quad \text{on } \partial\Omega.
\end{aligned}
\right.
\end{equation}

$\Sob{1}{q}$-solvability of the problem \eqref{eq:div-equation} is a classical result due to Bogovski\u{i} \cite{B80} by introducing an integral representation of solutions to the problem \eqref{eq:div-equation} on a star-shaped domain (see Galdi \cite{G11}).  This result was extended by Huber \cite{H11} to weighted Sobolev spaces as below.

\begin{theorem}\label{thm:div-equation}
Let $1<q<\infty$, $K_0\geq 1$, and $w\in A_q$ satisfying $[w]_{A_q}\leq K_0$. Then there exists a bounded linear operator
\[
   \mathcal{B} : \Leb{q,w,\#}(\Omega)\rightarrow \Sobd{1}{q,w}(\Omega)^d \]
such that $\Div(\mathcal{B}f)=f$ for all $f \in\Leb{q,w,\#}(\Omega)$, where $\Leb{q,w,\#}(\Omega)$ is the collection of all $f\in \Leb{q,w}(\Omega)$ with $\int_\Omega f \myd{x}=0$. Moreover, we have $\mathcal{B}f \in C_0^\infty(\Omega)^d$ if $f\in C_0^\infty(\Omega)$ with $\int_\Omega f \myd{x}=0$ and
\[
   \norm{\mathcal{B}f}{\Sob{1}{q,w}(\Omega)}\leq N \norm{f}{\Leb{q,w}(\Omega)}
\]
for all $f\in\Leb{q,w,\#}(\Omega)$, where the constant $N$ depends only on $d$, $q$, $K_0$, and $\Omega$.
\end{theorem}
\begin{remark}
The operator $\mathcal{B}$ is the same Bogovski\u{\i} operator introduced in \cite{B80}. If $\Omega$ is bounded star-shaped with respect to an open ball $B_R$ with $\overline{B_R}\subset \Omega$, then there exists a constant $N=N(d,q,K_0,\mathrm{diam }\, \Omega /R)>0$ such that
\[
       \norm{D (\mathcal{B}g)}{\Leb{q,w}(\Omega)}\leq N\norm{g}{\Leb{q,w}(\Omega)}
\]
for all $g \in \Leb{q,w,\#}(\Omega)$.
\end{remark}

\subsection{Potential estimates and solvability of parabolic equations on weighted spaces}\label{subsec:ep-weights}

In this subsection, we give some potential estimates on weighted $\Leb{q}$-spaces and state the solvability of elliptic and parabolic equations in weighted Sobolev spaces. We also state weighted a priori $\Leb{q}$-estimates for Poisson equations that will be used in this paper.

Let $\Phi$ be the fundamental solution of the Laplacian defined by
\begin{equation*}%\label{eq:fundamental-sol}
    \Phi(x)=\begin{dcases}
   \frac{1}{d(2-d)\omega_d}\frac{1}{|x|^{d-2}}& \quad \text{if } d\geq 3,\\
    \frac{1}{2\pi}\ln |x| & \quad \text{if } d=2,
    \end{dcases}
\end{equation*}
where $\omega_d$ is the volume of the unit ball in $\mathbb{R}^d$.

The following lemma will be used to prove the existence of weak and strong solutions to Stokes equations with simple coefficients. The proof is almost identical to that of Lemma 4.1 in \cite{DKP22}. The key difference is to apply weighted $\Leb{q}$-boundedness of singular integral operators (see e.g. \cite[\S 4.2, Chapter V]{S93})   instead of the unweighted version when we prove \eqref{eq:potential-2} and \eqref{eq:potential-5}. We omit its proof.
\begin{lemma}\label{lem:potential}
Let $1<q<\infty$, $1<q_0<d$, and $K\geq 1$. Fix $w\in A_q(\mathbb{R}^d,dx)$ with $[w]_{A_q}\leq K_0$. For each $h\in \Leb{q_0}(\mathbb{R}^d)\cap \Leb{q,w}(\mathbb{R}^d)$, define
\[	v_k(x)=\int_{\mathbb{R}^d} D_k \Phi(x-y)h(y)\myd{y}\quad \text{in } \mathbb{R}^d,\quad 1\leq k\leq d. \]
\begin{enumerate}[label=\textnormal{(\roman*)}]
\item $v_k \in \Leb{q_0^*}(\mathbb{R}^d)$ and $D v_k \in \Leb{q,w}(\mathbb{R}^d)$ with the estimate
\begin{align}
\norm{v_k}{\Leb{q_0^*}(\mathbb{R}^d)}&\leq N_1(d,q_0) \norm{h}{\Leb{q_0}(\mathbb{R}^d)},\notag\\%\label{eq:potential-1} \\
\norm{D v_k}{\Leb{q,w}(\mathbb{R}^d)}&\leq N_2(d,q,K_0)\norm{h}{\Leb{q,w}(\mathbb{R}^d)},\label{eq:potential-2}
\end{align}
where $q_0^*=dq_0/(d-q_0)$. We also have
\begin{equation*}%\label{eq:potential-3}
   \sum_{k=1}^d D_k v_k = h\quad \text{in } \mathbb{R}^d.
\end{equation*}
\item If $D h \in \Leb{q,w}(\mathbb{R}^d)$ in addition, then $D^2 v_k \in \Leb{q,w}(\mathbb{R}^d)$ with
\begin{equation*}%\label{eq:potential-4}
    \Delta v_k =D_k h\quad \text{in } \mathbb{R}^d
\end{equation*}
and
\begin{equation}\label{eq:potential-5}
      \norm{D^2 v_k}{\Leb{q,w}(\mathbb{R}^d)}\leq N(d,q,K_0) \norm{D_k h}{\Leb{q,w}(\mathbb{R}^d)}
\end{equation}
holds.
\item If $D h\in \Leb{q,w}(\mathbb{R}^d)\cap \Leb{q_0}(\mathbb{R}^d)$ in addition, then
\begin{equation*}
    Dv_k(x)=\int_{\mathbb{R}^d} D_k \Phi(x-y) Dh(y)\,dy\quad \text{in } \mathbb{R}^d.
\end{equation*}
\end{enumerate}
\end{lemma}

We will use the following weighted $\Leb{s,q}$-results that can be found in \cite[Theorem 5.2]{DK18-weight}.
\begin{theorem}\label{thm:classical-Lp}
  Let $0<T<\infty$, $K_0 \geq 1$, $1<s,q<\infty$, $w(t,x)=w_1(x)w_2(t)$, $[w_1]_{A_q(\mathbb{R}^d,dx)}\leq K_0$, and $[w_2]_{A_s(\mathbb{R},dt)}\leq K_0$.
\begin{enumerate}[label=\textnormal{(\roman*)}]
\item For every $f\in \Leb{s,q,w}(\mathbb{R}^d_T)$, there exists a unique $u\in \oSob{1,2}{s,q,w}(\mathbb{R}^d_T)$ satisfying
\begin{equation*}%\label{eq:ndiv-simple}
	\partial_t u -a^{ij}(t)D_{ij} u =f\quad \text{in } \mathbb{R}^d_T.
\end{equation*}
Moreover, we have
\[	\norm{D^2 u}{\Leb{s,q,w}(\mathbb{R}^d_T)}\leq N_1 \norm{f}{\Leb{s,q,w}(\mathbb{R}^d_T)} \]
and
\[	\norm{u}{\Sob{1,2}{s,q,w}(\mathbb{R}^d_T)}\leq N_2 \norm{f}{\Leb{s,q,w}(\mathbb{R}^d_T)} \]
for some constants  $N_1$ depending only on $d,s,q,K_0,\nu$ and $N_2$ depending only on $d,s,q,K_0,\nu,T$.
\item For every $F=(F^1,\dots,F^d)\in \Leb{s,q,w}(\mathbb{R}^d_T)^d$, there exists a unique $u\in \odSob{1}{s,q,w}(\mathbb{R}^d_T)$ satisfying
\begin{equation*}%\label{eq:div-simple}
\partial_t u -D_i (a^{ij}(t)D_j u)=\Div F\quad \text{in } \mathbb{R}^d_T,
\end{equation*}
i.e., $u\in \odSob{1}{s,q,w}(\mathbb{R}^d_T)$ satisfies
\[
-\int_0^T\int_{\mathbb{R}^d} u \partial_t \phi -a^{ij}(t)D_ju D_i \phi \myd{x}dt=-\int_0^T\int_{\mathbb{R}^d} F \cdot \nabla \phi \myd{x}dt
\]
for all $\phi \in C_0^\infty([0,T)\times\mathbb{R}^d)$.
Moreover, we have
\[	\norm{D u}{\Leb{s,q,w}(\mathbb{R}^d_T)}\leq N_1 \norm{F}{\Leb{s,q,w}(\mathbb{R}^d_T)} \]
and
\[	\norm{u}{\dSob{1}{s,q,w}(\mathbb{R}^d_T)}\leq N_2 \norm{F}{\Leb{s,q,w}(\mathbb{R}^d_T)} \]
for some constants  $N_1$ depending only on $d,s,q,K_0,\nu$ and $N_2$ depending only on $d,s,q,K_0,\nu,T$.
\end{enumerate}
\end{theorem}

We will also use the following regularity results which can be easily proved by using Theorem \ref{thm:classical-Lp}.
\begin{corollary}\label{cor:a-priori-regularity}
Let $K_0 \geq 1$, $1<q<\infty$, and $w\in A_q (\mathbb{R}^d,dx)$ satisfying $[w]_{A_q}\leq K_0$.
\begin{enumerate}[label=\textnormal{(\roman*)}]
\item If $u\in \Sob{2}{q,w}(\mathbb{R}^d)$ satisfies
\[    -\Delta u = f\quad \text{in } \mathbb{R}^d \]
for some $f\in \Leb{q,w}(\mathbb{R}^d)$, then there exists a constant $N=N(d,q,K_0)>0$ such that
\[     \norm{D^2 u}{\Leb{q,w}(\mathbb{R}^d)}\leq N \norm{f}{\Leb{q,w}(\mathbb{R}^d)}. \]
\item If $u\in \Sob{1}{q,w}(\mathbb{R}^d)$ satisfies
\[    -\Delta u = \Div F\quad \text{in } \mathbb{R}^d \]
for some $F\in \Leb{q,w}(\mathbb{R}^d)^d$, then there exists a constant $N=N(d,q,K_0)>0$ such that
\[     \norm{Du}{\Leb{q,w}(\mathbb{R}^d)}\leq N \norm{F}{\Leb{q,w}(\mathbb{R}^d)}. \]
\end{enumerate}
\end{corollary}

\section{Stokes equations with simple coefficients}\label{sec:4}

In this section, we consider the Cauchy problem for Stokes equations with simple coefficients, that is, for $0<T<\infty$,
\begin{equation}\label{eq:Stokes-nondiv-simple}
\left\{
\begin{alignedat}{2}
\partial_t u-a^{ij}(t)D_{ij} u+\nabla p &=f&&\quad \text{in } (0,T)\times\mathbb{R}^d,\\
\Div u &=g&&\quad \text{in } (0,T)\times\mathbb{R}^d,\\
u&=0&&\quad \text{on } \{t=0\}\times \mathbb{R}^d,
\end{alignedat}
\right.
\end{equation}
where the viscosity coefficient $a^{ij}$ is merely measurable in $t$ and  satisfies uniform ellipticity condition \eqref{eq:elliptic}. We also consider the Cauchy problem for Stokes equations in divergence form:
\begin{equation}\label{eq:Stokes-div-simple}
\left\{
\begin{alignedat}{2}
\partial_t u-D_i(a^{ij}(t)D_{j} u)+\nabla p &=\Div \boldF&&\quad \text{in } (0,T)\times\mathbb{R}^d,\\
\Div u &=g&&\quad \text{in } (0,T)\times\mathbb{R}^d,\\
u&=0&&\quad \text{on } \{t=0\}\times \mathbb{R}^d.
\end{alignedat}
\right.
\end{equation}

We first state the $\Sob{1,2}{s,q,w}$-solvability for the problem \eqref{eq:Stokes-nondiv-simple} in $\mathbb{R}^d_T$. {The argument is nearly the same as in the proof of Theorem 1.4 in \cite{DKP22}.} We will explain the main difference in  \ref{app:B} for the sake of completeness.
\begin{theorem}\label{thm:mixed-simple}
Let $1<s,q<\infty$, $0<T<\infty$, and let $K_0\geq 1$ be constant, $w(t,x)=w_1(x)w_2(t)$, $[w_1]_{A_q(\mathbb{R}^d,dx)}\leq K_0$, and $[w_2]_{A_s(\mathbb{R},dt)}\leq K_0$. Then for every $f\in \Leb{s,q,w}(\mathbb{R}^d_T)^d$, $g\in\odSob{1}{s,q,w}(\mathbb{R}^d_T)$, and $g_t=\Div G$ for some vector field $G =(G^1,\dots,G^d)\in \Leb{s,q,w}(\mathbb{R}^d_T)^d$ in the sense that
\begin{equation*}
\int_{\mathbb{R}^d_T} g\varphi_t \myd{x}dt=\int_{\mathbb{R}^d_T} G \cdot \nabla \varphi \myd{x}dt
\end{equation*}
for all $\varphi \in C_0^\infty([0,T)\times\mathbb{R}^d)$,  there exists a unique strong solution $(u,p)$ to \eqref{eq:Stokes-nondiv-simple} satisfying
\[
u\in\oSob{1,2}{s,q,w}(\mathbb{R}^d_T)^d,\quad \nabla p \in\Leb{s,q,w}(\mathbb{R}^d_T)^d,\quad
{(p(t,\cdot))_{B_1}=0\,\,\text{for all } t\in (0,T)}\]
Moreover, we have
\begin{align}
\norm{D^2 u}{\Leb{s,q,w}(\mathbb{R}^d_T)}&\leq N_1 (\norm{f}{\Leb{s,q,w}(\mathbb{R}^d_T)}+\norm{Dg}{\Leb{s,q,w}(\mathbb{R}^d_T)}),\label{eq:simple-Hessian} \\
\norm{\nabla p}{\Leb{s,q,w}(\mathbb{R}^d_T)}&\leq N_1 (\norm{f}{\Leb{s,q,w}(\mathbb{R}^d_T)}+\norm{G}{\Leb{s,q,w}(\mathbb{R}^d_T)}+\norm{Dg}{\Leb{s,q,w}(\mathbb{R}^d_T)}),\notag\\%\label{eq:simple-gradient-pressure}\\
\norm{u_t}{\Leb{s,q,w}(\mathbb{R}^d_T)}&\leq N_1 (\norm{f}{\Leb{s,q,w}(\mathbb{R}^d_T)}+\norm{G}{\Leb{s,q,w}(\mathbb{R}^d_T)}),\notag%\label{eq:simple-u-t-estimate}
\end{align}
and
\begin{equation}\label{eq:simple-norm-estimate}
\begin{aligned}
&\relphantom{=}\norm{u}{\Sob{1,2}{s,q,w}(\mathbb{R}^d_T)}+\norm{\nabla p}{\Leb{s,q,w}(\mathbb{R}^d_T)}\\
&\leq N_2 (\norm{f}{\Leb{s,q,w}(\mathbb{R}^d_T)}+\norm{Dg}{\Leb{s,q,w}(\mathbb{R}^d_T)}+\norm{G}{\Leb{s,q,w}(\mathbb{R}^d_T)}),
\end{aligned}
\end{equation}
where $N_1=N_1(d,s,q,K_0,\nu)>0$ and $N_2=N_2(d,s,q,K_0,\nu,T)>0$.
\end{theorem}
{
\begin{remark}
One may ask whether we can extend this theorem to a more generalized Stokes system 
\[ \partial_t u^k-a^{ij}_{kl}(t)D_{ij}u^k +\partial_k p=f_k\quad\text{and}\quad\Div u =g,\quad k=1,\dots,d,\]
where $a^{ij}_{kl}$ satisfies the uniform Legendre-Hadamard condition. Our argument is limited to proving such results since it is difficult to find an equation for the vorticity. 
\end{remark}}

By using Theorem \ref{thm:mixed-simple}, we also have the existence and uniqueness of weak solutions to \eqref{eq:Stokes-div-simple} in $\mathbb{R}^d_T$ as follows, which can be deduced from Theorem \ref{thm:mixed-simple} and a duality argument based on Theorem \ref{thm:classical-Lp}. We give {its} proof in  \ref{app:C} for the sake of completeness.

\begin{theorem}\label{thm:mixed-simple-div}
Let $1<s,q<\infty$, $0<T<\infty$, and let $K_0\geq 1$ be constant, $w(t,x)=w_1(x)w_2(t)$, $[w_1]_{A_q(\mathbb{R}^d,dx)}\leq K_0$, and $[w_2]_{A_s(\mathbb{R},dt)}\leq K_0$. For every $\boldF\in \Leb{s,q,w}(\mathbb{R}^d_T)^{d\times d}$ and $g\in \Leb{s,q,w}(\mathbb{R}^d_T)$ satisfying $g_t=\Div\Div \boldG$ for some matrix field $\mathbf{G} \in \Leb{s,q,w}(\mathbb{R}^d_T)^{d\times d}$ in the sense that
\begin{equation}\label{eq:compatibility-divergence}
 \int_{\mathbb{R}^d_T} g \varphi_t \myd{x}dt=-\int_{\mathbb{R}^d_T} \boldG:\nabla^2 \varphi \myd{x}dt
\end{equation}
for all $\varphi \in C_0^\infty([0,T)\times\mathbb{R}^d)$, there exists a unique weak solution $(u,p)$ of \eqref{eq:Stokes-div-simple} satisfying
\[
u\in\odSob{1}{s,q,w}(\mathbb{R}^d_T)^d,
\quad  p \in\Leb{s,q,w}(\mathbb{R}^d_T).
\]
Moreover, we have
\begin{align*}%\label{eq:a-priori-estimate-Du}
\norm{Du}{\Leb{s,q,w}(\mathbb{R}^d_T)}&\leq N_1 \left(\norm{g}{\Leb{s,q,w}(\mathbb{R}^d_T)}+\norm{\boldF}{\Leb{s,q,w}(\mathbb{R}^d_T)} \right),\\
%\label{eq:a-priori-estimate-Du-2}
\norm{p}{\Leb{s,q,w}(\mathbb{R}^d_T)}&\leq N_1 \left(\norm{g}{\Leb{s,q,w}(\mathbb{R}^d_T)}+\norm{\boldF}{\Leb{s,q,w}(\mathbb{R}^d_T)}+\norm{\boldG}{\Leb{s,q,w}(\mathbb{R}^d_T)} \right),
\end{align*}
and
\[
  \norm{u}{\dSob{1}{s,q,w}(\mathbb{R}^d_T)}+\norm{p}{\Leb{s,q,w}(\mathbb{R}^d_T)}\leq N_2 \left(\norm{\boldF}{\Leb{s,q,w}(\mathbb{R}^d_T)}+\norm{g}{\Leb{s,q,w}(\mathbb{R}^d_T)}+\norm{\boldG}{\Leb{s,q,w}(\mathbb{R}^d_T)}\right),
\]
where $N_1=N_1(d,s,q,K_0,\nu)>0$ and  $N_2=N_2(d,s,q,K_0,\nu,T)>0$.
\end{theorem}

Recall that for $0<\alpha\leq 1$ and each parabolic cylinder $Q$ in $\mathbb{R}^{d+1}$, we write
\[
   [u]_{C^{\alpha/2,\alpha}(Q)} :=\sup_{(t,x),(s,y) \in Q, (t,x)\neq (s,y)} \frac{|u(t,x)-u(s,y)|}{|t-s|^{\alpha/2}+|x-y|^\alpha}
\]
and we define
\[
  \norm{u}{C^{\alpha/2,\alpha}(Q)}:=\norm{u}{\Leb{\infty}(Q)}+[u]_{C^{\alpha/2,\alpha}(Q)}.
\]

We have the following local H\"older estimates for the vorticity of a solution to nonstationary Stokes equations.
\begin{lemma}\label{lem:Holder-estimate-Domega}
Let $1<q_0<\infty$ and let $u \in \Sob{0,1}{q_0}(Q_1)^d$ be a weak solution of
\begin{equation}\label{eq:Stokes-Holder-2}
\partial_t u -D_i(a^{ij}(t)D_{j} u)+\nabla p =0,\quad \Div u =g(t)
\end{equation}
in $Q_1$ {for some measurable function $g:(-1,0)\rightarrow\mathbb{R}$.} There exists a constant $N=N(d,\nu,q_0)>0$ such that
\[	\norm{\omega}{C^{1/2,1}(Q_{1/2})} \leq N \norm{\omega}{\Leb{q_0}(Q_1)}, \]
where $\omega=\nabla \times u$.
\end{lemma}
\begin{proof}
By taking mollification in $x$, we may assume that $u$ is smooth in $x$. For $\psi \in C_0^\infty(Q_1)$ and $k,l=1,\dots,d$, define {$\phi=(D_k \psi)e_l-(D_l \psi)e_k$,} where $\{e_k\}$ is the standard basis for $\mathbb{R}^d$. Then it is easy to see that $\Div \phi(t,\cdot)=0$ in $B_1$ for $t\in (-1,0)$. If we use $\phi$ as a test function in the definition of the weak solution, then one can show that {$\omega=\nabla \times u$} is a very weak solution to the vorticity equation
\[
\partial_t \omega^{kl}-D_i(a^{ij}(t)D_{j} \omega^{kl})=0\quad \text{in } Q_1.\]
Then the desired result follows from interior estimates for parabolic equations with coefficients measurable in $t$, Sobolev embedding theorem, and a standard iteration argument. We omit the details (see e.g. \cite[Chapter 2]{K08} and \cite{DK11}).
\end{proof}

Since $a^{ij}$ depends only on $t$,  by using Lemma \ref{lem:Holder-estimate-Domega} {and} a standard scaling argument, we have the following mean oscillation estimate for vorticity and its gradient. See e.g. \cite[Lemma 4]{DK11} for the proof.

\begin{lemma}\label{lem:mean-oscillation}
Let $\kappa \geq 8$, $0<r<\infty$, and $1<q_0<\infty$.
\begin{enumerate}[label=\textnormal{(\roman*)}]
\item  If $(u,p)$ is a weak solution of \eqref{eq:Stokes-Holder-2}  in $Q_{\kappa r}(X_0)$ and $\omega = \nabla\times u$, then there exists a constant $N=N(d,q_0,\nu)>0$ such that
\begin{equation*}%\label{eq:mean-oscillation-omega-2}
	(|\omega-(\omega)_{Q_r(X_0)}|)_{Q_r(X_0)} \leq N \kappa^{-1} (|\omega|^{q_0})_{Q_{\kappa r}(X_0)}^{1/q_0}.
\end{equation*}
\item If $(u,p)$ is a strong solution of \eqref{eq:Stokes-Holder-2}  in $Q_{\kappa r}(X_0)$ and $\omega = \nabla\times u$, then there exists a constant $N=N(d,q_0,\nu)>0$ such that
\begin{equation*}%\label{eq:mean-oscillation-omega}
	(|D\omega-(D\omega)_{Q_r(X_0)}|)_{Q_r(X_0)} \leq N \kappa^{-1} (|D\omega|^{q_0})_{Q_{\kappa r}(X_0)}^{1/q_0}.
\end{equation*}
\end{enumerate}
\end{lemma}

\section{Stokes equations in nondivergence form}\label{sec:5}
This section is devoted to proving Theorem \ref{thm:A}. We first obtain a mean oscillation estimate for the gradient of vorticity $\omega=\nabla \times u$ of a strong solution $u$ to \eqref{eq:Stokes-nondiv}.

\begin{lemma}\label{lem:mixed-norm}
Let $\kappa \geq 8$, $\delta \in (0,1)$, $q,\mu,\mu'\in(1,\infty)$, $1/\mu+1/\mu'=1$, and $a^{ij}$ satisfy Assumption \ref{assump:VMO} $(\delta)$. Then for any $0< r\leq R_0/\kappa $, $(t_0,x_0) \in \mathbb{R}^{d+1}$, and $(u,p)\in \Sobloc{1,2}{q\mu}(\mathbb{R}^{d+1})^d\times \Sobloc{0,1}{1}(\mathbb{R}^{d+1})$ satisfying
\begin{equation}\label{eq:Stokes-div-mean}
  	\partial_t u-a^{ij}(t,x)D_{ij} u+\nabla p=f,\quad \Div u=g\quad \text{in } Q_{\kappa r}(t_0,x_0),
\end{equation}
where $f\in \Lebloc{q}(\mathbb{R}^{d+1})^d$ and $g\in \Sobloc{0,1}{q}(\mathbb{R}^{d+1})$, we have
\begin{align*}
&\left(|D\omega-(D\omega)_{Q_r(t_0,x_0)}|\right)_{Q_r(t_0,x_0)}\\
&\leq N \kappa^{-1} \left[(|D^2 u|^q)^{1/q}_{Q_{\kappa r}(t_0,x_0)} +(|f|^{q})^{1/q}_{Q_{\kappa r}(t_0,x_0)}+(|Dg|^{q})^{1/q}_{Q_{\kappa r}(t_0,x_0)}\right]\\
&\relphantom{=}+N\kappa^{(d+2)/q} \left[ (|f|^{q})_{Q_{\kappa r}(t_0,x_0)}^{{1}/{q}} + \delta^{1/(q\mu')}(|D^2 u|^{q\mu})_{Q_{\kappa r}(t_0,x_0)}^{{1}/{(q\mu)}} +\left(|Dg|^q \right)^{1/q}_{Q_{\kappa r}(t_0,x_0)}\right]
\end{align*}
 for some constant $N=N(d,q,\nu)>0$.
\end{lemma}

\begin{proof}
For an integrable function $h$ defined on $Q_r$, define
\begin{equation}\label{eq:h-mollification}
   h^{(\varepsilon)}(t,x)=\int_{Q_\varepsilon} h(t+s,x+y)\eta_\varepsilon(s,y)\myd{y}ds,\quad (t,x) \in (-r^2+\varepsilon^2,0)\times B_{r-\varepsilon},
\end{equation}
where $\eta \in C_0^\infty(Q_1)$, $\eta_\varepsilon(t,x)=\varepsilon^{-d-2} \eta(t/\varepsilon^2,x/\varepsilon)$, and $\int_{Q_1} \eta \myd{x}dt=1$.

By mollifying equation \eqref{eq:Stokes-div-mean}, we get
\[
     \partial_t u^{(\varepsilon)}-a^{ij}D_{ij}u^{(\varepsilon)}+\nabla p^{(\varepsilon)}=f^{(\varepsilon)}+(a^{ij}D_{ij} u)^{(\varepsilon)}-a^{ij}D_{ij}u^{(\varepsilon)}
\]
in $Q_{r'}(t_0,x_0)$ for $0<r'<\kappa r$ {and for sufficiently small $\varepsilon$}. If we prove the estimate in the lemma for $u^{(\varepsilon)}$, we get the desired result by letting $\varepsilon\rightarrow 0$. Hence we may assume that $u$ and $p$ are infinitely differentiable. Since $(u,p)$ satisfies \eqref{eq:Stokes-div-mean}, it follows that $g\in \dSobloc{1}{q}(\mathbb{R}^{d+1})$.

By translation invariance, we may assume that $(t_0,x_0)=(0,0)$. Let  $\zeta_r(x)$ and $\psi_r(t)$ be infinitely differentiable functions defined on $\mathbb{R}^d$ and $\mathbb{R}$ satisfying
\[
  \zeta_r(x)=1\quad\text{on } B_{2r/3},\quad \zeta_r(x)=0\quad \text{on } \mathbb{R}^d\setminus B_r,
\]
\[
  \psi_r(t)=1\quad\text{on } t\in (-4r^2/9,4r^2/9),\quad \psi_r(t)=0\quad \text{on } \mathbb{R}\setminus (-r^2,r^2).
\]
Set $\phi_r(t,x)=\psi_r(t)\zeta_r(x)$. Then $\phi_r=1$ on $Q_{2r/3}$ and $|D\phi_r|\leq 4/r$.

Consider the following initial-value problem for Stokes equations:
\begin{equation}\label{eq:Stokes-w-part}
   \left\{\begin{alignedat}{2}
\partial_t u_1-\overline{a}^{ij}(t)D_{ij} u_1 + \nabla p_1&=1_{Q_{\kappa r}} h&&\quad \text{in } (-(\kappa r)^2,0)\times \mathbb{R}^d,\\
\Div u_1 &=\tilde{g} &&\quad \text{in } (-(\kappa r)^2,0)\times \mathbb{R}^d,\\
u_1&=0&&\quad \text{on } \{t=-\kappa r^2\}\times \mathbb{R}^d,
\end{alignedat}
\right.
\end{equation}
where
\begin{gather*}
h(t,x)=[f+ (a^{ij}-\overline{a}^{ij}(t))D_{ij} u],\quad \tilde{g}(t,x):=(g-[g(t,\cdot)]_{{\zeta_{\kappa r,B_{\kappa r}}}})\phi_{\kappa r},\\
  [f]_{\zeta_{r,B_r}}:= \frac{1}{\left(\int_{B_r} \zeta_r dx\right)}\int_{B_r} f \zeta_r \,dx.
\end{gather*}

By using the Poincar\'e inequality, it is easily seen that
\begin{equation}\label{eq:Poincare-g-tilde}
    \norm{D\tilde{g}}{\Leb{q}((-(\kappa r)^2,0)\times \mathbb{R}^d)}\leq N(d,q)\norm{Dg}{\Leb{q}(Q_{\kappa r})}.
\end{equation}
Also, since $g\in \dSobloc{1}{q}(\mathbb{R}^{d+1})$, it follows that
\[
\tilde{g} \in \odSob{1}{q}((-(\kappa r)^2,0)\times \mathbb{R}^d).
\]
It remains to show the compatibility condition, i.e., there exists $\tilde{G} \in \Leb{q}((-(\kappa r)^2,0)\times \mathbb{R}^d)^d$ such that
\[
  \partial_t \tilde{g}=\Div \tilde{G}\quad \text{in } (-(\kappa r)^2,0)\times \mathbb{R}^d
\]
in the sense of  \eqref{eq:compatibility-nondiv}.

Note that
\[
  \partial_t \tilde{g}=(\partial_t g-[\partial_t g(t,\cdot)]_{{\zeta_{\kappa r,B_{\kappa r}}}})\phi_{\kappa r}+( g-[g(t,\cdot)]_{{\zeta_{\kappa r,B_{\kappa r}}}})\partial_t\phi_{\kappa r}.
\]
By integrating it over $B_{\kappa r}$, we have
\[
  \int_{B_{\kappa r}} \partial_t \tilde{g} \myd{x} =0.
\]
Hence by Theorem \ref{thm:div-equation}, there exists $G \in \Leb{q}(-(\kappa r)^2,0;\Sobd{1}{q}(B_{\kappa r})^d)$ such that
\[   \Div G = \partial_t\tilde{g}\quad \text{in }(-(\kappa r)^2,0)\times B_{\kappa r},\quad G=0\quad \text{on } (-(\kappa r)^2,0) \times \partial B_{\kappa r}.\]
Extend $G$ to be zero outside $(-(\kappa r)^2,0)\times B_{\kappa r}$ and denote this extension by $\tilde{G}$. Since $\tilde{g}$ has compact support on $(-(\kappa r)^2,0)\times B_{\kappa r}$ and $G(t,\cdot)=0$ on $(-(\kappa r)^2,0) \times \partial B_{\kappa r}$ {for $t\in (-(\kappa r)^2,0)$}, we see that
\[
  \Div \tilde{G}=\partial_t\tilde{g}\quad \text{in } (-(\kappa r)^2,0)\times\mathbb{R}^d
\]
in the sense of \eqref{eq:compatibility-nondiv}.
Since $u\in \Sobloc{1,2}{q}(\mathbb{R}^{d+1})^d$ satisfies \eqref{eq:Stokes-div-mean}, it follows from Theorem \ref{thm:mixed-simple} that there exists a unique strong solution $(u_1,p_1)$ to \eqref{eq:Stokes-w-part} satisfying
\[
u_1\in \oSob{1,2}{q}((-(\kappa r)^2,0)\times\mathbb{R}^d)^d,\quad \nabla p_1\in \Leb{q}((-(\kappa r)^2,0)\times\mathbb{R}^d)^d.
\]
 Moreover, it follows from \eqref{eq:simple-Hessian} and \eqref{eq:Poincare-g-tilde} that
 \begin{equation}\label{eq:Lq-estimate-perturbed}
  \begin{aligned}
	\norm{D^2 u_1}{\Leb{q}((-(\kappa r)^2,0)\times\mathbb{R}^d)}&\leq N \norm{1_{Q_{\kappa r}} [f+(a^{ij}-\overline{a}^{ij})D_{ij}u]}{\Leb{q}((-(\kappa r)^2,0)\times\mathbb{R}^d)}\\
  &\relphantom{=}+N\norm{Dg}{\Leb{q}(Q_{\kappa r})},
  \end{aligned}
\end{equation}
where $N=N(d,q,\nu)>0$.

Define $(u_2,p_2)=(u-u_1,p-p_1)$. Then $(u_2,p_2)$ satisfies
\[ \left\{\begin{alignedat}{1}
\partial_t u_2-\overline{a}^{ij}(t)D_{ij} u_2 + \nabla p_2&=0\\
\Div u_2 &={[g(t,\cdot)]_{\zeta_{\kappa r},B_{\kappa r}}}
\end{alignedat}
\right.\,\,\,\text{in } Q_{2\kappa r/3}.
\]
Write $\omega=\nabla \times u$, $\omega_1=\nabla \times u_1$, and $\omega_2=\nabla \times u_2$. By Lemma \ref{lem:mean-oscillation} (i), we have
\begin{equation}\label{eq:omega1-oscillation}
\begin{aligned}
\left( |D\omega_2-(D\omega_2)_{Q_r}| \right)_{Q_r}&\leq N\kappa^{-1} (|D\omega_2|^{q})^{1/q}_{Q_{2\kappa r/3}}\\
&\leq N(d,q,\nu)\kappa^{-1}\left[ (|D \omega|^{q})_{Q_{2\kappa r/3}}^{1/q}+(|D \omega_1|^{q})_{Q_{2\kappa r/3}}^{1/q}\right].
\end{aligned}
\end{equation}
Since $a^{ij}$ satisfies Assumption \ref{assump:VMO} $(\delta)$ and $a^{ij}$, $\overline{a}^{ij}$ are bounded by $\nu^{-1}$,  it follows from \eqref{eq:Lq-estimate-perturbed} and H\"older's inequality that
\begin{equation}\label{eq:omega2-oscillation}
(|D\omega_1|^{q})_{Q_{\kappa r}}^{1/q} \leq N(d,q,\nu)\left[(|f|^{q})_{Q_{\kappa r}}^{1/q} +  \delta^{1/(q\mu')} \left( |D^2 u|^{q\mu} \right)^{1/(q\mu)}_{Q_{\kappa r}}+\left(|Dg|^q \right)^{1/q}_{Q_{\kappa r}}\right].
\end{equation}

By \eqref{eq:omega1-oscillation} and {\eqref{eq:omega2-oscillation}},  we get
\begin{align*}
	&(|D\omega_2-(D\omega_2)_{Q_r}|)_{Q_r} \\
&\leq N(d,q,\nu) \kappa^{-1}\left[ (|D\omega|^{q})^{1/q}_{Q_{\kappa r}} +  (|f|^{q})^{1/q}_{Q_{\kappa r}} +(|Dg|^{q})^{1/q}_{Q_{\kappa r}}+  (|D^2 u|^q)_{Q_{\kappa r}}^{1/q}\right].
\end{align*}
Since $\omega=\omega_1+\omega_2$,  {the above inequality} and \eqref{eq:omega2-oscillation} imply
\begin{align*}
&\fint_{Q_r} |D\omega-(D\omega)_{Q_r}|\myd{xdt}\\
&\leq \fint_{Q_r} |D\omega_2-(D\omega_2)_{Q_r}|\myd{xdt}+2\fint_{Q_r} |D\omega_1 |\myd{xdt}\\
&\leq N \kappa^{-1} \left[(|D^2 u|^q)^{{1}/{q}}_{Q_{\kappa r}} +(|f|^{q})^{{1}/{q}}_{Q_{\kappa r}}+(|Dg|^{q})_{Q_{\kappa r}}^{{1}/{q}}\right]\\
&\relphantom{=}+N\kappa^{{(d+2)}/{q}} \left[ (|f|^{q})_{Q_{\kappa r}}^{{1}/{q}} + \delta^{{1}/{(q\mu')}}(|D^2 u|^{q\mu})_{Q_{\kappa r}}^{{1}/{(q\mu)}}+(|Dg|^q)^{1/q}_{Q_{\kappa r}} \right]
\end{align*}
for some constant $N=N(d,q,\nu)>0$. This completes the proof of Lemma \ref{lem:mixed-norm}.
\end{proof}

The following proposition does not require the compatibility condition on $g$ since it only involves a priori estimates for $D^2u$.

\begin{proposition}\label{prop:Stokes-partition}
Let $0<T<\infty$, $K_0\geq 1$, $1<s,q<\infty$, $t_1 \in \mathbb{R}$, $w(t,x)=w_1(x)w_2(t)$, $[w_1]_{A_q(\mathbb{R}^d,dx)}\leq K_0$, and $[w_2]_{A_s(\mathbb{R},dt)}\leq K_0$. Then there exist $\delta>0$ and $R_1>0$ such that under Assumption \ref{assump:VMO} $(\delta)$, for any $u\in \oSob{1,2}{s,q,w}(\mathbb{R}^{d}_T)^d$ vanishing outside $(t_1-(R_0R_1)^2,t_1)\times \mathbb{R}^d$ and $p\in \Sobloc{0,1}{1}(\mathbb{R}^{d}_T)$ satisfying
\[	\partial_t u-a^{ij}(t,x)D_{ij} u+\nabla p=f,\quad \Div u=g\quad \text{in } \mathbb{R}^{d}_T,\]
 where $f\in \Leb{s,q,w}(\mathbb{R}^{d}_T)^d$ and $g\in \Sob{0,1}{s,q,w}(\mathbb{R}^{d}_T)$,  there exists a constant $N=N(d,s,q,K_0,\nu)>0$ such that
\[	\norm{D^2 u}{\Leb{s,q,w}(\mathbb{R}^{d}_T)}\leq N \left(\norm{f}{\Leb{s,q,w}(\mathbb{R}^{d}_T)}+\norm{Dg}{\Leb{s,q,w}(\mathbb{R}^{d}_T)}\right).\]
\end{proposition}
\begin{proof}
Since $w_1 \in A_q(\mathbb{R}^{d},dx)$ and $w_2 \in A_s(\mathbb{R},dt)$, it follows from Proposition \ref{prop:weight-property} (iii) that there exist {
$\sigma_1>0$ and $\sigma_2>0$} such that $q-\sigma_1>1$, $s-\sigma_2>1$, and
\[	w_1 \in A_{q-\sigma_1}(\mathbb{R}^{d},dx),\quad w_2\in A_{s-\sigma_2}(\mathbb{R},dt).\]

Choose $q_0,\mu \in (1,\infty)$ so that
\[	q_0\mu = \min\left\{\frac{q}{q-\sigma_1},\frac{s}{s-\sigma_2} \right\}>1.\]
By Proposition \ref{prop:weight-property} (ii), we also have
\begin{equation*}
\begin{aligned}
w_1 \in A_{q-\sigma_1}\subset A_{q/(q_0\mu)}&\subset A_{q/q_0}(\mathbb{R}^{d},dx),\\
w_2 \in A_{s-\sigma_2}\subset A_{q/(q_0\mu)}&\subset A_{s/q_0}(\mathbb{R},dt).
\end{aligned}
\end{equation*}
Then by H\"older's inequality (see e.g. \cite[Lemma 5.10]{DK18-weight}), we have
\[	u \in \Sobloc{1,2}{q_0\mu}(\mathbb{R}^{d}_T)^d,\quad f\in \Lebloc{q_0\mu}(\mathbb{R}^{d}_T)^d,\quad \text{and}\quad g \in \Sobloc{0,1}{q_0\mu}(\mathbb{R}^{d}_T).\]

Let $\kappa \geq 8$, $0<\delta<1$, and $R_1>0$ be constants to be specified below. For each $(t,x)\in \mathbb{R}^{d}_T$ and $Q^n \in\mathbb{C}_n$ such that $(t,x)\in Q^n$, $n\in \mathbb{Z}$, find $(t_0,x_0) \in \mathbb{R}^{d}_T$ and the smallest $r\in (0,\infty)$ so that $Q^n \subset Q_r(t_0,x_0)$ and
\begin{equation*}%\label{eq:dyadic-sharp}
	\fint_{Q^n} |f(s,y)-f|_n(t,x)|\myd{sdy} \leq N \fint_{Q_r(t_0,x_0)} |f(s,y)-(f)_{Q_r(t_0,x_0)}|\myd{sdy},
\end{equation*}
where $N$ depends only on $d$.

On one hand, if $r>R_0/\kappa$, since $u$ vanishes outside $(t_1-(R_0R_1)^2,t_1)\times \mathbb{R}^d$, we have
\begin{align*}
&\fint_{Q^n} |D\omega(s,y)-(D\omega)|_n (t,x)|\myd{sdy}\\
&\leq N \fint_{Q_r(t_0,x_0)}|D\omega(s,y)-(D\omega)_{Q_r(t_0,x_0)}| \myd{sdy} \\
&\leq N\left(\fint_{Q_r(t_0,x_0)} I_{(t_1-(R_0R_1)^2,t_1)}\myd{sdy} \right)^{1-{1}/{q_0}}\left(\fint_{Q_r(t_0,x_0)} |D\omega|^{q_0}\myd{sdy}\right)^{{1}/{q_0}}\\
&\leq N\kappa^{2\left(1-{1}/{q_0}\right)} R_1^{2\left(1-{1}/{q_0}\right)} [{M_T}(|D\omega|^{q_0})]^{{1}/{q_0}}(t,x)
\end{align*}
for some constant $N=N(d)>0$.

On the other hand, if $0<r\leq R_0/\kappa$, then it follows from Lemma \ref{lem:mixed-norm} that
\begin{align*}
&\fint_{Q^n} |D\omega(s,y)-(D\omega)|_n (t,x)|\myd{sdy}\\
&\leq N \fint_{Q_r(t_0,x_0)} |D\omega-(D\omega)_{Q_r(t_0,x_0)} | \myd{sdy}\\
&\leq N \kappa^{-1}\left[ (|D^2 u|^{q_0})^{1/q_0}_{Q_{\kappa r}(t_0,x_0)} +(|f|^{q_0})^{1/q_0}_{Q_{\kappa r}(t_0,x_0)}+(|Dg|^{q_0})^{1/q_0}_{Q_{\kappa r}(t_0,x_0)}\right]\\
&\relphantom{=}+N\kappa^{{(d+2)}/{q_0}} \left[ (|f|^{q_0})_{Q_{\kappa r}(t_0,x_0)}^{1/q_0} + \delta^{{1}/{(q_0\mu')}}(|D^2 u|^{q_0\mu})_{Q_{\kappa r}(t_0,x_0)}^{{1}/{(q_0\mu)}}+(|Dg|^{q_0})^{{1}/{q_0}}_{Q_{\kappa r}(t_0,x_0)} \right]\\
&\leq N\kappa^{-1} {M_T}(|D^2 u|^{q_0})^{{1}/{q_0}}(t,x)+ N(\kappa^{-1}+\kappa^{{(d+2)}/{q_0}}) {M_T}(|f|^{q_0})^{1/q_0}(t,x)\\
&\relphantom{=} + N\kappa^{{(d+2)}/{q_0}} \delta^{{1}/{(q_0\mu')}} {M_T}(|D^2 u|^{q_0\mu})^{{1}/{(q_0\mu)}}(t,x)\\
&\relphantom{=}+N(\kappa^{-1}+\kappa^{{(d+2)}/{q_0}})  {M_T}(|Dg|^{q_0})^{{1}/{q_0}}(t,x)
\end{align*}
for {$(t,x) \in \mathbb{R}^d_T$} and  some constant $N=N(d,q_0,\nu)>0$. Hence by taking the supremum with respect to all $Q^n \ni (t,x)$, $n\in \mathbb{Z}$, we see that
\begin{equation}\label{eq:sharp-function-estimate}
\begin{aligned}
&(D \omega)^{\#}_{dy}(t,x)\\
&\leq N\left(\kappa^{-1}+\kappa^{(d+2)/q_0}\delta^{{1}/{(q_0\mu')}}+\kappa^{2(1-1/q_0)} R_1^{2(1-1/q_0)} \right){M_T}(|D^2 u|^{q_0\mu})^{{1}/{(q_0\mu)}}(t,x)\\
&\relphantom{=}+N(\kappa^{-1}+\kappa^{(d+2)/q_0}) {M_T}(|f|^{q_0})^{1/q_0}(t,x)\\
&\relphantom{=}+N(\kappa^{-1}+\kappa^{{(d+2)}/{q_0}}) {M_T}(|Dg|^{q_0})^{{1}/{q_0}}(t,x)
\end{aligned}
\end{equation}
for some constant $N=N(d,q_0,\nu)>0$ and for all $(t,x)\in \mathbb{R}^{d}_T$.

By Lemma \ref{lem:HL}, we have
\begin{equation}\label{eq:HL-control}
\begin{aligned}
  \norm{[M_T(|D^2 u|^{q_0\mu})]^{1/(q_0\mu)}}{\Leb{s,q,w}(\mathbb{R}^{d}_T)}&=\norm{M_T(|D^2 u|^{q_0\mu})}{\Leb{s/(q_0\mu),q/(q_0\mu),w}(\mathbb{R}^{d}_T)}^{1/(q_0\mu)}\\
  &\leq N \norm{|D^2 u|^{q_0\mu}}{\Leb{s/(q_0\mu),q/(q_0\mu),w}(\mathbb{R}^{d}_T)}^{1/(q_0\mu)}\\
  &=N\norm{D^2 u}{\Leb{s,q,w}(\mathbb{R}^{d}_T)}
\end{aligned}
\end{equation}
for some constant $N=N(d,s/(q_0\mu),q/(q_0\mu),K_0)>0$ and hence $N=N(d,s,q,K_0)>0$.
Hence it follows from \eqref{eq:sharp-function-estimate}, the generalized Fefferman-Stein theorem (Lemma \ref{lem:FS}), and \eqref{eq:HL-control} that
\begin{equation}\label{eq:Domega-mixed-weight-estimate}
\begin{aligned}
&\norm{D\omega}{\Leb{s,q,w}(\mathbb{R}^{d}_T)}\\
&\leq N\left(\kappa^{-1}+\kappa^{(d+2)/q_0}\delta^{{1}/{(q_0\mu)}}+\kappa^{2(1-1/q_0)} R_1^{2(1-1/q_0)} \right)\norm{D^2 u}{\Leb{s,q,w}(\mathbb{R}^{d}_T)}\\
&\relphantom{=}+N(\kappa^{-1}+\kappa^{(d+2)/q_0})(\norm{f}{\Leb{s,q,w}(\mathbb{R}^{d}_T)}+\norm{Dg}{\Leb{s,q,w}(\mathbb{R}^{d}_T)})
\end{aligned}
\end{equation}
for some constant $N=N(d,s,q,K_0,\nu)>0$.

Since $\Div u=g$ in $\mathbb{R}^d_T$, we get
\[	\Delta u^i=D_i g+\sum_{k\neq i} D_k(D_ku^i-D_i u^k)\quad  \text{in } \mathbb{R}^d_T,\quad 1\leq i\leq d.\]
Hence it follows from Corollary \ref{cor:a-priori-regularity} (i) and \eqref{eq:Domega-mixed-weight-estimate} that
\begin{align*}
&\norm{D^2 u}{\Leb{s,q,w}(\mathbb{R}^{d}_T)}\\
&\leq N(\kappa^{-1}+\kappa^{(d+2)/q_0} \delta^{{1}/{(q_0\mu)}} +\kappa^{2(1-1/q_0)}R_1^{2(1-1/q_0)})\norm{D^2 u}{\Leb{s,q,w}(\mathbb{R}^{d}_T)}\\
&\relphantom{=}+N(\kappa^{-1}+\kappa^{(d+2)/q})(\norm{f}{\Leb{s,q,w}(\mathbb{R}^{d}_T)}+\norm{Dg}{\Leb{s,q,w}(\mathbb{R}^{d}_T)})
\end{align*}
for some constant $N=N(d,\nu,s,q,K_0)>0$. Choose $\kappa \geq 8$ large so that $N\kappa^{-1} \leq 1/6$ and choose $0<\delta<1$ so that $N\kappa^{(d+2)/q_0} \delta^{{1}/{(q_0\mu)}} \leq 1/6$. Finally, choose $R_1>0$ so that $N\kappa^{2(1-1/q_0)}R_1^{2(1-1/q_0)} \leq 1/6$. Then we get
\[	\norm{D^2 u}{\Leb{s,q,w}(\mathbb{R}^{d}_T)}\leq N (\norm{f}{\Leb{s,q,w}(\mathbb{R}^{d}_T)}+\norm{Dg}{\Leb{s,q,w}(\mathbb{R}^{d}_T)})\]
for some constant $N=N(d,s,q,\nu,K_0)>0$. This completes the proof of Proposition \ref{prop:Stokes-partition}.
 \end{proof}

 Now we apply ``partition of unity in time'' argument {to remove the assumption that $u$ has compact support in time as in \cite[Lemma 6.5]{DK20}.}

\begin{theorem}\label{thm:mixed-norm-a-priori}
Let $0<T<\infty$, $K_0 \geq 1$, $1<s,q<\infty$, $w=w_1(x)w_2(t)$, where
  \[
    w_1\in A_q(\mathbb{R}^{d},dx),\quad w_2\in A_s(\mathbb{R},dt),\quad [w]_{A_{s,q}}\leq K_0.
  \]
There exists
\[	\delta=\delta(d,\nu,s,q,K_0) \in (0,1)  \]
such that under Assumption \ref{assump:VMO} $(\delta)$, if $(u,p)$ is a strong solution to \eqref{eq:Stokes-nondiv} in $\mathbb{R}^{d}_T$ with $u(0,\cdot)=0$ on $\mathbb{R}^d$ satisfying $u \in \oSob{1,2}{s,q,w}(\mathbb{R}^d_T)^d$, $f\in \Leb{s,q,w}(\mathbb{R}^{d}_T)^d$, and $g \in\Sob{0,1}{s,q,w}(\mathbb{R}^{d}_T)$, then we have
\begin{align*}
\norm{D^2 u}{\Leb{s,q,w}(\mathbb{R}^d_T)}    &\leq N_1\left( \norm{f}{\Leb{s,q,w}(\mathbb{R}^d_T)}+\norm{Dg}{\Leb{s,q,w}(\mathbb{R}^d_T)}\right)+N_2\norm{u}{\Leb{s,q,w}(\mathbb{R}^d_T)},
\end{align*}
where $N_1=N_1(d,s,q,K_0,\nu)>0$ and $N_2=N_2(d,s,q,K_0,\nu,R_0)>0$. Moreover, if $\nabla p \in \Leb{s,q,w}(\mathbb{R}^d_T)^d$, and $g\in\odSob{1}{s,q,w}(\mathbb{R}^d_T)$ and $g_t=\Div G$ for some $G \in \Leb{s,q,w}(\mathbb{R}^d_T)^d$ in the sense of \eqref{eq:compatibility-nondiv}, then there exists a constant $N=N(d,\nu,s,q,K_0,R_0)>0$ such that
\begin{equation}\label{eq:estimate-Stokes-VMO}
\begin{aligned}
  &\norm{\partial_t u}{\Leb{s,q,w}(\mathbb{R}^d_{T})}+\norm{D^2 u}{\Leb{s,q,w}(\mathbb{R}^d_{T})}+\norm{\nabla p}{\Leb{s,q,w}(\mathbb{R}^d_{T})}\\
  &\leq N\left(\norm{f}{\Leb{s,q,w}(\mathbb{R}^d_{T})}+ \norm{Dg}{\Leb{s,q,w}(\mathbb{R}^d_{T})}+\norm{G}{\Leb{s,q,w}(\mathbb{R}^d_{T})}+\norm{u}{\Leb{s,q,w}(\mathbb{R}^d_T)}\right).
\end{aligned}
\end{equation}
\end{theorem}

\begin{proof}
{Take $\delta>0$ and $R_1>0$ given in Proposition \ref{prop:Stokes-partition}.  Choose sequences $t_k\in \mathbb{R}$ and $\{\eta_k(t)\}$ so that $\eta_k\geq 0$, $\eta_k \in C_0^\infty(\mathbb{R})$, $\supp \eta_k\subset (t_k-(R_0 R_1)^2,t_k)$ and
\begin{equation}\label{eq:summation-eta}
1 \leq \sum_{k=1}^\infty |\eta_k(t)|^s \leq \chi_0,\quad \sum_{k=1}^\infty |\eta_k'(t)|^s \leq \chi_1\quad \text{for all } t\in  (0,T),
\end{equation}
where $\chi_0$ depends only on $s$, and $\chi_1$ depends only on $d$, $s$, $R_0$, and $R_1$. }

{
Note that $u_k(t,x):=u(t,x)\eta_k(t)$ and $p_k(t,x):=p(t,x)\eta_k(t)$ satisfies
\[
\left\{\begin{aligned}
(u_k)_t-a^{ij}D_{ij}u_k+\nabla p_k&=\eta_kf+\eta_k'u,\\
\Div u_k &=\eta_k g
\end{aligned}\right.\]
in $\mathbb{R}^d_T$.}
{Then it follows from Proposition \ref{prop:Stokes-partition} that }
\begin{equation}\label{eq:Stokes-partition-unity-1}
\begin{aligned}
\norm{D^2 u_k}{\Leb{s,q,w}(\mathbb{R}^{d}_T)}&\leq N\left( \norm{f\eta_k}{\Leb{s,q,w}(\mathbb{R}^d_T)}+\norm{(Dg)\eta_k}{\Leb{s,q,w}(\mathbb{R}^{d}_T)}+\norm{u\eta_k'}{\Leb{s,q,w}(\mathbb{R}^{d}_T)}\right)
\end{aligned}
\end{equation}
{for some constant $N=N(d,\nu,s,q,K_0)>0$.}

{By summing  \eqref{eq:Stokes-partition-unity-1}  over $k$ and using \eqref{eq:summation-eta}, we get }
\begin{equation}\label{eq:hessian}
	\norm{D^2 u}{\Leb{s,q,w}(\mathbb{R}^d_T)} \leq N_1\left( \norm{f}{\Leb{s,q,w}(\mathbb{R}^d_T)}+\norm{Dg}{\Leb{s,q,w}(\mathbb{R}^d_T)}\right)+N_2 \norm{u}{\Leb{s,q,w}(\mathbb{R}^d_T)},
\end{equation}
where $N_1=N_1(d,s,q,K_0,\nu)>0$ and $N_2=N_2(d,s,q,K_0,\nu,R_0)>0$.

To show \eqref{eq:estimate-Stokes-VMO}, since $(u,p)$ satisfies $u \in \oSob{1,2}{s,q,w}(\mathbb{R}^d_T)^d$, $\nabla p \in \Leb{s,q,w}(\mathbb{R}^d_T)^d$, and
\[
  \partial_t u -\Delta u +\nabla p =f+[a^{ij}(t,x)-\delta^{ij}]D_{ij}u,\quad \Div u=g\quad \text{in } \mathbb{R}^d_T,
\]
where $g_t=\Div G$ in the sense of \eqref{eq:compatibility-nondiv} for some $G\in \Leb{s,q,w}(\mathbb{R}^d_T)^d$, it follows from Theorem \ref{thm:mixed-simple} and \eqref{eq:hessian} that
\begin{align*}
&\norm{\partial_t u}{\Leb{s,q,w}(\mathbb{R}^d_T)}+\norm{D^2 u}{\Leb{s,q,w}(\mathbb{R}^d_T)}+\norm{\nabla p}{\Leb{s,q,w}(\mathbb{R}^d_T)}\\
&\leq N (\norm{f}{\Leb{s,q,w}(\mathbb{R}^d_T)}+\norm{Dg}{\Leb{s,q,w}(\mathbb{R}^d_T)}+\norm{G}{\Leb{s,q,w}(\mathbb{R}^d_T)}+ \norm{D^2 u}{\Leb{s,q,w}(\mathbb{R}^d_T)})\\
&\leq N( \norm{f}{\Leb{s,q,w}(\mathbb{R}^d_T)}+\norm{Dg}{\Leb{s,q,w}(\mathbb{R}^d_T)}+\norm{G}{\Leb{s,q,w}(\mathbb{R}^d_T)}+\norm{u}{\Leb{s,q,w}(\mathbb{R}^d_T)}),
\end{align*}
where $N=N(d,s,q,\nu,K_0,R_0)>0$. This completes the proof of Theorem \ref{thm:mixed-norm-a-priori}.
\end{proof}

The following lemma helps us  absorb the term $\norm{u}{\Leb{s,q,w}(\mathbb{R}^d_T)}$ on the right-hand side in \eqref{eq:estimate-Stokes-VMO} into the left-hand side of \eqref{eq:estimate-Stokes-VMO} which can be easily proved by using fundamental theorem of calculus and Lemma \ref{lem:HL}{, so the} proof is omitted.
\begin{lemma}\label{lem:time-control}
Let $T>0$, $1<s,q<\infty$, $K_0\geq 1$, and let $w(t,x)=w_1(x)w_2(t)$, where $[w_1]_{A_q(\mathbb{R}^d,dx)}\leq K_0$, and $[w_2]_{A_s(\mathbb{R},dt)}\leq K_0$. Then there exists a constant $N=N(d,s,q,K_0)>0$ such that
\[
    \norm{u}{\Leb{s,q,w}(\mathbb{R}^d_T)}\leq N T\norm{\partial_t u}{\Leb{s,q,w}(\mathbb{R}^d_T)}
\]
for all $u\in \oSob{1,2}{s,q,w}(\mathbb{R}^d_T)$.
\end{lemma}

Now we are ready to prove Theorem \ref{thm:A}.
\begin{proof}%[Proof of Theorem \ref{thm:A}]
Since $u\in\oSob{1,2}{s,q,w}(\mathbb{R}^d_T)^d$, $f\in\Leb{s,q,w}(\mathbb{R}^d_T)^d$, and $g\in \odSob{1}{s,q,w}(\mathbb{R}^d_T)$, {we can} extend $u$, $f$, and $g$ to be zero for $t<0$ {so that $u\in\Sob{1,2}{s,q,w}((-\infty,T)\times \mathbb{R}^d)^d$, $f\in\Leb{s,q,w}((-\infty,T)\times \mathbb{R}^d)^d$, and $g\in \dSob{1}{s,q,w}((-\infty,T)\times \mathbb{R}^d)$}. Take a positive integer $m$ to be specified below and set
\[   s_j=\frac{j T}{m},\quad j=-1,0, 1, 2, \dots,m-1 \]
and $\eta_j \in C^\infty(\mathbb{R})$ satisfying
\[     \eta_j(t)=1 \quad \text{for } t \geq s_j,\quad \eta_j(t)=0\quad \text{for } t \leq s_{j-1}\text{ with } |\eta_j'|\leq \frac{2m}{T}.\]

{Note that
\[
\left\{
\begin{aligned}
\partial_t(\eta_j u)-a^{ij}(t,x)D_{ij}(\eta_j u)+\nabla (\eta_j p)&=\eta_j f+\eta_j' u \\
\Div(\eta_j u)&=\eta_j \Div u= \eta_j g
\end{aligned}
\right.
\]
in $\mathbb{R}^d_T$ and $(\eta_j u)(s_{j-1},\cdot)=0$ for $j=0,1,2,\dots,m-1$. }Moreover, $(\eta_j g)_t=\Div \tilde{G}$ for some $\tilde{G} \in \Leb{s,q,w}(\mathbb{R}^d_T)^d$ satisfying
\[
     \norm{\tilde{G}}{\Leb{s,q,w}(\mathbb{R}^d_T)}\leq (1+2m[w_2]_{A_s}^{1/s}) \norm{G}{\Leb{s,q,w}(\mathbb{R}^d_T)}.
\]
Indeed, by the compatibility condition \eqref{eq:compatibility-nondiv}, we have
\begin{equation}\label{eq:g-varphi}
\int_{\mathbb{R}^d_T} g \varphi_t \myd{x}dt=\int_{\mathbb{R}^d_T} G_i D_i\varphi \myd{x}dt
\end{equation}
for all $\varphi \in C_0^\infty([0,T)\times\mathbb{R}^d)$. For $\psi \in C_0^\infty([0,T)\times\mathbb{R}^d)$, put $\varphi(t,x)=\int_t^T \psi(s,x)\myd{s}$ in \eqref{eq:g-varphi}. Then the Fubini theorem gives
\begin{align*}
\int_{\mathbb{R}^d_T} g \eta_j \psi_t \myd{x}dt &=\int_{\mathbb{R}^d} \eta_jG_i D_i \psi\myd{x}dt+\int_{\mathbb{R}^d_T} \left(\int_0^t G_i(s,x)\myd{s}\right) (\eta_j)_t\psi \myd{x}dt
\end{align*}
for all $\psi \in C_0^\infty([0,T)\times\mathbb{R}^d)$. Then we can define
\[    \tilde{G}_i^j(t,x) = \eta_j(t) G_i(t,x)+ \eta_j'(t)\left(\int_0^t G_i(s,x)\myd{s}\right). \]

By the Minkowski integral inequality, H\"older's inequality, and the definition of $A_s$ weights, we have
\begin{align}\label{eq:G-tilde}
    &\norm{\tilde{G}^j_i}{\Leb{s,q,w}(\mathbb{R}^d_T)}\leq  \norm{\eta_j G_i}{\Leb{s,q,w}(\mathbb{R}^d_T)}+\norm*{(\eta_j)_t \int_0^t G_i(\tau,\cdot)\myd{\tau}}{\Leb{s,q,w}(\mathbb{R}^d_T)}\\
    &\leq \norm{G}{\Leb{s,q,w}(\mathbb{R}^d_T)} + \frac{2m}{T} \norm{G}{\Leb{s,q,w}(\mathbb{R}^d_T)} T\left(\fint_0^T w_2^{-1/(s-1)} \myd{\tau} \right)^{{(s-1)}/{s}} \left(\fint_0^T w_2 \myd{\tau}\right)^{1/s}\nonumber\\
    &\leq (1+2m [w_2]_{A_s}^{1/s})\norm{G_i}{\Leb{s,q,w}(\mathbb{R}^d_T)}\nonumber
  \end{align}
  for all $j=0,1,\dots, m$.

To proceed further, for simplicity, we write
\[    \norm{u}{(s_j,s_{j+1})}:=\norm{u}{\Leb{s,q,w}((s_j,s_{j+1})\times\mathbb{R}^d)}.\]
Then by Theorem \ref{thm:mixed-norm-a-priori} and Lemma \ref{lem:time-control}, we have
\begin{align}
&\norm{u}{(s_j,s_{j+1})}\leq \norm{u\eta_j}{(s_{j-1},s_{j+1})}\label{eq:u-iteration}\\
&\leq N \left(\frac{T}{m}\right) \norm{\partial_t (\eta_j u)}{(s_{j-1},s_{j+1})}\nonumber\\
&\leq N \left(\frac{T}{m}\right)\left( \norm{f\eta_j}{(s_{j-1},s_{j+1})}+\norm{\eta_j(Dg)}{(s_{j-1},s_{j+1})}+\norm{\tilde{G}}{(s_{j-1},s_{j+1})}+\norm{u}{(s_j,s_{j+1})}\right),\nonumber
\end{align}
where  $N=N(d,\nu,s,q,K_0,R_0)>0$. Choose $m$ sufficiently large integer $m$ so that $N(T/m)\leq 1/2$. Then by \eqref{eq:G-tilde} and \eqref{eq:u-iteration}, we have
\begin{align*}
  \norm{u}{(s_j,s_{j+1})}& \leq N \left(\norm{f}{(0,s_{j+1})}+\norm{Dg}{(0,s_{j+1})}+\norm{G}{(0,s_{j+1})}+\norm{u}{(0,s_{j})} \right),
\end{align*}
where $N=N(d,\nu,s,q,K_0,R_0,T)>0$ and $j=0,1,\dots,m-1$. By induction and noting that  $\norm{u}{\Leb{s,q,w}((0,s_0)\times\mathbb{R}^d)}=0$, we get
\begin{equation*}%\label{eq:estimate-Stokes-VMO-2}
\begin{aligned}
  &\norm{\partial_t u}{\Leb{s,q,w}(\mathbb{R}^d_{T})}+\norm{u}{\Leb{s,q,w}(\mathbb{R}^d_{T})}+\norm{D^2 u}{\Leb{s,q,w}(\mathbb{R}^d_{T})}+\norm{\nabla p}{\Leb{s,q,w}(\mathbb{R}^d_{T})}\\
  &\leq N\left(\norm{f}{\Leb{s,q,w}(\mathbb{R}^d_{T})}+ \norm{Dg}{\Leb{s,q,w}(\mathbb{R}^d_{T})}+\norm{G}{\Leb{s,q,w}(\mathbb{R}^d_{T})}\right),
\end{aligned}
\end{equation*}
where $N=N(d,\nu,s,q,K_0,R_0,T)>0$. Hence by the method of continuity with Theorem \ref{thm:mixed-simple}, we get the desired solvability results for the problem \eqref{eq:Stokes-nondiv}. This completes the proof of Theorem \ref{thm:A}.
\end{proof}

\section{Stokes equations in divergence form}\label{sec:6}
In this section, we prove Theorem \ref{thm:C}. The proof of Theorem \ref{thm:C} is similar to that of Theorem \ref{thm:A} with some modification. Rather than giving full details of the proof, we highlight the essential differences from the proof of Theorem \ref{thm:A}.

We first obtain a mean oscillation estimate for vorticity $\omega =\nabla \times u$ of weak solutions $u$ to \eqref{eq:Stokes-div}.

\begin{lemma}\label{lem:mixed-norm-div}
Let $\kappa \geq 8$, $1<q<\infty$, $\mu,\mu' \in (1,\infty)$, $1/\mu+1/\mu'=1$, and $a^{ij}$ satisfy Assumption \ref{assump:VMO} $(\delta)$.  Then for any $0< r\leq R_0/\kappa $, $(t_0,x_0) \in \mathbb{R}^{d+1}$, and $u\in \dSobloc{1}{q\mu}(\mathbb{R}^{d+1})^d$ satisfying
\begin{equation*}%\label{eq:Stokes-div-mean2}
  	\partial_t u-D_i(a^{ij}(t,x)D_{j} u)+\nabla p=\Div \boldF,\quad \Div u=g\quad \text{in } Q_{\kappa r}(t_0,x_0),
\end{equation*}
where $\boldF\in \Lebloc{q}(\mathbb{R}^{d+1})^{d\times d}$ and $g\in \Lebloc{q}(\mathbb{R}^{d+1})$, we have
\begin{align*}
&\left(|\omega-(\omega)_{Q_r(t_0,x_0)}|\right)_{Q_r(t_0,x_0)}\\
&\leq N \kappa^{-1} \left[(|D u|^q)^{1/q}_{Q_{\kappa r}(t_0,x_0)} +(|\boldF|^{q})^{1/q}_{Q_{\kappa r}(t_0,x_0)}\right]\\
&\relphantom{=}+N\kappa^{(d+2)/q} \left[ (|\boldF|^{q})_{Q_{\kappa r}(t_0,x_0)}^{{1}/{q}} +(|g|^{q})_{Q_{\kappa r}(t_0,x_0)}^{{1}/{q}}+ \delta^{1/(q\mu')}(|D u|^{q\mu})_{Q_{\kappa r}(t_0,x_0)}^{{1}/{(q\mu)}}\right]
\end{align*}
 for some constant $N=N(d,q,\nu)>0$.
\end{lemma}
\begin{proof}
The proof is essentially the same as that of Lemma \ref{lem:mixed-norm}. The main difference is to apply Theorem \ref{thm:mixed-simple-div} instead of Theorem \ref{thm:mixed-simple}.

By translation invariance, we may assume that $(t_0,x_0)=(0,0)$. For a locally integrable function $h$ defined on $Q_{\kappa r}$, let $h^{(\varepsilon)}$ denote the standard mollification defined in \eqref{eq:h-mollification}. Fix $0<r'<\kappa r$ and let $\varphi \in C_{0}^\infty(Q_{r'})$. Then for small ${\varepsilon>0}$ satisfying $\varepsilon^2<(\kappa r)^2-(r')^2$, it is easy to verify that
\[
\int_{Q_{r'}} h^{(\varepsilon)} (t,x) \partial_t \varphi(t,x) \myd{t}dx=\int_{Q_{\kappa r}} h(s,y)\partial_t \tilde{\varphi}^{(\varepsilon)}(s,y)\myd{s}dy,
\]
where
\begin{equation}\label{eq:forward-mollification}
     \tilde{\varphi}^{(\varepsilon)}(s,y)=\int_{Q_{r'}} \eta_\varepsilon(s-t,y-x) \varphi(t,x)\myd{tdx}.
\end{equation}
Then for $\varphi \in C_{0}^\infty(Q_{r'})^d$ satisfying $\Div \varphi(t,\cdot)=0$ in $B_{r'}$ for all $t\in (-(r')^2,0)$, using $\tilde{\varphi}^{(\varepsilon)}$ as a test function in the definition of weak solutions, we get
\begin{align*}
&{-}\int_{Q_{r'}} u^{(\varepsilon)} \cdot \partial_t \varphi \myd{z}+\int_{Q_{r'}} a^{ij} D_j u^{(\varepsilon)} \cdot D_i \varphi \myd{x}\\
&={-}\int_{Q_{r'}} \boldF^{(\varepsilon)}:\nabla \varphi \myd{z}{-}\int_{Q_{r'}} [a^{ij}D_j u^{(\varepsilon)}-(a^{ij}D_j u)^{(\varepsilon)}] \cdot D_i\varphi  \myd{z}.
\end{align*}
In other words, $u^{(\varepsilon)}$ is a weak solution to
\[
     \partial_t u^{(\varepsilon)}-D_i(a^{ij}D_{j}u^{(\varepsilon)})+\nabla p^{(\varepsilon)}=\Div (\boldF^{(\varepsilon)}+\boldH^{\varepsilon})\quad \text{in } Q_{r'},
\]
where $\boldH^{\varepsilon}=(H_1^\varepsilon,\dots,H_d^\varepsilon)$ and $H_i^\varepsilon:=(a^{ij}D_j u)^{(\varepsilon)}-a^{ij}(t,x)D_j u^{(\varepsilon)}$ in $Q_{r'}$. If we prove the estimate in the lemma for $u^{(\varepsilon)}$, we get the desired result by letting $\varepsilon\rightarrow 0$. Hence we may assume that $u$ and $g$ are infinitely differentiable.

Let  $\zeta_r(x)$ and $\psi_r(t)$ be  infinitely differentiable functions defined on $\mathbb{R}^d$ and $\mathbb{R}$ satisfying $0\leq \zeta_r \leq 1$, $0\leq \psi_r \leq 1$,
\[
  \zeta_r(x)=1\quad\text{on } B_{2r/3},\quad \zeta_r(x)=0\quad \text{on } \mathbb{R}^d\setminus B_r,
\]
\[
  \psi_r(t)=1\quad\text{on } t\in (-4r^2/9,4r^2/9),\quad \psi_r(t)=0\quad \text{on } \mathbb{R}\setminus (-r^2,r^2).
\]
Set $\phi_r(t,x)=\psi_r(t)\zeta_r(x)$. Then $\phi_r=1$ on $Q_{2r/3}$ and $|D\phi_r|\leq 4/r$.

Since
\[    \int_{B_{\kappa r}} (g-[g(t,\cdot)]_{\zeta_{\kappa r},B_{\kappa r}} )\phi_{\kappa r } \myd{x}=0,\]
for each $t\in [-(\kappa r)^2,0)$, it follows from Theorem \ref{thm:div-equation} that there exists $G(t,\cdot)\in \Sobd{1}{q}(B_{\kappa r})^d$ satisfying
\[
    \Div G =(g-[g(t,\cdot)]_{\zeta_{\kappa r},B_{\kappa r}} )\phi_{\kappa r } \quad \text{in } B_{\kappa r},\quad G=0\quad \text{on } \partial B_{\kappa r}
\]
and
\begin{equation}\label{eq:G-estimate}
\norm{DG(t,\cdot)}{\Leb{q}(B_{\kappa r})}\leq N(d,q) \norm{(g(t,\cdot)-[g(t,\cdot)]_{\zeta_{\kappa r},B_{\kappa r}})\phi_{\kappa r}(t,\cdot)}{\Leb{q}(B_{\kappa r})}
\end{equation}
for $t\in (-(\kappa r)^2,0)$ and $G(-(\kappa r)^2,\cdot)=0$ on $B_{\kappa r}$. Hence by \eqref{eq:G-estimate} and H\"older's inequality, we have
\begin{equation}\label{eq:DG-control}
    \norm{D G}{\Leb{q}(Q_{\kappa r})}\leq N(d,q) \norm{g}{\Leb{q}(Q_{\kappa r})}.
\end{equation}

Now we choose $\varphi \in C_0^\infty (B_2)$ so that $\varphi =0$ in $B_1$ and $\int_{B_2} \varphi \myd{x}=1$. Define $\overline{G}^i(t,x)$ by
\[
\overline{G}^i(t,x):=\begin{dcases}
G^i(t,x)&\quad \text{in } B_{\kappa r},\\
c^i(t)\varphi\left(\frac{x}{\kappa r}\right)&\quad \text{in } B_{2\kappa r}\setminus B_{\kappa r},\\
0 &\quad \text{in } B_{2\kappa r}^c,
\end{dcases}
\]
where $c^i(t){=-(\int_{B_{\kappa r}} G^i(t,x)\myd{x})(\int_{B_{2\kappa r}} \varphi\left(\frac{x}{\kappa r}\right)\myd{x})^{-1}}$  so that
\[        \int_{B_{2\kappa r}} \overline{G}^i(t,x)\myd{x}=0.\]
Define
\[    h(t,x):=\Div \overline{G}(t,x).\]
Note that $h=(g-[g]_{B_{\kappa r,\zeta_{\kappa r}}})\phi_{\kappa r}$ in $Q_{\kappa r}$. Since $g$ is infinitely differentiable in $t$, $h$ is also infinitely differentiable in $t$. Moreover, since $\int_{B_{2 \kappa r}} \overline{G}^i \myd{x}=0$, it follows that $\int_{ B_{2\kappa r}} \partial_t \overline{G}^i \myd{x}=0$. Hence by  Theorem \ref{thm:div-equation}, there exists $H^i \in \Sobd{1}{q}(B_{2\kappa r})$ satisfying
\[
      \Div H^i=\partial_t \overline{G}^i\quad \text{in } B_{2\kappa r},\quad H^i=0\quad \text{on } \partial B_{2\kappa r}
\]
for each $t$. Extend $H^i$ to be zero outside $B_{2\kappa r}$. Since $\partial_t \overline{G}^i$ has compact support in $B_{2\kappa r}$ for each $t \in (-(\kappa r)^2,0)$ and $\overline{G}^i(-(\kappa r)^2,\cdot)=0$, we have
\[
    \int_{-(\kappa r)^2}^0 \int_{\mathbb{R}^d} (\Div H^i)\phi \myd{x}dt= \int_{-(\kappa r)^2}^0 \int_{\mathbb{R}^d} (\partial_t \overline{G}^i)\phi \myd{x}dt=-\int_{-(\kappa r)^2}^0 \int_{\mathbb{R}^d} \overline{G}^i \partial_t \phi \myd{x}dt
\]
for all $\phi \in C_0^\infty([-(\kappa r)^2,0)\times \mathbb{R}^d)$. For all $\psi \in  C_0^\infty([-(\kappa r)^2,0)\times \mathbb{R}^d)$, taking $\phi=D_i \psi$, we have
\[
 -\int_{-(\kappa r)^2}^0 \int_{\mathbb{R}^d} H^i \cdot \nabla (D_i\psi)\myd{x}dt=\int_{-(\kappa r)^2}^0 \int_{\mathbb{R}^d} h \partial_t \psi \myd{x}dt.
\]
Hence $h$ satisfies the compatibility condition \eqref{eq:compatibility-divergence}.

By a change of variables, H\"older's inequality, and the Poincar\'e inequality, we get
\begin{equation}\label{eq:c-G-control}
|c^i(t)|\leq \frac{N(d,q)}{(\kappa r)^{d/q}} \norm{G^i(t,\cdot)}{\Leb{q}(B_{\kappa r})}\leq \frac{N(d,q)}{(\kappa r)^{d/q}} (\kappa r) \norm{DG^i(t,\cdot)}{\Leb{q}(B_{\kappa r})}
\end{equation}
Moreover, it follows from \eqref{eq:DG-control}, \eqref{eq:c-G-control}, and a change of variable that
\begin{align*}
     &\norm{h}{\Leb{q}((-(\kappa r)^2,0)\times \mathbb{R}^d)}\\
     &\leq N(d,q)\left(\norm{DG}{\Leb{q}(Q_{\kappa r})}+\frac{1}{\kappa r}\norm{c(t)D\varphi(x/(\kappa r)) }{\Leb{q}((-(\kappa r)^2,0)\times (B_{2\kappa r}\setminus B_{\kappa r}))}\right)\\
     &\leq N(d,q)\norm{g}{\Leb{q}(Q_{\kappa r})}.
\end{align*}

Now consider the following initial-value problem for Stokes equations:  for $l=1,\dots,d$,
\begin{equation*}%\label{eq:Stokes-w-part-div}
   \left\{\begin{alignedat}{2}
\partial_t u_1^l-D_{i}(\overline{a}^{ij}(t) D_ju_1^l) + D_l p_1&=D_i[1_{Q_{\kappa r}}H^{ij,l}]&&\quad \text{in } (-(\kappa r)^2,0)\times \mathbb{R}^d,\\
\Div u_1 &=h&&\quad \text{in } (-(\kappa r)^2,0)\times \mathbb{R}^d,\\
u_1&=0&&\quad \text{on } \{t=-\kappa r^2\}\times \mathbb{R}^d,
\end{alignedat}
\right.
\end{equation*}
where $H^{ij,l}=[F^{li}+ (a^{ij}-\overline{a}^{ij}(t))D_{j} u^l]$.
If we define $u_2=u-u_1$ and $p_2=p-p_1$, then $(u_2,p_2)$ satisfies
\begin{equation*}%\label{eq:Stokes-w-part-div}
   \left\{\begin{alignedat}{1}
\partial_t u_2-D_{i}(\overline{a}^{ij}(t) D_ju_2) + \nabla p_2&=0\\
\Div u_2 &={[g(t,\cdot)]_{\zeta_{\kappa r},B_{\kappa r}}}
\end{alignedat}
\right.\,\,\,\text{in } Q_{2\kappa r/3}.
\end{equation*}
Since $h$ satisfies the compatibility condition, following exactly the same argument as in {the proof of} Lemma \ref{lem:mixed-norm}, we can prove the desired result.
\end{proof}

Following exactly the same argument as in Proposition \ref{prop:Stokes-partition} using Lemma \ref{lem:mixed-norm-div} instead of Lemma \ref{lem:mixed-norm}, one can prove the following proposition of which proof is omitted. This proposition is necessary to perform a partition of unity in time argument.
\begin{proposition}\label{prop:Stokes-partition-2}
Let $0<T<\infty$, $K_0\geq 1$, $1<s,q<\infty$, $t_1 \in \mathbb{R}$, $w=w_1(x)w_2(t)$, where
  \[
    w_1 \in A_q(\mathbb{R}^{d},dx),\quad w_2\in A_s(\mathbb{R},dt),\quad [w]_{A_{s,q}}\leq K_0.
  \]
 Then there exist $\delta>0$ and $R_1>0$ such that under Assumption \ref{assump:VMO} $(\delta)$, for $u\in \odSob{1}{s,q,w}(\mathbb{R}^{d}_T)^d$ vanishing outside $(t_1-(R_0R_1)^2,t_1)\times \mathbb{R}^d$ and is a weak solution to
\[	\partial_t u-D_i(a^{ij}(t,x)D_{j} u)+\nabla p=\Div \boldF,\quad \Div u=g\quad \text{in } \mathbb{R}^{d}_T,\]
where $\boldF\in \Leb{s,q,w}(\mathbb{R}^{d}_T)^{d\times d}$ and $g\in \Leb{s,q,w}(\mathbb{R}^{d}_T)$, then  there exists a constant $N=N(d,s,q,K_0,\nu)>0$ such that
\[	\norm{Du}{\Leb{s,q,w}(\mathbb{R}^{d}_T)}\leq N \left(\norm{\boldF}{\Leb{s,q,w}(\mathbb{R}^{d}_T)}+\norm{g}{\Leb{s,q,w}(\mathbb{R}^{d}_T)}\right).\]
\end{proposition}

Using Proposition \ref{prop:Stokes-partition-2}, we can prove the gradient estimate of weak solutions  by following exactly the same argument as in Theorem \ref{thm:mixed-norm-a-priori}.
\begin{theorem}\label{thm:mixed-norm-a-priori-div}
Let $0<T<\infty$, $K_0 \geq 1$, $1<s,q<\infty$, $w(t,x)=w_1(x)w_2(t)$, $[w_1]_{A_q(\mathbb{R}^d,dx)}\leq K_0$, and $[w_2]_{A_s(\mathbb{R},dt)}\leq K_0$. There exists $\delta=\delta(d,\nu,s,q,K_0)>0$ such that under Assumption \ref{assump:VMO} $(\delta)$, if $u \in \odSob{1}{s,q,w}(\mathbb{R}^d_T)^d$ is a weak solution to \eqref{eq:Stokes-div} in $\mathbb{R}^d_T$ with $u(0,\cdot)=0$ on $\mathbb{R}^d$ satisfying $g \in \Leb{s,q,w}(\mathbb{R}^d_T)$ and $\boldF \in \Leb{s,q,w}(\mathbb{R}^d_T)^{d\times d}$, then
\begin{align*}
\norm{D u}{\Leb{s,q,w}(\mathbb{R}^d_T)}    &\leq N_1(\norm{\boldF}{\Leb{s,q,w}(\mathbb{R}^d_T)}+\norm{g}{\Leb{s,q,w}(\mathbb{R}^d_T)})+N_2\norm{u}{\Leb{s,q,w}(\mathbb{R}^d_T)},
\end{align*}
where $N_1=N_1(d,s,q,K_0,\nu)>0$ and $N_2=N_2(d,s,q,K_0,\nu,R_0)>0$. Moreover, if $p\in \Leb{s,q,w}(\mathbb{R}^d_T)$, $g_t = \Div\Div\boldG$ for some $\boldG \in \Leb{s,q,w}(\mathbb{R}^d_T)^{d\times d}$ in the sense of \eqref{eq:compatibility-div}, then \begin{equation}\label{eq:estimate-Stokes-VMO-3}
\begin{aligned}
&\norm{Du}{\Leb{s,q,w}(\mathbb{R}^d_T)} +\norm{p}{\Leb{s,q,w}(\mathbb{R}^d_T)}\\
&\leq N_1(\norm{\boldF}{\Leb{s,q,w}(\mathbb{R}^d_T)}+\norm{g}{\Leb{s,q,w}(\mathbb{R}^d_T)}+\norm{\boldG}{\Leb{s,q,w}(\mathbb{R}^d_T)})+N_2\norm{u}{\Leb{s,q,w}(\mathbb{R}^d_T)}.
\end{aligned}
\end{equation}
\end{theorem}

The key difference from the nondivergence form case is that it is hard to immediately absorb the term $\norm{u}{\Leb{s,q,w}(\mathbb{R}^d_T)}$ into the left-hand side of \eqref{eq:estimate-Stokes-VMO-3} by using the time derivative of $u$ as in Lemma \ref{lem:time-control} (see also the proof of  Theorem \ref{thm:mixed-norm-a-priori}). Nevertheless, using a mollification argument, we can absorb the term into the left-hand side. See, for instance, Dong-Liu \cite{DLiu22}. For the sake of completeness, we explain it in detail.

Choose {a radially symmetric} $\zeta \in C_0^\infty(\mathbb{R}^d)$ satisfying $\supp \zeta \subset B_1$ and $\int_{\mathbb{R}^d} \zeta \myd{x}=1$, and for $\varepsilon>0$, let $\zeta^\varepsilon(x):=\varepsilon^{-d} \zeta(x/\varepsilon)$. The following lemma can be easily proved by the definition of mollification, fundamental theorem of calculus, and Lemma \ref{lem:HL} (see \cite[Lemma A.2]{DLiu22}).
\begin{lemma}\label{lem:mollification-weight}
Let  $1<q<\infty$, $K_0 \geq 1$, and $w_1\in A_q(\mathbb{R}^d,dx)$ with $[w_1]_{A_q}\leq K_0$. For $\varepsilon>0$ and $v\in \Sob{1}{q,w_1}(\mathbb{R}^d)$, define
\[  v^{(\varepsilon)}(x):=(v*\zeta^\varepsilon)(x)=\int_{\mathbb{R}^d} \zeta^\varepsilon(x-y)v(y)\myd{y}.\]
 Then there exists a constant $N=N(d,q,K_0)>0$ such that
\[
      \norm{v^{(\varepsilon)}-v}{\Leb{q,w_1}(\mathbb{R}^d)}\leq N \varepsilon \norm{Dv}{\Leb{q,w_1}(\mathbb{R}^d)}.
\]
\end{lemma}

\begin{lemma}\label{lem:mollification-x-to-get-t}
Let $1<s,q<\infty$, $K_0\geq 1$, and let $w(t,x)=w_1(x)w_2(t)$, $[w_1]_{A_q(\mathbb{R}^d,dx)}\leq K_0$, and $[w_2]_{A_s(\mathbb{R},dt)}\leq K_0$.  If $u\in \dSob{1}{s,q,w}(\mathbb{R}^d_T)$ satisfies
\begin{equation}\label{eq:u-distribution}
  \partial_t u=f + \Div F\quad \text{in } \mathbb{R}^d_T
\end{equation}
for some $f \in \Leb{s,q,w}(\mathbb{R}^d_T)$ and $F=(F^1,\dots,F^d)\in\Leb{s,q,w}(\mathbb{R}^d_T)^d$, then $u^{(\varepsilon)} := u*\zeta^\varepsilon \in \Sob{1,2}{s,q,w}(\mathbb{R}^d_T)$ and for any $0<\varepsilon<1$, we have
\begin{equation}\label{eq:mollification-t-estimate}
   \norm{\partial_t u^{(\varepsilon)}}{\Leb{s,q,w}(\mathbb{R}^d_T)}\leq\frac{N}{\varepsilon} \norm{F}{\Leb{s,q,w}(\mathbb{R}^d_T)}+N\norm{f}{\Leb{s,q,w}(\mathbb{R}^d_T)},
\end{equation}
where $N=N(d,s,q,K_0)>0$.
\end{lemma}
\begin{proof}
Since $u$ satisfies \eqref{eq:u-distribution}, we have
\begin{equation}\label{eq:u-phi}
-\int_{\mathbb{R}^d_T}u \partial_t \phi \myd{x}dt=\int_{\mathbb{R}^d_T} f\phi \myd{x}dt-\int_{\mathbb{R}^d_T} F \cdot \nabla \phi\myd{x}dt
\end{equation}
for all $\phi \in C_0^\infty(\mathbb{R}^d_T)$. Put $\phi=\psi*\zeta^\varepsilon$ in \eqref{eq:u-phi}, where $\psi \in C_0^\infty(\mathbb{R}^d_T)$. Then the Fubini theorem gives
\begin{equation}
				\label{eq11.40}
-\int_{\mathbb{R}^d_T}u^{(\varepsilon)} \partial_t {\psi} \myd{x}dt=\int_{\mathbb{R}^d_T} f^{(\varepsilon)}{\psi} \myd{x}dt-\int_{\mathbb{R}^d_T} F^{(\varepsilon)} \cdot \nabla {\psi}\myd{x}dt. 	
\end{equation}
Note that
\begin{align*}
|f*\zeta^\varepsilon(t,x)|&=\left|\int_{B_1} f(t,x-\varepsilon y) \zeta(y)\myd{y} \right|\\
&\leq \norm{\zeta}{\Leb{\infty}}|B_1| \fint_{B_1} |f(t,x-\varepsilon y)|\myd{y}\\
&\leq \norm{\zeta}{\Leb{\infty}} |B_1| {\mathcal{M}^xf(t,x)},
\end{align*}
for $(t,x)\in\mathbb{R}^d_T$, where $\mathcal{M}^x$ is the Hardy-Littlewood maximal function in $x$.

By Lemma \ref{lem:HL}, we have
\[
     \norm{f*\zeta^\varepsilon}{\Leb{s,q,w}(\mathbb{R}^d_T)}\leq N(d,s,q,K_0) \norm{f}{\Leb{s,q,w}(\mathbb{R}^d_T)}
\]
for $0<\varepsilon<1$. Similarly, we have
\[
     \norm{F^i*(D_i\zeta^\varepsilon)}{\Leb{s,q,w}(\mathbb{R}^d_T)}\leq \frac{N(d,s,q,K_0)}{\varepsilon} \norm{D_i \zeta}{\Leb{\infty}(\mathbb{R}^d)} \norm{F}{\Leb{s,q,w}(\mathbb{R}^d_T)}
\]
for $0<\varepsilon<1$. Hence it follows from {\eqref{eq11.40}}  and H\"older's inequality that
\begin{align*}
\left|\int_{\mathbb{R}^d_T} u^{(\varepsilon)}\partial_t\psi\myd{x}dt \right|&\leq N(d,s,q,K_0)\left(\norm{f}{\Leb{s,q,w}(\mathbb{R}^d_T)}+\frac{1}{\varepsilon}\norm{F}{\Leb{s,q,w}(\mathbb{R}^d_T)}\right)\norm{\psi}{\Leb{s',q',\tilde{w}}(\mathbb{R}^d_T)},
\end{align*}
where $\tilde{w}=w_1^{-1/(q-1)}w_2^{-1/(s-1)}$ and for all ${\psi}\in C_0^\infty(\mathbb{R}^d_T)$. Therefore by duality, $\partial_t u^{(\varepsilon)}$ satisfies  \eqref{eq:mollification-t-estimate}. This completes the proof of Lemma \ref{lem:mollification-x-to-get-t}.
\end{proof}

Now we are ready to prove Theorem \ref{thm:C}.
\begin{proof}
By Theorem \ref{thm:mixed-norm-a-priori-div} and the method of continuity together with Theorem \ref{thm:mixed-simple-div}, it suffices to show that there exists a constant $N=N(d,s,q,\nu,K_0,T)>0$ such that
\[
    \norm{u}{\Leb{s,q,w}(\mathbb{R}^d_T)} \leq N (\norm{\boldF}{\Leb{s,q,w}(\mathbb{R}^d_T)}+\norm{g}{\Leb{s,q,w}(\mathbb{R}^d_T)}+\norm{\boldG}{\Leb{s,q,w}(\mathbb{R}^d_T)}).
\]

Since $u\in \odSob{1}{s,q,w}(\mathbb{R}^d_T)$, $|\boldF|,g\in\Leb{s,q,w}(\mathbb{R}^d_T)$, {we can} extend $u$, $\boldF$, and $g$ to be zero for $t<0$. Take a positive integer $m$ to be specified below and set
\[   s_j=\frac{j T}{m},\quad j=-1,0, 1, 2, \dots,m-1 \]
and $\eta_j \in C^\infty(\mathbb{R})${, where $\eta_j$ is defined in the proof of Theorem \ref{thm:A}.}  It is easy to see that $\eta_k u \in \odSob{1}{s,q,w}((s_{k-1},T)\times\mathbb{R}^d)^d$ for $k=0,\dots,m-1$.

{
Note that
\begin{equation}\label{eq:Stokes-cutoff-div}
\left\{
\begin{aligned}
\partial_t (\eta_k u)-D_i (a^{ij}D_j(\eta_k u))+\nabla (\eta_k p)&=\Div (\eta_k \boldF)+\eta_k' u,\\
\Div(\eta_k u)&=\eta_k \Div u= \eta_k g
\end{aligned}
\right.
\end{equation}
in $\mathbb{R}^d_T$ and $(\eta_k u)(s_{k-1},\cdot)=0$ on $\mathbb{R}^d$ for $k=0,1,\dots,m-1$.} Similar to the nondivergence case, one can show that there exists $\tilde{\mathbf{G}}^k \in \Leb{s,q,w}(\mathbb{R}^d_T)^{d\times d}$ satisfying $(\eta_k g)_t=\Div\Div \tilde{\mathbf{G}}^k$ in the sense of \eqref{eq:compatibility-div} and
\begin{equation}\label{eq:compatibility-div-control}
     \norm{\tilde{\mathbf{G}}^k}{\Leb{s,q,w}(\mathbb{R}^d_T)}\leq (1+2m[w_2]_{A_s}^{1/s}) \norm{\mathbf{G}}{\Leb{s,q,w}(\mathbb{R}^d_T)}.
\end{equation}

  For simplicity, we write
  \[    \norm{u}{(s_k,s_{k+1})}=\norm{u}{\Leb{s,q,w}((s_k,s_{k+1})\times \mathbb{R}^d)}.\]
 Since $(u,p)$ satisfies \eqref{eq:Stokes-cutoff-div}, it follows from Lemma \ref{lem:mollification-x-to-get-t} that
\begin{equation}\label{eq:mollification-time-control-last}
\begin{aligned}
 \norm{\partial_t(\eta_k u)^{(\varepsilon)}}{(0,T)} &\leq \frac{N}{\varepsilon}\left(\norm{\eta_k \boldF}{(0,T)}+\norm{D_j(\eta_k u)}{(0,T)}+\norm{\eta_k p}{(0,T)}\right)\\
   &\relphantom{=}+N \norm{\eta_k' u}{(0,T)}
\end{aligned}
\end{equation}
{for some $N=N(d,s,q,K_0,\nu)>0$.}

By Theorem \ref{thm:mixed-norm-a-priori-div}, \eqref{eq:Stokes-cutoff-div}, and \eqref{eq:compatibility-div-control}, there exists a constant $N=N(d,s,q,\nu,K_0,R_0)>0$ such that
\begin{equation}\label{eq:ueta-k}
\begin{aligned}
&\norm{D(u\eta_k)}{(s_k,s_{k+1})}+\norm{\eta_k p}{(s_k,s_{k+1})}\\
&\leq \norm{Du}{(0,s_{k+1})}+\norm{p}{(0,s_{k+1})}\\
&\leq N \left( \norm{\boldF}{(0,T)}+(1+2m[w_2]_{A_s}^{1/s})\norm{\boldG}{(0,T)}+\norm{g}{(0,T)}+\norm{u}{(0,s_{k+1})}\right)
\end{aligned}
\end{equation}
for all $k=0,\dots,m-1$. Hence it follows from Lemmas \ref{lem:mollification-weight}, \ref{lem:time-control}, \eqref{eq:ueta-k}, and \eqref{eq:mollification-time-control-last} that {
\begin{align*}
&\norm{u}{(s_k,s_{k+1})}\\
&\leq \norm{(u\eta_k)^{(\varepsilon)}-u\eta_k}{(s_{k-1},s_{k+1})}+\norm{(u\eta_k)^{(\varepsilon)}}{(s_{k-1},s_{k+1})}\\
&\leq N\varepsilon \norm{D(u\eta_k)}{(s_{k-1},s_{k+1})}+N\left(\frac{T}{m}\right)\norm{\partial_t (u\eta_k)^{(\varepsilon)}}{(s_{k-1},s_{k+1})}\\
&\leq N\varepsilon\left(\norm{\boldF}{(0,T)}+(1+2m[w_2]_{A_s}^{1/s})\norm{\boldG}{(0,T)}+\norm{g}{(0,T)}+\norm{u}{(0,s_{k+1})} \right)\\
&\relphantom{=}+N\left(\frac{T}{m\varepsilon}\right)\left(\norm{D_j(\eta_k u)}{(s_{k-1},s_{k+1})}+\norm{\eta_k p}{(s_{k-1},s_{k+1})}+\norm{\eta_k \boldF}{(0,T)}\right)\\
&\relphantom{=}+N\left(\frac{T}{m}\right)\norm{\eta_k' u}{(0,T)}\\
&\leq N\varepsilon\left(\norm{\boldF}{(0,T)}+(1+2m[w_2]_{A_s}^{1/s})\norm{\boldG}{(0,T)}+\norm{g}{(0,T)}+\norm{u}{(0,s_{k})}+\norm{u}{(s_k,s_{k+1})} \right)\\
&\relphantom{=}+N\left(\frac{T}{m\varepsilon}\right)\left(\norm{\boldF}{(0,T)}+\norm{u}{(0,s_{k})}+\norm{u}{(s_k,s_{k+1})}\right)+N\norm{u}{({s_{k-1},s_{k}})}\\
&\leq N\left(\varepsilon +\frac{T}{m\varepsilon}  \right)\norm{u}{(s_{k},s_{k+1})}\\
&\relphantom{=}+N\varepsilon\left(\norm{\boldF}{(0,T)}+(1+2m[w_2]_{A_s}^{1/s})\norm{\boldG}{(0,T)}+\norm{g}{(0,T)}+\norm{u}{(0,s_k)} \right)\\
&\relphantom{=}+N\left(\frac{T}{m\varepsilon} \right)\left( \norm{\boldF}{(0,T)}+\norm{u}{(0,s_k)}\right)
{+N\norm{u}{(s_{k-1},s_{k})}}
\end{align*} }
for some constant $N=N(d,s,q,\nu,K_0,R_0)>0$. {Choose} $\varepsilon>0$ sufficiently small and then choose $m$ sufficiently large so that $\norm{u}{(s_k,s_{k+1})}$ is absorbed into the left-hand side. Then we have
\[
    \norm{u}{(s_k,s_{k+1})}\leq N\left(\norm{\boldF}{(0,T)}+\norm{\boldG}{(0,T)}+\norm{g}{(0,T)}+\norm{u}{(0,s_k)}\right)+{N\norm{u}{(s_{k-1},s_{k})}}
\]
for some constant $N=N(d,s,q,\nu,K_0,T,R_0)>0$ and $k=0,\dots,m-1$. By induction and noting that $\norm{u}{(0,s_0)}=0$, we get
\[
    \norm{u}{(0,T)}\leq N\left(\norm{\boldF}{(0,T)}+\norm{\boldG}{(0,T)}+\norm{g}{(0,T)}\right)
\]
for some constant $N=N(d,s,q,\nu,K_0,T,R_0)>0$. This completes the proof of Theorem \ref{thm:C}.
   \end{proof}

\section{Interior mixed-norm derivative estimates for Stokes equations}\label{sec:7}

This section is devoted to proving Theorems \ref{thm:B} and \ref{thm:D}, which will be given in Subsections \ref{subsec:Hessian} and \ref{subsec:gradient}, respectively.

\subsection{Interior Hessian estimates for Stokes equations in nondivergence form}\label{subsec:Hessian}
To prove Theorem \ref{thm:B}, we use the following theorem that was implicitly proved in Dong-Phan \cite[Theorem 1.11]{DP21}.

\begin{theorem}\label{thm:DP}
Let $1<s,q<\infty$, $\nu \in (0,1)$, and $1/2\leq r<R\leq 1$. There exists $\delta=\delta(d,\nu,s,q)\in (0,1)$ such that under Assumption \ref{assump:VMO} $(\delta)$, if $(u,p)\in \tilde{W}^{1,2}_{s,q}(Q_R)^d\times\Sob{0,1}{1}(Q_R)$ is a strong solution to \eqref{eq:Stokes-nondiv} in $Q_R$, $f\in \Leb{s,q}(Q_R)^d$ and $g\in \Sob{0,1}{s,q}(Q_R)$, then there exists a constant $N=N(d,s,q,\nu,r,R,R_0)>0$ such that
\begin{align*}
\norm{D^2u}{\Leb{s,q}(Q_{r})}\leq N\left[\frac{1}{(R-r)^b}\norm{u}{\Leb{s,1}(Q_R)}+\norm{f}{\Leb{s,q}(Q_R)}+\norm{Dg}{\Leb{s,q}(Q_R)} \right]
\end{align*}
for some constant $b=b(d,q)>2$.
\end{theorem}

For the sake of completeness, we give a proof of Theorem \ref{thm:DP} {by using following lemma}.
\begin{lemma}[{\cite[Lemma 4.13]{DP21}}]\label{lem:DP}
Let $1/2\leq R< 1$, $R_1\in (0,R_0)$, $R_1 \in (0,R_0)$, $\delta \in (0,1)$, $\kappa \in (0,1/4)$, $1<s,q<\infty$, $q_1 \in (1,\min\{s,q\})$, and $1<q_0<q_1$. Suppose that $(u,p)\in \tilde{W}^{1,2}_{s,q}(Q_{R+R_1})^d\times \Sob{0,1}{1}(Q_{R+R_1})$ is a strong solution to \eqref{eq:Stokes-nondiv} in $Q_{R+R_1}$ {for some $f\in \Leb{s,q}(Q_{R+R_1})^d$ and $g \in \Sob{0,1}{s,q}(Q_{R+R_1})$}. Then
\begin{align*}
\norm{D^2 u}{\Leb{s,q}(Q_R)}&\leq N \kappa^{-(d+2)/q_0}\norm{f}{\Leb{s,q}(Q_{R+R_1})}+N\kappa^{-(d+2)/q_0}\norm{Dg}{\Leb{s,q}(Q_{R+R_1})} \\
&\relphantom{=}+N\left(\kappa^{-(d+2)/q_0}\delta^{1/q_0-1/q_1}+\kappa \right)\norm{D^2 u}{\Leb{s,q}(Q_{R+R_1})}\\
&\relphantom{=}+ N\kappa^{-{(d+2)}/{q_0}}R_1^{-1} \norm{Du}{\Leb{s,q}(Q_{R+R_1})}
\end{align*}
for some constant $N=N(d,s,q,\nu)>0$.
\end{lemma}	

\begin{proof}
Fix $1/2\leq r<R\leq 1$. For $k=1,2,\dots,$ we write
\[   r_k=R-\frac{R-r}{2^{k-1}},\quad k=1,2,\dots.\]
Then $r_1=r$ and $r_k$ is increasing satisfying $\lim_{k\rightarrow\infty} r_k=R$. Let $k_0$ be the smallest positive integer $k$ such that $2^{-k}{(R-r)}\leq R_0$. For $k\geq k_0$, we apply Lemma \ref{lem:DP} with $R=r_k$ and $R_1=2^{-k}(R-r)$. Since $r_k+R_1=r_{k+1}$, we get
\begin{align*}
\norm{D^2 u}{\Leb{s,q}(Q_{r_k})}&\leq N \kappa^{-(d+2)/q_0}\norm{f}{\Leb{s,q}(Q_{r_{k+1}})}+N\kappa^{-(d+2)/q_0}\norm{Dg}{\Leb{s,q}(Q_{r_{k+1}})} \\
&\relphantom{=}+N\left(\kappa^{-(d+2)/q_0}\delta^{1/q_0-1/q_1}+\kappa \right)\norm{D^2 u}{\Leb{s,q}(Q_{r_{k+1}})}\\
&\relphantom{=}+ N\kappa^{-{(d+2)}/{q_0}} \frac{2^k}{R-r} \norm{Du}{\Leb{s,q}(Q_{r_{k+1}})}
\end{align*}
for some constant $N=N(d,s,q,\nu)>0$.
By the Gagliardo-Nirenberg interpolation inequality, we have
\begin{align*}
     \norm{Du}{\Leb{s,q}(Q_r)}&\leq N(d,q)\norm{D^2 u}{\Leb{s,q}(Q_r)}^{\theta}\norm{u}{\Leb{s,1}(Q_r)}^{1-\theta}\\
     &\relphantom{=}+N(d,q){r^{-d+d/q-1} }\norm{u}{\Leb{s,1}(Q_r)},
\end{align*}
where
\[  \frac{1}{2}< \theta=\frac{1+{1}/{d}-{1}/{q}}{1+{2}/{d}-{1}/{q}}<1.\]
 Then by Young's inequality, we get
\begin{align*}
\norm{D^2 u}{\Leb{s,q}(Q_{r_k})}&\leq N \kappa^{-(d+2)/q_0}\norm{f}{\Leb{s,q}(Q_{r_{k+1}})}+N\kappa^{-(d+2)/q_0}\norm{Dg}{\Leb{s,q}(Q_{r_{k+1}})} \\
&\relphantom{=}+N\left(\kappa^{-(d+2)/q_0}\delta^{1/q_0-1/q_1}+\kappa \right)\norm{D^2 u}{\Leb{s,q}(Q_{r_{k+1}})}\\
&\relphantom{=}+ N\kappa^{-(d+2)/q_0}\left(\frac{2^{k/(1-\theta)}}{(R-r)^{1/(1-\theta)}}+{r^{-d+d/q-1}_{k+1}}\frac{2^k}{R-r} \right)\norm{u}{\Leb{s,1}(Q_{r_{k+1}})},
\end{align*}
where the constant $N$ depends only on $d$ and $q$. {Observe that $r_{k+1}^{-d+d/q-1}\leq r^{-d+d/q-1}$ for all $k$.}

Choose $\kappa$ sufficiently small and then $\delta$ sufficiently small so that
\[
   N\left(\kappa^{-(d+2)/q_0}\delta^{1/q_0-1/q_1}+\kappa \right)\leq 3^{-1/(1-\theta)}.
\]
Then multiply both sides of the inequality by $3^{-k/(1-\theta)}$ and sum over $k=k_0,k_0+1$, ... to obtain
\[
    \norm{D^2 u}{\Leb{s,q}(Q_r)}\leq N\left[\left(\frac{1}{(R-r)}+\frac{1}{(R-r)^b}\right)\norm{u}{\Leb{s,1}(Q_R)}+\norm{f}{\Leb{s,q}(Q_R)}+\norm{Dg}{\Leb{s,q}(Q_R)}\right],
\]
where $N=N(d,s,q,\nu,r,R,R_0)>0$ and $b=1/(1-\theta)$. Since $0<R-r<1/2$, we get the desired estimate. This completes the proof of Theorem \ref{thm:DP}.
\end{proof}

Another ingredient to prove Theorem \ref{thm:B} is the following regularity result for Stokes equations when the exterior force $f$ is bounded and has compact support.

\begin{lemma}\label{lem:regularity-Stokes}
Let $1<q_0,s,q<\infty$ and $0<T<\infty$. There exists $\delta>0$ such that under Assumption \ref{assump:VMO} $(\delta)$, if $(u,p)$ is a strong solution to \eqref{eq:Stokes-nondiv} in $\mathbb{R}^d_T$ with $u(0,\cdot)=0$ on $\mathbb{R}^d$ satisfying
\[
u \in \oSob{1,2}{q_0}(\mathbb{R}^d_T)^d,\quad \nabla p\in \Leb{q_0}(\mathbb{R}^d_T)^d,\quad
\]$f\in \Leb{\infty}(\mathbb{R}^d_T)^d$ having compact support in $\mathbb{R}^d_T${,} and $g=0$, then $(u,\nabla p)\in \Sob{1,2}{s,q}(\mathbb{R}^d_T)^d \times \Leb{s,q}(\mathbb{R}^d_T)^d$.
\end{lemma}
\begin{proof}
Choose $s_1 \in (\max\{s,q_0\},\infty)$ and $q_1 \in (\max\{q,q_0\},\infty)$. Let $q_*=\min\{q,q_0\}$ and define $w(t,x)=w(x)=(1+|x|)^\alpha$, where $d(q_1/q_*-1)<\alpha<d(q_1-1)$. Then by Proposition \ref{prop:weight-property} (iv), $w \in A_{s,q_1}$. By H\"older's inequality, we have
\begin{align*}
\int_{\mathbb{R}^d} |g|^r \myd{x}&\leq \left(\int_{\mathbb{R}^d} |g|^{q_1} (1+|x|)^\alpha\myd{x}\right)^{r/q_1} \left(\int_{\mathbb{R}^d} (1+|x|)^{-\frac{\alpha r}{q_1-r}} \myd{x}\right)^{1-r/q_1}
\end{align*}
for $r\in \{q,q_0\}$. Note that the integral
\[\int_{\mathbb{R}^d} (1+|x|)^{-\frac{\alpha r}{q_1-r}} \myd{x} \]
is finite if and only if $\alpha>d(q_1/r-1)$, which is satisfied by the choice of $q_1$. Since $(0,T)$ has a finite Lebesgue measure, this implies that $\Leb{s_1,q_1,w}(\mathbb{R}^d_T)\subset \Leb{s,q}(\mathbb{R}^d_T)\cap \Leb{q_0}(\mathbb{R}^d_T)$  for our specific weight $w$. Similarly, $\Sob{1,2}{s_1,q_1,w}(\mathbb{R}^d_T)\subset \Sob{1,2}{s,q}(\mathbb{R}^d_T)\cap \Sob{1,2}{q_0}(\mathbb{R}^d_T)$ for our specific weight $w$.

By Theorem \ref{thm:A}, there exist $\delta_1>0$ and strong solutions $(v_1,p_1)$, $(v_2,p_2)$ to \eqref{eq:Stokes-nondiv} in $\mathbb{R}^d_T$ with $v_1(0,\cdot)=v_2(0,\cdot)=0$ on $\mathbb{R}^d$ satisfying $(v_1,\nabla p_1) \in \oSob{1,2}{s_1,q_1,w}(\mathbb{R}^d_T)^d \times \Leb{s_1,q_1,w}(\mathbb{R}^d_T)^d$ and  $(v_2,\nabla p_2) \in \oSob{1,2}{s,q}(\mathbb{R}^d_T)^d \times \Leb{s,q}(\mathbb{R}^d_T)^d$, where $a^{ij}$ satisfies Assumption \ref{assump:VMO} $(\delta_1)$.

Since $\oSob{1,2}{s_1,q_1,w}(\mathbb{R}^d_T)\subset \oSob{1,2}{s,q}(\mathbb{R}^d_T)$, by the uniqueness assertion of Theorem \ref{thm:A}, we conclude that $v_1=v_2$. Choose $0<\delta_2<\delta_1$ so that the uniqueness assertion in $\oSob{1,2}{q_0}(\mathbb{R}^d_T)$ of Theorem \ref{thm:A} holds for $a^{ij}$ satisfying Assumption \ref{assump:VMO} $(\delta_2)$. Since $(v_1,\nabla p_1) \in \oSob{1,2}{q_0}(\mathbb{R}^d_T)^d \times \Leb{q_0}(\mathbb{R}^d_T)^d$ and $(u,p)$ is a strong solution to \eqref{eq:Stokes-nondiv} in $\mathbb{R}^d_T$ satisfying $(u,\nabla p) \in \oSob{1,2}{q_0}(\mathbb{R}^d_T)^d \times \Leb{q_0}(\mathbb{R}^d_T)^d$, it follows from the uniqueness assertion that $v_1=u$. Therefore, $v_1=v_2=u$, which proves that  $u$ belongs to $\Sob{1,2}{s,q}(\mathbb{R}^d_T)^d$. This completes the proof of Lemma \ref{lem:regularity-Stokes}.
\end{proof}

Now we are ready to prove Theorem \ref{thm:B}.
\begin{proof}
By taking mollification in $(t,x)$, we have
\[   \partial_t u^{(\varepsilon)}-a^{ij}(t,x)D_{ij} u^{(\varepsilon)}+\nabla p^{(\varepsilon)}=f^{(\varepsilon)}+h^\varepsilon,\quad \Div u^{(\varepsilon)}=g^{(\varepsilon)}\quad \text{in } Q_{3/4} \]
with $0<\varepsilon<1/4$, where
\[   h^\varepsilon(t,x)=[a^{ij}(t,x)D_{ij}u]^{(\varepsilon)}(t,x)-a^{ij}(t,x)D_{ij}u^{(\varepsilon)}(t,x).\]

By Theorem \ref{thm:A}, there exist $\delta_1>0$ and a unique strong solution $(u_1^\varepsilon,p_1^\varepsilon)$ to
\begin{equation*}%\label{eq:u1-perturb-problem}
    \partial_t u_1^\varepsilon -a^{ij}D_{ij} u_1^\varepsilon + \nabla p_1^\varepsilon = h^\varepsilon 1_{Q_{3/4}},\quad \Div u_1^\varepsilon =0\quad{\text{in $(-1,0)\times\mathbb{R}^d$}}
\end{equation*}
with $u_1^{\varepsilon}(-1,\cdot)=0$ on $\mathbb{R}^d$ satisfying
\[
u_1^\varepsilon \in \oSob{1,2}{q_0}((-1,0)\times \mathbb{R}^d)^d\quad\text{and}\quad \nabla p_1^\varepsilon \in \Leb{q_0}((-1,0)\times \mathbb{R}^d)^d.
\]
Moreover, we have
\begin{equation}\label{eq:u-1-varepsilon}
   \norm{u_1^\varepsilon}{\Sob{1,2}{q_0}((-1,0)\times \mathbb{R}^d)}+ \norm{\nabla p_1^\varepsilon}{\Leb{q_0}((-1,0)\times \mathbb{R}^d)}\leq N \norm{h^\varepsilon}{\Leb{q_0}(Q_{3/4})},
\end{equation}
where $N$ is independent of $\varepsilon$. By Lemma \ref{lem:regularity-Stokes}, there exists $\delta_2>0$ such that under Assumption \ref{assump:VMO} $(\delta_2)$, $u_1^\varepsilon \in \Sob{1,2}{s,q}((-1,0)\times \mathbb{R}^d)^d$. Moreover, if we define $u_2^\varepsilon=u^{(\varepsilon)}-u_1^\varepsilon$ and $p_2^\varepsilon=p^{(\varepsilon)}-p_1^\varepsilon$, then $u_2^{\varepsilon}\in \tilde{W}^{1,2}_{s,q}(Q_{3/4})^d$ and $(u_2^\varepsilon,p_2^\varepsilon)$ is a solution to
\[
   \partial_t u_2^\varepsilon - a^{ij}D_{ij} u_2^\varepsilon + \nabla p_2^\varepsilon = f^{(\varepsilon)},\quad \Div u_2^\varepsilon =g^{(\varepsilon)}\quad \text{in } Q_{3/4}.
\]
Hence by Theorem \ref{thm:DP}, there exists $0<\delta_3<\min\{\delta_1,\delta_2\}$ such that under Assumption \ref{assump:VMO} $(\delta_3)$, we have
\begin{equation}\label{eq:epsilon-inequality}
\begin{aligned}
\norm{D^2 u_2^\varepsilon}{\Leb{s,q}(Q_{1/2})}&\leq N \left(\norm{u_2^\varepsilon}{\Leb{s,1}(Q_{3/4})}+\norm{f^{(\varepsilon)}}{\Leb{s,q}(Q_{3/4})}+\norm{Dg^{(\varepsilon)}}{\Leb{s,q}(Q_{3/4})} \right)\\
&\leq N\left(\norm{u_1^\varepsilon}{\Leb{s,1}(Q_{3/4})}+\norm{u^{(\varepsilon)}}{\Leb{s,1}(Q_{3/4})}\right.\\
&\relphantom{=}\hspace{2em}\left.+\norm{f^{(\varepsilon)}}{\Leb{s,q}(Q_{3/4})}+\norm{Dg^{(\varepsilon)}}{\Leb{s,q}(Q_{3/4})} \right)
\end{aligned}
\end{equation}
for some constant $N=N(d,s,q,\nu,R_0)>0$.

Since $u\in \Leb{s,1}(Q_1)^d$, $f\in \Leb{s,q}(Q_1)^d$, and $g\in \Sob{0,1}{s,q}(Q_1)$, we have $u^{(\varepsilon)}\rightarrow u$ in $\Leb{s,1}(Q_{3/4})$, $f^{(\varepsilon)}\rightarrow f$ in $\Leb{s,q}(Q_{3/4})$, and $Dg^{(\varepsilon)}\rightarrow Dg$ in $\Leb{s,q}(Q_{3/4})$.

Note that $h^\varepsilon \rightarrow 0$ in $\Leb{q_0}(Q_{3/4})$ as $\varepsilon \rightarrow 0$. Hence by \eqref{eq:u-1-varepsilon} and Sobolev embedding theorem, we have $\norm{u_1^\varepsilon}{\Leb{s,1}(Q_{3/4})}\rightarrow 0$ as $\varepsilon \rightarrow 0$. This implies that there exists a constant $N$ independent of $\varepsilon$ such that
\[
   \sup_{\varepsilon>0} \norm{D^2 u_2^\varepsilon}{\Leb{s,q}(Q_{1/2})}\leq N.
\]
Hence by the weak compactness in $\Leb{s,q}(Q_{1/2})$, there exists a subsequence $\{D^2 u_2^{\varepsilon_j}\}$ of $\{D^2 u_2^\varepsilon\}$ which converges weakly to a function $v$ in $\Leb{s,q}(Q_{1/2})$.

On the other hand, since $D^2 u^{(\varepsilon)}\rightarrow D^2 u$ strongly in $\Leb{q_0}(Q_{3/4})$ and $D^2 u_1^\varepsilon \rightarrow 0$ strongly in $\Leb{q_0}(Q_{3/4})$ by \eqref{eq:u-1-varepsilon} as $\varepsilon \rightarrow 0^+$, it follows that $D^2 u_2^\varepsilon\rightarrow D^2 u$ strongly in $\Leb{q_0}(Q_{3/4})$. Hence it follows that $D^2 u=v$ in $Q_{1/2}$. Therefore, by taking liminf in \eqref{eq:epsilon-inequality}, we get the desired result. This completes the proof of Theorem \ref{thm:B}.
\end{proof}

\subsection{Interior gradient estimates for Stokes equations in divergence form}\label{subsec:gradient}

The following theorem is an analog of Theorem \ref{thm:DP}, which was implicitly proved in \cite[Theorem 1.9]{DP21}.

\begin{theorem}\label{thm:DP-2}
Let $1<s,q<\infty$, $\nu \in (0,1)$, and $1/2\leq r<R\leq 1$. There exists $\delta=\delta(d,\nu,s,q)\in (0,1)$ such that under Assumption \ref{assump:VMO} $(\delta_1)$, if $u\in \Sob{0,1}{s,q}(Q_R)^d$ is a weak solution to \eqref{eq:Stokes-nondiv} in $Q_R$, $\boldF\in \Leb{s,q}(Q_R)^{d\times d}$ and $g\in \Leb{s,q}(Q_R)$, then there exists a constant $N=N(d,s,q,\nu,r,R,R_0)>0$ such that
\begin{align*}
\norm{D u}{\Leb{s,q}(Q_{r})}\leq N\left[\frac{1}{(R-r)^b}\norm{u}{\Leb{s,1}(Q_R)}+\norm{\boldF}{\Leb{s,q}(Q_R)}+\norm{g}{\Leb{s,q}(Q_R)} \right]
\end{align*}
for some $b=b(d,q)>2$.
\end{theorem}

We omit the proof of Theorem \ref{thm:DP-2} since it is almost identical to that of Theorem \ref{thm:DP} by using \cite[Lemma 3.11]{DP21} stated below.

\begin{lemma}
Let $1/2\leq R< 1$, $R_1\in (0,R_0)$, $R_1 \in (0,R_0)$, $\delta \in (0,1)$, $\kappa \in (0,1/4)$, $1<s,q<\infty$, $q_1 \in (1,\min\{s,q\})$, and $1<q_0<q_1$. Suppose that $u\in \Sob{0,1}{s,q}(Q_{R+R_1})^d$  is a weak solution to \eqref{eq:Stokes-div} in $Q_{R+R_1}$ {for some $\boldF \in \Leb{s,q}(Q_{R+R_1})^{d\times d}$ and $g\in \Leb{s,q}(Q_{R+R_1})$}. Then
\begin{align*}
\norm{D u}{\Leb{s,q}(Q_R)}&\leq N \kappa^{-(d+2)/q_0}\norm{\boldF}{\Leb{s,q}(Q_{R+R_1})}+N\kappa^{-(d+2)/q_0}\norm{g}{\Leb{s,q}(Q_{R+R_1})} \\
&\relphantom{=}+N\left(\kappa^{-(d+2)/q_0}\delta^{1/q_0-1/q_1}+\kappa \right)\norm{D u}{\Leb{s,q}(Q_{R+R_1})}\\
&\relphantom{=}+ N\kappa^{-{(d+2)}/{q_0}}R_1^{-1} \norm{u}{\Leb{s,q}(Q_{R+R_1})}.
\end{align*}
\end{lemma}

\begin{remark}
The conditions $p\in \Leb{1}(Q_{R+R_1})$ and $u_t \in \mathbb{H}^{-1}_1(Q_{R+R_1})$ are not essentially used in the proof in \cite[Lemma 3.11]{DP21}. It suffices to derive a vorticity equation from Stokes equations with simple coefficients. For simplicity, we assume that $\boldF=0$. For $k,l=1,\dots,d$ and $\psi \in C_0^\infty (Q_1)$, define $\phi=(D_k \psi)e_l-(D_l \psi)e_k$. Then it is easy to see that $\Div \phi(t,\cdot)=0$ in $B_1$ for $t\in (-1,0)$. For $u=(u^1,\dots,u^d)$,  define $\omega=\nabla\times u$. If we use $\phi$ as a test function in the definition of weak solutions, then it is easy to check that $\omega_{kl}$ is a very weak solution of
\[  \partial_t \omega_{kl}-D_i(a^{ij}D_{j} \omega_{kl})=0\quad \text{in } Q_1. \]
\end{remark}

Another ingredient for proving Theorem \ref{thm:D} is the following regularity lemma similar to Lemma \ref{lem:regularity-Stokes}, which can be proved by using Theorem \ref{thm:C} instead of Theorem \ref{thm:A} of which proof is omitted.

\begin{lemma}\label{lem:regularity-Stokes-div}
Let $1<q_0,s,q<\infty$. There exists $\delta>0$ such that under Assumption \ref{assump:VMO} $(\delta)$, if $(u,p)\in \odSob{1}{q_0}(\mathbb{R}^d_T)^d\times \Leb{q_0}(\mathbb{R}^d_T)$ is a weak solution of the problem \eqref{eq:Stokes-div} with $u(0,\cdot)=0$ on $\mathbb{R}^d$, $\boldF \in \Leb{\infty}(\mathbb{R}^d_T)^{d\times d}$ having compact support in $\mathbb{R}^d_T$ and $g=0$, then $u\in \dSob{1}{s,q}(\mathbb{R}^d_T)^d$.
\end{lemma}

Now we are ready to prove Theorem \ref{thm:D}.
\begin{proof}
We may assume that $u\in \Leb{s,1}(Q_1)^d$ {because otherwise the desired inequality is trivial}.
Suppose first that $s>q_0$. Choose $1<s_0,s_1,s_2,q_1,q_2<\infty$ so that
\[
     \frac{1}{s_k} \leq \frac{1}{2}+\frac{1}{s_{k+1}},\quad q_0=s_0<s_1\leq  s_2=s,\quad \text{and}\,\,
    q_1=q_2=q.
\]

We first show that $D u\in \Leb{s_1,q}(Q_{3/4})$. {Fix $\psi \in C_0^\infty(Q_{7/8})^d$ with $\Div \psi(t,\cdot)=0$  in $B_{7/8}$ for $t\in (-(7/8)^2,0)$ and $0<\varepsilon<1/8$. Then if we use $\phi=\tilde{\psi}^{(\varepsilon)}$ as a test function in the definition of weak solutions, where $\tilde{\psi}^{(\varepsilon)}$ is defined in \eqref{eq:forward-mollification}, then one can check that } $u^{(\varepsilon)}$ is a weak solution of
\[
  \partial_t u^{(\varepsilon)} -D_i (a^{ij}D_j u^{(\varepsilon)})+\nabla p^{(\varepsilon)}=\Div( \boldF^{(\varepsilon)} + \boldH^{\varepsilon})\quad \text{in } Q_{7/8}
\]
and
\[    \Div u^{(\varepsilon)}= g^{(\varepsilon)}\quad \text{in } Q_{7/8}, \]
where $\boldH^{\varepsilon}=(H_1^\varepsilon,\dots,H_d^\varepsilon)$ and $H_i^\varepsilon=(a^{ij}D_j u)^{(\varepsilon)}-a^{ij}D_ju^{(\varepsilon)}$. Then by Theorem \ref{thm:C}, there exist $\delta_1>0$ and a unique $(u_1^\varepsilon,p_1^\varepsilon) \in \odSob{1}{q_0}((-1,0)\times\mathbb{R}^d)^d \times \Leb{q_0}((-1,0)\times\mathbb{R}^d)$ satisfying
\begin{equation}\label{eq:divergence-form-perturbation}
\partial_t u_1^\varepsilon-D_i (a^{ij}D_j u_1^\varepsilon)+\nabla p_1^\varepsilon=\Div (\boldH^\varepsilon 1_{Q_{7/8}}),\quad \Div u_1^\varepsilon=0\quad {\text{in } (-1,0)\times \mathbb{R}^d}
\end{equation}
with  $u_1^\varepsilon(-1,\cdot)=0$ on $\mathbb{R}^d$. Moreover, we have
\begin{equation}\label{eq:u1-perturb-3}
\norm{u_1^\varepsilon}{\dSob{1}{q_0}((-1,0)\times\mathbb{R}^d)}+\norm{p_1^\varepsilon}{\Leb{q_0}((-1,0)\times \mathbb{R}^d)}\leq N \norm{\boldH^\varepsilon}{\Leb{q_0}(Q_{7/8})},
\end{equation}
where the constant $N$ is independent of $\varepsilon$.

Define $u_2^\varepsilon=u^{(\varepsilon)}-u_1^\varepsilon$. By Lemma \ref{lem:regularity-Stokes-div}, there exists $\delta_2>0$ such that under Assumption \ref{assump:VMO} $(\delta_2)$, $u_2^\varepsilon \in {\mathcal{H}}^1_{s_{1},q}(Q_{7/8})^d$ is a weak solution to
\[
\partial_t u_2^\varepsilon-D_i(a^{ij}D_j u^\varepsilon_2)+\nabla p_2^\varepsilon=\Div \boldF^{(\varepsilon)}\quad\text{and}\quad \Div u_2^\varepsilon=g^{(\varepsilon)}\quad \text{in } Q_{7/8}.
\]
By Theorem \ref{thm:DP-2}, there exists $0<\delta_3<\min\{\delta_1,\delta_2\}$ such that under Assumption \ref{assump:VMO} ($\delta_3$), we have
\begin{align*}
\norm{Du_2^\varepsilon}{\Leb{s_{1},q}(Q_{3/4})}&\leq N\left(\norm{u_2^{\varepsilon}}{\Leb{s_{1},1}(Q_{7/8})}+\norm{\boldF^{(\varepsilon)}}{\Leb{s_{1},q}(Q_{7/8})}+\norm{g^{(\varepsilon)}}{\Leb{s_{1},q}(Q_{7/8})}\right)\\
&\leq N \left(\norm{u^{(\varepsilon)}}{\Leb{s_{1},1}(Q_{7/8})}+\norm{u_1^{\varepsilon}}{\Leb{s_{1},1}(Q_{7/8})} \right.\\
&\relphantom{=}\quad \left.+\norm{\boldF^{(\varepsilon)}}{\Leb{s_{1},q}(Q_{7/8})}+\norm{g^{(\varepsilon)}}{\Leb{s_{1},q}(Q_{7/8})}\right)
\end{align*}
for some constant $N=N(d,s,q,\nu,R_0)>0$. Since $u \in \Leb{s,1}(Q_1)^d$, we see that $u^{(\varepsilon)}\rightarrow u$ in $\Leb{s_1,1}(Q_{7/8})$, $\boldF^{(\varepsilon)}\rightarrow \boldF$, and $g^{(\varepsilon)}\rightarrow g$ in $\Leb{s_1,q}(Q_{7/8})$. Since $Du\in \Leb{q_0}(Q_1)$, it follows that $\boldH^\varepsilon \rightarrow 0$ in $\Leb{q_0}(Q_{7/8})$ as $\varepsilon\rightarrow 0$.  By \eqref{eq:u1-perturb-3} and Lemma \ref{lem:KRW}, we have
{$\norm{Du_1^{\varepsilon}}{\Leb{q_{0}}(Q_{7/8})}\rightarrow 0$ and}
$\norm{u_1^{\varepsilon}}{\Leb{s_{1},q_{0}}(Q_{7/8})}\rightarrow 0$ as $\varepsilon\rightarrow 0+$, and hence $\norm{u_1^\varepsilon}{\Leb{s_1,1}(Q_{7/8})}\rightarrow 0$ as $\varepsilon\rightarrow 0$. This implies that
\[
\sup_{\varepsilon>0} \norm{D u_2^\varepsilon}{\Leb{s_{1},q}(Q_{3/4})}\leq N
\]
for some constant $N>0$. Hence by the weak compactness in $\Leb{s_{1},q}(Q_{3/4})$, there  exists a convergent subsequence $\{D u_2^{\varepsilon_j}\}$ of $\{D u_2^{\varepsilon}\}$  which converges weakly to a function $v$ in $\Leb{s_{1},q}(Q_{3/4})$.

On the other hand, since $D u^{(\varepsilon)}\rightarrow D u$ strongly in $\Leb{q_0}(Q_{7/8})$ and $D u_1^\varepsilon \rightarrow 0$ in $\Leb{q_0}(Q_{7/8})$ by \eqref{eq:u1-perturb-3} as $\varepsilon\rightarrow 0^+$, it follows that $D u_2^\varepsilon\rightarrow D u$ strongly in $\Leb{q_0}(Q_{7/8})$. Hence we conclude that $v=D u$ in $Q_{3/4}$. Therefore, under Assumption \ref{assump:VMO} $(\delta_3)$, $Du \in \Leb{s_1,q}(Q_{3/4})$ and
\begin{align*}
\norm{D u}{\Leb{s_{1},q}(Q_{3/4})}&\leq N \left(\norm{u}{\Leb{s_{1},1}(Q_{3/4})}+\norm{\boldF}{\Leb{s_{1},q}(Q_{7/8})}+\norm{g}{\Leb{s_{1},q}(Q_{7/8})} \right)
\end{align*}
for some constant $N=N(d,s,q,\nu,R_0)>0$.

Since $Du \in \Leb{s_1,q}(Q_{3/4})$, we see that $\boldH^\varepsilon\rightarrow 0$ in $\Leb{s_1,q}(Q_{5/8})$ as $\varepsilon\rightarrow 0$. Then by Theorem \ref{thm:C}, there exists $0<\delta_4<\delta_3$ and a unique $(u_1^\varepsilon,p_1^\varepsilon) \in \odSob{1}{s_1,q}((-1,0)\times\mathbb{R}^d)^d\times \Leb{s_1,q}((-1,0)\times\mathbb{R}^d)$ satisfying
\[
\partial_t u_1^\varepsilon-D_i (a^{ij}D_j u_1^\varepsilon)+\nabla p_1^\varepsilon=\Div (\boldH^\varepsilon 1_{Q_{5/8}}),\quad \Div u_1^\varepsilon=0,\quad u_1^\varepsilon(-1,\cdot)=0\quad\text{on } \mathbb{R}^d,
\]
where $a^{ij}$ satisfies Assumption \ref{assump:VMO} ($\delta_4$).
Moreover, we have
\begin{equation}\label{eq:u1-perturb-4}
\norm{u_1^\varepsilon}{\dSob{1}{s_1,q}((-1,0)\times\mathbb{R}^d)}+\norm{p_1^\varepsilon}{\Leb{s_1,q}((-1,0)\times \mathbb{R}^d)}\leq N \norm{\boldH^\varepsilon}{\Leb{s_1,q}(Q_{5/8})},
\end{equation}
where the constant $N$ is independent of $\varepsilon$. By \eqref{eq:u1-perturb-4} and Lemma \ref{lem:KRW}, $\norm{u_1^\varepsilon}{\Leb{s,q}(Q_{5/8})}\rightarrow 0$ as $\varepsilon\rightarrow 0$. Hence following the above compactness argument, we can show that $Du \in \Leb{s,q}(Q_{1/2})$. Moreover, we have
\begin{align*}
\norm{D u}{\Leb{s,q}(Q_{1/2})}&\leq N \left(\norm{u}{\Leb{s,1}(Q_{5/8})}+\norm{\boldF}{\Leb{s,q}(Q_{5/8})}+\norm{g}{\Leb{s,q}(Q_{5/8})} \right)
\end{align*}
for some constant $N=N(d,s,q,\nu,R_0)>0$. {In a similar way, we can also prove the case $s\leq q_0$.}  This completes the proof of Theorem \ref{thm:D}.
\end{proof}

\begin{remark}\label{rem:t-dependent}
(i) If the viscosity coefficient $a^{ij}$ depends only on $t$, then we can show that Theorem \ref{thm:D} holds if $u\in \Leb{s,1}(Q_1)^d$ is a very weak solution to \eqref{eq:Stokes-div} in $Q_1$ for some $\boldF \in \Leb{s,q}(Q_1)^{d}$ and $g\in \Leb{s,q}(Q_1)$. We say that $u \in \Leb{s,1}(Q_1)^d$ is a \emph{very weak solution} to \eqref{eq:Stokes-div} in $Q_1$ if
\begin{equation*}%\label{eq:very-weak}
\int_{Q_1} u\cdot(\partial_t \phi + a^{ij}D_{ij} \phi)\myd{x}dt=-\int_{Q_1} \boldF:\nabla\phi \myd{x}dt
\end{equation*}
for all $\phi \in C_0^\infty(Q_1)^d$ with $\Div \phi(t)=0$ in $B_1$ for all $t\in (-1,0)$, and
\begin{equation*}
-\int_{B_1} u \cdot \nabla \varphi \myd{x}=\int_{B_1} g \varphi \myd{x}
\end{equation*}
for all $\varphi \in C_0^\infty(B_1)$ and a.e. $t\in (-1,0)$.

Let $\phi \in C_0^\infty(\mathbb{R})$ and $\zeta \in C_0^\infty(B_1)^d$, where $\phi=0$ if $t\geq 0$, $\int_{-1}^0 \phi \,dt=1$, $\Div \zeta =0$ in $B_1$, and $\int_{B_1} \zeta \myd{x}=1$. Define $\phi_\eta(t)=\eta^{-2}\varphi(t/\eta^2)$ and $\zeta_\varepsilon(x)=\varepsilon^{-d}\zeta(x/\varepsilon)$. For $(t,x) \in (-1+\eta^2,0)\times B_{1-\varepsilon}$, define
\begin{align*}
    u^{(\eta,\varepsilon)}(t,x)=(u^{(\varepsilon)})^{(\eta)}(t,x)&=\int_{-\eta^2}^0\int_{B_\varepsilon} u(t+s,x+y) \varphi_{\eta}(s)\zeta_\varepsilon(y)\myd{y}ds\\
    &=\int_{Q_1} u(s,y) \phi_{\eta}(s-t)\zeta_{\varepsilon}(y-x)\myd{s}dy.
\end{align*}
Then it is easy to verify that for small $\varepsilon,\eta>0$, $u^{(\eta,\varepsilon)}$ is a weak solution to
\[   \partial_t u^{(\eta,\varepsilon)}-D_i(a^{ij}(t)D_{j} u^{(\eta,\varepsilon)})+\nabla p^{(\eta,\varepsilon)}=\Div (\boldF^{(\eta,\varepsilon)}+\boldH^{\eta,\varepsilon})\quad \text{in } Q_{3/4}, \]
where $\boldH^{\eta,\varepsilon}=(H_1^{\eta,\varepsilon},\dots,H_d^{\eta,\varepsilon})$,
\[
H^{\eta,\varepsilon}_i(t,x)=(a^{ij}D_{j} u^{(\varepsilon)})^{(\eta)}-a^{ij}D_{j} u^{(\eta,\varepsilon)},\quad i=1,\dots,d,\]
and
 \[   \Div u^{(\eta,\varepsilon)}=g^{(\eta,\varepsilon)}\quad\text{in } Q_{3/4}. \]
Following the argument as in the proof of Theorem \ref{thm:D}, we can prove the desired result. We give a sketch of the proof.

Since $u \in \Leb{s,1}(Q_1)^d$,  $\boldH^{\eta,\varepsilon} \in \Leb{s,q}(Q_1)^{d\times d}$. By Theorem \ref{thm:mixed-simple-div}, there exists $(u_1^{\eta,\varepsilon},p_1^{\eta,\varepsilon})\in \odSob{1}{s,q}((-(3/4)^2,0)\times\mathbb{R}^d)^d\times \Leb{s,q}((-(3/4)^2,0)\times\mathbb{R}^d)$ satisfying \eqref{eq:divergence-form-perturbation} with $u_1^{\eta,\varepsilon}(-(3/4)^2,\cdot)=0$ on $\mathbb{R}^d$, where  $H^\varepsilon$ is replaced with $H^{\eta,\varepsilon}$. Moreover, we have
\[
    \norm{u_1^{\eta,\varepsilon}}{\dSob{1}{s,q}((-(3/4)^2,0)\times\mathbb{R}^d)}\leq N \norm{H^{\eta,\varepsilon}}{\Leb{s,q}(Q_{3/4})}
\]
for some constant $N$ independent of $\eta$. Define $u_2^{\eta,\varepsilon}=u^{(\eta,\varepsilon)}-u_1^{\eta,\varepsilon}$. Then following the argument as in the proof of Theorem \ref{thm:D}, we have
  \begin{align*}
 \norm{D u_2^{\eta,\varepsilon}}{\Leb{s,q}(Q_{1/2})}&\leq N \left(\norm{u_2^{\eta,\varepsilon}}{\Leb{s,1}(Q_{3/4})}+\norm{\boldF^{(\eta,\varepsilon)}}{\Leb{s,q}(Q_{3/4})}+\norm{g^{(\eta,\varepsilon)}}{\Leb{s,q}(Q_{3/4})}\right)\\
 &\leq N \left(\norm{u_1^{\eta,\varepsilon}}{\Leb{s,1}(Q_{3/4})}+\norm{u^{(\eta,\varepsilon)}}{\Leb{s,1}(Q_{3/4})}\right.\\
&\relphantom{=}\quad\left.+\norm{\boldF^{(\eta,\varepsilon)}}{\Leb{s,q}(Q_{3/4})}+\norm{g^{(\eta,\varepsilon)}}{\Leb{s,q}(Q_{3/4})}\right)
 \end{align*}
 for some constant $N=N(d,s,q,\nu)>0$.  Note that for fixed $\varepsilon>0$, $H_i^{\eta,\varepsilon} \rightarrow 0$ in $\Leb{s,q}(Q_{3/4})$ as $\eta \rightarrow 0$. Hence it follows that
$\norm{u_1^{\eta,\varepsilon}}{\Leb{s,1}(Q_{3/4})}\rightarrow 0$ as $\eta \rightarrow 0$ and
\[
    \sup_{\eta} \norm{D u^{\eta,\varepsilon}_2}{\Leb{s,q}(Q_{1/2})}\leq N(\varepsilon),
\]
where $N$ is independent of $\eta$.

Note that $a^{ij}D_{j} u^{(\eta,\varepsilon)}\rightarrow a^{ij}D_j u^{(\varepsilon)}$ {a.e} as $\eta\rightarrow0$. Also, it follows that
\[
    |D_j u^{(\eta,\varepsilon)}(t,x)-D_j u^{(\varepsilon)}(t,x)| \leq N \mathcal{M}^t (D_j u^{(\varepsilon)})(t,x)
\]
{for some constant $N=N(d)>0$,} where $\mathcal{M}^t$ denotes the one-dimensional maximal function in $t$. Hence by the Hardy-Littlewood maximal function theorem and the dominated convergence theorem, we can show that $\boldH^{\eta,\varepsilon}\rightarrow 0$ in $\Leb{s,q}(Q_{3/4})$ as $\eta\rightarrow 0$ for fixed $\varepsilon>0$. By a  compactness argument as in the proof of Theorem \ref{thm:B}, we get
 \[
 \norm{Du^{(\varepsilon)}}{\Leb{s,q}(Q_{1/2})}\leq N\left(\norm{u^{(\varepsilon)}}{\Leb{s,1}(Q_{3/4})}+\norm{\boldF^{(\varepsilon)}}{\Leb{s,q}(Q_{3/4})}+\norm{g^{(\varepsilon)}}{\Leb{s,q}(Q_{3/4})} \right)
 \]
 for some constant $N=N(d,s,q,\nu)>0$. Since $u^{(\varepsilon)}\rightarrow u$ in $\Leb{s,1}(Q_{3/4})$, $\boldF^{(\varepsilon)}\rightarrow \boldF$, and $g^{(\varepsilon)}\rightarrow g$ in $\Leb{s,q}(Q_{3/4})$, it follows that
\[
\sup_{\varepsilon>0}   \norm{Du^{(\varepsilon)}}{\Leb{s,q}(Q_{1/2})} \leq N.
\]
Hence by a previous compactness argument, it follows that $D u$ exists in $Q_{1/2}$ and is in $\Leb{s,q}(Q_{1/2})$. Moreover, we have
 \[
 \norm{Du}{\Leb{s,q}(Q_{1/2})}\leq N\left(\norm{u}{\Leb{s,1}(Q_{3/4})}+\norm{\boldF}{\Leb{s,q}(Q_{3/4})}+\norm{g}{\Leb{s,q}(Q_{3/4})} \right)
 \]
for some constant $N=N(d,s,q,\nu)>0$.

(ii) Similarly, if $u\in \Leb{s,1}(Q_1)^d$ satisfies
\begin{equation*}%\label{eq:very-weak-2}
-\int_{Q_1} u\cdot(\partial_t \phi + a^{ij}D_{ij} \phi)\myd{x}dt=-\int_{Q_1} f\cdot \phi \myd{x}dt
\end{equation*}
for all $\phi \in C_0^\infty(Q_1)^d$ with $\Div \phi(t)=0$ for $t\in (-1,0)$ and
\begin{equation*}
-\int_{B_1} u \cdot \nabla \varphi \myd{x}=\int_{B_1} g \varphi \myd{x}
\end{equation*}
for all $\varphi \in C_0^\infty(B_1)$ and a.e. $t\in (-1,0)$ with $f\in\Leb{s,q}(Q_1)^{d}$ and $g\in \Sob{0,1}{s,q}(Q_{1})$, then $D^2 u \in \Leb{s,q}(Q_{1/2})$ and
\[
\norm{D^2 u}{\Leb{s,q}(Q_{1/2})}\leq N\left(\norm{u}{\Leb{s,1}(Q_1)}+\norm{f}{\Leb{s,q}(Q_1)}+\norm{g}{\Sob{0,1}{s,q}(Q_1)} \right)
 \]
 for some constant $N=N(d,s,q,\nu)>0$.

Indeed, by well-known mixed norm solvability results of the heat equation in a bounded cylindrical domain (see e.g. \cite{DK18-weight}), there exists $v \in \Sob{1,2}{s,q}(Q_1)^d$ such that $\partial_t v -\Delta v=f$ in $Q_1$ and  $v=0$ on $\partial_p Q_1$. Moreover, we have
\begin{equation}\label{eq:heat-estimate-very-weak}
\norm{v}{\Sob{1,2}{s,q}(Q_1)}\leq N \norm{f}{\Leb{s,q}(Q_1)}
\end{equation}
 for some constant $N=N(d,s,q)>0$. Define $w=u-v$. Then it is easy to show that $w$ is a very weak solution to
\[
   w_t-a^{ij}(t)D_{ij}w+\nabla p =D_i((a^{ij}-\delta^{ij})D_j v)\quad\text{and}\quad \Div w = g-\Div v\quad \text{in } Q_1.
\]
Hence it follows from the previous result and \eqref{eq:heat-estimate-very-weak} that $Dw\in \Leb{s,q}(Q_{3/4})$ and
\begin{equation}\label{eq:gradient-w-estimate-very-weak}
\begin{aligned}
    \norm{Dw}{\Leb{s,q}(Q_{3/4})}&\leq N \left(\norm{w}{\Leb{s,1}(Q_1)}+\norm{Dv}{\Leb{s,q}(Q_1)} +\norm{g}{\Leb{s,q}(Q_1)}\right)\\
    &\leq N  \left(\norm{w}{\Leb{s,1}(Q_1)}+\norm{f}{\Leb{s,q}(Q_1)}+\norm{g}{\Leb{s,q}(Q_1)} \right)
\end{aligned}
\end{equation}
for some constant $N=N(d,s,q,\nu)>0$.
Since $u=v+w$ and $Dv \in\Leb{s,q}(Q_{3/4})$, we have $Du \in \Leb{s,q}(Q_{3/4})$. Moreover, it follows from \eqref{eq:heat-estimate-very-weak} and \eqref{eq:gradient-w-estimate-very-weak} that
\begin{align*}
    \norm{Du}{\Leb{s,q}(Q_{3/4})}&\leq N \left(\norm{w}{\Leb{s,1}(Q_1)}+\norm{v}{\Leb{s,1}(Q_1)}+\norm{f}{\Leb{s,q}(Q_1)}+\norm{g}{\Leb{s,q}(Q_1)} \right)\\
    &\leq N  \left(\norm{u}{\Leb{s,1}(Q_1)}+\norm{f}{\Leb{s,q}(Q_1)}+\norm{g}{\Leb{s,q}(Q_1)} \right)
\end{align*}
for some constant $N=N(d,s,q,\nu)>0$.

For $1\leq k\leq d$, observe that $D_k u$ is a very weak solution to
\[
\partial_t v -D_i(a^{ij}D_j v)+\nabla p = \Div \boldF,\quad \Div v= D_k g\quad \text{in } Q_{3/4}, \]
where $\boldF^{ij}=f^i \delta_{jk}$. Hence it follows from the previous result that $D(D_k u) \in \Leb{s,q}(Q_{1/2})$ and
\begin{equation*}%\label{eq:hessian-very-weak}
\begin{aligned}
   \norm{D^2 u}{\Leb{s,q}(Q_{1/2})}&\leq N \left(\norm{D_k u}{\Leb{s,q}(Q_{3/4})}+\norm{f}{\Leb{s,q}(Q_{3/4})}+\norm{D_k g}{\Leb{s,q}(Q_{3/4})} \right)
\end{aligned}
\end{equation*}
for some constant $N=N(d,s,q,\nu)>0$. Then by using interpolation inequality on $Du$ and a standard iteration argument as in Lemma \ref{lem:DP},  one can prove that
\[
    \norm{D^2u}{\Leb{s,q}(Q_{1/2})}\leq N \left(\norm{u}{\Leb{s,1}(Q_{3/4})}+\norm{f}{\Leb{s,q}(Q_{3/4})}+\norm{Dg}{\Leb{s,q}(Q_{3/4})} \right)
\]
for some constant $N=N(d,s,q,\nu)>0$. We omit the details.
\end{remark}

\section{Boundary mixed-norm Hessian estimates for Stokes equations}\label{sec:8}

In this section, we briefly sketch how to obtain boundary mixed-norm Hessian estimates under the Lions boundary conditions. The details of this proof are omitted for the sake of brevity, but essentially only involve the same procedures in Sections \ref{sec:4}, \ref{sec:5}, and \ref{sec:7}.  As usual, we may assume that $a^{ij}$ is symmetric.

We first obtain a weighted mixed-norm estimates for Stokes equations in nondivergence form with simple coefficients, a weighted version of Dong-Kim-Phan \cite[Theorem 1.4]{DKP22}. This result can be obtained by following an argument in Theorem \ref{thm:mixed-simple} and an extension argument given in Dong-Kim-Phan \cite[Theorem 1.4]{DKP22}. Then we obtain mean oscillation estimates for $D\omega$ similar to Lemma \ref{lem:mixed-norm} using the H\"older estimate for $D\omega$ which was proved in \cite[Lemma 3.2]{DKP22} and following an argument in Dong-Kim  \cite[Lemma 5.13]{DK18-weight} and \cite[Lemmas 5.1 and 5.3]{DKP22}. Then under Assumption \ref{assump:VMO} $(\delta)$, we can prove weighted mixed-norm solvability results for Stokes equations in nondivergence form under the Lions boundary conditions following the proof of Theorem \ref{thm:A} and the method of continuity.

In summary, we have the following theorem.
\begin{theorem}\label{thm:boundary-solvability}
Let $1<s,q<\infty$, $0<T<\infty$, and let $K_0\geq 1$ be constant, $w=w_1(x)w_2(t)$, where $[w_1]_{A_q(\mathbb{R}^d,dx)}\leq K_0$ and $[w_2]_{A_s(\mathbb{R},dt)}\leq K_0$. 
 There exists $0<\delta<1$ depending only on $d$, $\nu$, $s$, $q$,  and $K_0$ such that under Assumption \ref{assump:VMO} $(\delta)$, for every $f\in \Leb{s,q,w}((0,T)\times\mathbb{R}^d_+)^d$, $g\in\odSob{1}{s,q,w}((0,T)\times\mathbb{R}^d_+)$, and $g_t=\Div G$ for some vector field $G =(G_1,\dots,G_d)\in \Leb{s,q,w}((0,T)\times\mathbb{R}^d_+)^d$ in the sense that
\begin{equation*}
\int_{(0,T)\times\mathbb{R}^d_+} g\varphi_t \myd{x}dt=\int_{(0,T)\times\mathbb{R}^d_+} G \cdot \nabla \varphi \myd{x}dt
\end{equation*}
for any $\varphi \in C_0^\infty([0,T)\times\mathbb{R}^d_+)$, there exists a unique strong solution $(u,p)$ to \eqref{eq:Stokes-nondiv} in $(0,T)\times\mathbb{R}^d_+$ with $u(0,\cdot)=0$ on $\mathbb{R}^d_+$ satisfying
\[ u\in\oSob{1,2}{s,q,w}((0,T)\times\mathbb{R}^d_+)^d,\quad \nabla p \in\Leb{s,q,w}((0,T)\times\mathbb{R}^d_+)^d,\]
and
\[   D_d u^k=u^d=0\quad \text{on } [0,T)\times \mathbb{R}^{d-1}\times\{0\},\quad k=1,2,\dots,d-1.\]
 Moreover, we have
\begin{align*}
  &\norm{u}{\Sob{1,2}{s,q,w}((0,T)\times\mathbb{R}^d_+)}+\norm{\nabla p}{\Leb{s,q,w}((0,T)\times\mathbb{R}^d_+)}\\
  &\leq N \left(\norm{f}{\Leb{s,q}((0,T)\times\mathbb{R}^d_+)}+\norm{Dg}{\Leb{s,q,w}((0,T)\times\mathbb{R}^d_+)}+\norm{G}{\Leb{s,q,w}((0,T)\times\mathbb{R}^d_+)}\right),
\end{align*}
where $N=N(d,s,q,K_0,\nu,R_0,T)>0$.
\end{theorem}

To prove Theorem \ref{thm:boundary-estimate}, let $\tilde{u}^k$ be the even extensions of $u^k$ with respect to $x_d$, $k=1,\dots,d-1$, $\tilde{u}^d$ be the odd extensions of $u^d$ with respect to $x_d$.  Let $\tilde{f}^k(t,\cdot)$ be the even extension of $f^k(t,\cdot)$ for $k=1,\dots,d-1$, and $\tilde{f}^d(t,\cdot)$ be the odd extension of $f^d(t,\cdot)$. Similarly, let $\tilde{g}(t,\cdot)$ be the even extension of $g(t,\cdot)$ with respect to $x_d$. Let $\tilde{p}$ be the even extension of $p$ in $x_d$. By \eqref{eq:Lions-boundary}, $\tilde{u} \in \Sob{1,2}{q_0}(Q_1)^d$. Also, it is easy to verify that $\tilde{p} \in \Sob{0,1}{1}(Q_1)$, $\tilde{f} \in \Leb{s,q}(Q_1)$, and $\tilde{g}\in \tilde{W}^{1,2}_{s,q}(Q_1)$. Moreover, $\tilde{u}|_{Q_1^+} = u$, $\tilde{p}|_{Q_1^+}=p$, $\tilde{f}|_{Q_1^+}=f$, and $\tilde{g}|_{Q_1^+}=g$.

Define $\overline{a}^{ij}(t,x',x_d)=a^{ij}(t,x',x_d)$ if $x_d>0$. For $x_d<0$, define $\overline{a}^{ij}(t,x',x_d)$ to be $a^{ij}(t,x',-x_d)$ if $i,j=1,\dots,d-1$, $\overline{a}^{id}(t,x',x_d)=\overline{a}^{di}(t,x',x_d):=-a^{id}(t,x',-x_d)$ for $i=1,\dots,d-1$. Finally, we define $\overline{a}^{dd}(t,x)=a^{dd}(t,x',-x_d)$. By a direct computation, $(\tilde{u},\tilde{p})$ satisfies
\[
\partial_t \tilde{u}-\overline{a}^{ij}D_{ij} \tilde{u}+\nabla \tilde{p}=\tilde{f},\quad \Div\tilde{u}=\tilde{g}\quad \text{in } Q_1,\]
and
\[    D_d \tilde{u}^k=\tilde{u}^d=0\quad \text{on } (-1,0]\times B_1' \times \{0\}.\]

Choose a mollifier $\varphi_\varepsilon$ which is symmetric with respect to the $x_d$ variable, i.e., $\varphi(s,y',-y_d)=\varphi(s,y',y_d)$. Define
\[
     \tilde{u}^{(\varepsilon)}(t,x',x_d)=\int_{Q_\varepsilon} \varphi_\varepsilon(s,y',y_d)\tilde{u}(t-s,x-y',x_d-y_d)\myd{y}{ds}
\]
for $(t,x) \in (-1+\varepsilon^2,0)\times B_{1-\varepsilon}$, $0<\varepsilon<1$.

By a change of variables, {we have $(\tilde{u}^d)^{(\varepsilon)}(t,x',0)=0$} since $\varphi_\varepsilon$ is a symmetric mollifier with respect to the $x_d$ variable.

Similarly, we have $D_d (\tilde{u}^{(\varepsilon)})^k(t,x',0)=0$ for $k=1,2,\dots,d-1$. Hence $\tilde{u}^{(\varepsilon)}$ satisfies the Lions boundary conditions. Now we mollify the equation and write
\[
   \partial_t \tilde{u}^{(\varepsilon)}-\overline{a}^{ij}(t,x)D_{ij} \tilde{u}^{(\varepsilon)}+\nabla \tilde{p}^{(\varepsilon)}=\tilde{f}^{(\varepsilon)} + h^\varepsilon \quad \text{on } Q_1,
\]
where
\[   h^\varepsilon(t,x)=[\overline{a}^{ij}(t,x)D_{ij} \tilde{u}]^{(\varepsilon)}-\overline{a}^{ij}(t,x)D_{ij}\tilde{u}^{(\varepsilon)}.\]

To apply Theorem \ref{thm:boundary-solvability}, we need to extend $\overline{a}^{ij}$ to the whole space so that the extended one satisfies Assumption \ref{assump:VMO} $(\delta)$. Since $a^{ij}$ satisfies Assumption \ref{assump:VMO-boundary} $(\delta)$, there exists $0<R_0<1/4$ such that for each $(t_0,x_0) \in \overline{Q_2^+}$ and for all $0<r<R_0$, there exists $\hat{a}^{ij}(t)$ satisfying uniform ellipticity conditions \eqref{eq:elliptic} such that
\[
   \fint_{Q_r^+(t_0,x_0)} |\overline{a}^{ij}(t,x)-\hat{a}^{ij}(t)|\myd{x}dt\leq \delta
\]
 Choose $\eta \in C_0^\infty(B_{7/4})$ satisfying $\eta=1$ in $B_{5/4}$ and define
\[
     \tilde{a}^{ij}(t,x)=\overline{a}^{ij}(t,x)\eta(x)+\delta^{ij}(1-\eta(x)).
\]
Then $\tilde{a}^{ij}$ is bounded and uniformly elliptic. By extending $\tilde{a}^{ij}$ periodically in $t$ if necessary, a direct computation shows that there exists $0<R_1<R_0$ depending only on $d$, $\delta$, $\nu$, $R_0$ such that for any $(t_0,x_0)\in\mathbb{R}^{d+1}$ and for all $0<r<R_1$, we have
\[
    \fint_{Q_r(t_0,x_0)} |\tilde{a}^{ij}-(\tilde{a}^{ij})_{B_r(x_0)}(t)|\myd{x}dt\leq 4\delta.
\]

By Theorem \ref{thm:boundary-solvability}, there exists $(u_1^\varepsilon, p_1^\varepsilon)$ satisfying $(u_1^\varepsilon,\nabla p_1^\varepsilon)\in\oSob{1,2}{q_0}((-1,0)\times\mathbb{R}^d_+)^d\times \Leb{q_0}((-1,0)\times\mathbb{R}^d_+)^d$, with $u_1^\varepsilon(-1,\cdot)=0$ on $\mathbb{R}^d_+$,
\[
     \partial_t u_1^{\varepsilon}-\tilde{a}^{ij}D_{ij} {u}_1^\varepsilon + \nabla {p}_1^\varepsilon=h^{\varepsilon}1_{Q_{3/4}^+},\quad \Div \tilde{u}_1^\varepsilon=0\quad \text{in } (-1,0)\times \mathbb{R}^d_+{,}
\]
and
\[
    D_d (u_1^\varepsilon)^k = (u_1^\varepsilon)^d=0\quad \text{on } (-1,0)\times \mathbb{R}^{d-1}\times\{0\}
\]
for $k=1,\dots, d-1$. Moreover, we have
\[  \norm{u_1^\varepsilon}{\Sob{1,2}{q_0}((-1,0)\times\mathbb{R}^d_+)}\leq N \norm{h^\varepsilon}{\Leb{q_0}(Q_{3/4}^+)}\]
for some constant $N$ independent of $\varepsilon$.

Define $u_2^\varepsilon=\tilde{u}^{(\varepsilon)}-u_1^\varepsilon$ and $p_2^\varepsilon=\tilde{p}^{(\varepsilon)}-p_1^\varepsilon$. Then using Theorem {\ref{thm:boundary-solvability}} as in Lemma \ref{lem:regularity-Stokes}, one can prove that $(u_2^\varepsilon,p_2^\varepsilon) \in \tilde{W}^{1,2}_{s,q}(Q_{3/4}^+)^d\times \Sob{0,1}{1}(Q_{3/4}^+)$ is a strong solution to \eqref{eq:Stokes-nondiv} in $Q_{3/4}^+$ satisfying the Lions boundary conditions on $(-(3/4)^2,0]\times B_{3/4}' \times\{0\}$, $\tilde{f}^{(\varepsilon)}$, and $\tilde{g}^{(\varepsilon)}$ instead of $f$ and $g$ by following the proof of Lemma \ref{lem:regularity-Stokes}. Then using a similar idea as in the proof of Theorem \ref{thm:B}, we can prove the desired estimates by using Dong-Kim-Phan \cite[Theorem 1.2]{DKP22} instead of Dong-Phan (Theorem \ref{thm:DP}). We leave the details to interested readers.

\appendix

\section{Proof of Theorem \ref{thm:mixed-simple} }\label{app:B}
This section is devoted to a proof of Theorem \ref{thm:mixed-simple} which concerns the solvability of Stokes equations in nondivergence form with simple coefficients. 

\begin{proof}
The proof is almost identical to that of Theorem 1.4 in \cite{DKP22} { which uses Proposition 4.2 in the same paper.} The proof of existence part is split into five steps. From Step 1 to Step 3, the key differences {between ours and Proposition 4.2 in \cite{DKP22}} are using Lemma \ref{lem:potential} instead of Lemma 4.1 in \cite{DKP22} and Theorem \ref{thm:classical-Lp} when we construct vorticity from given {external} force.  Following the construction of pressure in Step 4 of {Proposition 4.2} in \cite{DKP22}, if we define $h=f-\partial_t u+a^{ij} D_{ij}u$, then one can show that
\[  |\nabla p^\varepsilon(t,x)|\leq N(Mh)(t,x) \]
for some constant $N=N(d)>0$ and for all $\varepsilon>0$ and $(t,x)\in\mathbb{R}^d_T$. Hence it follows from Lemma \ref{lem:HL} that
\[	\norm{\nabla p^\varepsilon}{\Leb{s,q,w}(\mathbb{R}^d_T)}\leq N(d,s,q,K_0)\norm{h}{\Leb{s,q,w}(\mathbb{R}^d_T)} \]
is bounded uniformly in $\varepsilon>0$. {By subtracting a function of $t$, we may assume that $(p^\varepsilon(\cdot,t))_{B_1}={0}$ for $t\in (0,T)$.
Then} for each $R>1$, it follows from the Poincar\'e inequality (Lemma \ref{lem:weighted-Poincare}) that
\[
              \norm{ {p}^\varepsilon(t,\cdot)}{\Leb{q,w_1}(B_R)}\leq N(d,q,K_0,R)\norm{\nabla  {p}^\varepsilon(t,\cdot)}{\Leb{q,w_1}(B_R)}
\]
for each $t\in[0,T]$. By taking $\Leb{s,w_2}$-norm, we get
\[
         \norm{ {p}^\varepsilon}{\Leb{s,q,w}((0,T)\times B_R)}\leq N(d,s,q,K_0,R) \norm{h}{\Leb{s,q,w}((0,T)\times B_R)},
\]
which is uniformly bounded in $\varepsilon$. Hence by weak compactness results in weighted $\Leb{q}$-spaces, one can conclude that there exists a locally integrable function $p : \mathbb{R}^d_T\rightarrow\mathbb{R}$ such that $\nabla p \in \Leb{s,q,w}(\mathbb{R}^d_T)^d$ and $(u,p)$ satisfies equation \eqref{eq:Stokes-nondiv-simple} in $\mathbb{R}^d_T$.

{Although Step 5 is also similar to that of the proof of Theorem 1.4 in \cite{DKP22}, we give a detailed proof for the sake of convenience.}

\noindent \emph{Step 5}. Since $C_0^\infty(\mathbb{R}^d_T)$ is dense in $\Leb{s,q,w}(\mathbb{R}^d_T)$, we need to show that there exist $g^m$ and $G^m$ vanish for large $|x|$ uniformly in $t \in [0,T]$,
\[	
g^m \in \odSob{1}{s,q,w}(\mathbb{R}^d_T),\quad G^m \in \Leb{s,q,w}(\mathbb{R}^d_T)^d,\quad \partial_t g^m = \Div G^m\quad \text{in } \mathbb{R}^d_T, \]
and
\[	\norm{g-g^m}{\Leb{s,q,w}(\mathbb{R}^d_T)}+\norm{Dg-Dg^m}{\Leb{s,q,w}(\mathbb{R}^d_T)}+\norm{G-G^m}{\Leb{s,q,w}(\mathbb{R}^d_T)}\rightarrow 0\]
as $m\rightarrow\infty$.

Choose a sequence of smooth functions $\{\chi_m\}$ on $\mathbb{R}^d$ such that $\chi_m=1$ on $B_{m/2}$ and $\chi_m=0$ outside $B_m$, $m=1,2,3,\dots$. Define
\[
   c_m(t):=\frac{\int_{B_m} \nabla \chi_m(y) \cdot G(t,y) dy}{\int_{B_m} \chi_m(y) dy}.
\]
Note that
\[  \int_{B_m} (-\nabla \chi_m \cdot G+c_m(t) \chi_m(x))\myd{x}=0 \]
for a.e. $t\in (0,T)$. Hence by Theorem \ref{thm:div-equation}, using the integral representation of solutions, we can find $H^m$ in $(0,T)\times B_m$ such that
\[
  \left\{
  \begin{alignedat}{2}
    \Div H^m& =-\nabla \chi_m \cdot G+c_m(t) \chi_m(x)&&\quad \text{in } (0,T)\times B_m,\\
   H^m& =0&&\quad \text{on } (0,T)\times \partial B_m,
   \end{alignedat}
  \right.
\]
and
\begin{equation}\label{eq:DxH}
\begin{aligned}
&\norm{D H^m}{\Leb{s,q,w}((0,T)\times B_m)}\\
&\leq N(d,s,q,K_0)\left(\norm{\nabla \chi_m \cdot G}{\Leb{s,q,w}((0,T)\times B_m)}+\norm{c_m(t)\chi_m(x)}{\Leb{s,q,w}((0,T)\times B_m)} \right).
\end{aligned}
\end{equation}
By H\"older's inequality and the $A_q$-condition, we have
\begin{align*}
&|c_m(t)|\norm{\chi_m}{\Leb{q,w_1}(B_m)} \\
&\leq \frac{1}{\left|\int_{B_m} \chi_m(y)dy\right|} \left(\int_{B_m} |\nabla \chi_m(y)| |G(t,y)|dy\right)\left(\int_{B_m} |\chi_m(x)|^q w_1(x)\myd{x}\right)^{1/q}\\
&\leq \frac{\norm{(\nabla \chi_m)  G(t)}{\Leb{q,w_1}(B_m)}}{\left|\int_{B_m} \chi_m(y)dy\right|} \left(\int_{B_m} w_1(y)^{-1/(q-1)} dy \right)^{(q-1)/q}\left(\int_{B_m} |\chi_m(x)|^q w_1(x)\myd{x}\right)^{1/q}\\
&\leq \frac{|B_m|}{\left|\int_{B_m} \chi_m(y)dy\right|}[w_1]_{A_q}^{1/q}\norm{(\nabla \chi_m) G(t)}{\Leb{q,w_1}(B_m)}\leq N [w_1]_{A_q}^{1/q}\norm{(\nabla \chi_m) G(t)}{\Leb{q,w_1}(B_m)},
\end{align*}
where $N$ is independent of $m$. This implies that
\begin{equation}\label{eq:gradG}
\begin{aligned}
&\norm{\nabla \chi_m \cdot G}{\Leb{s,q,w}((0,T)\times B_m)}+\norm{c_m(t) \chi_m(x)}{\Leb{s,q,w}((0,T)\times B_m)}\\
&\leq \norm{\nabla \chi_m \cdot G}{\Leb{s,q,w}((0,T)\times B_m)} + N[w_1]_{A_q}^{1/q} \norm{\nabla \chi_m G}{\Leb{s,q,w}((0,T)\times B_m)}\\
& \leq N m^{-1} \norm{1_{B_m\setminus B_{m/2}} G}{\Leb{s,q,w}(\mathbb{R}^d_T)},
\end{aligned}
\end{equation}
where $N$ is independent of $m$. Hence it follows from \eqref{eq:DxH}, \eqref{eq:gradG}, the weighted Poincar\'e inequality (Lemma \ref{lem:weighted-Poincare}) on $B_m$ and the fact that $\norm{1_{B_m\setminus B_{m/2}} G}{\Leb{s,q,w}(\mathbb{R}^d_T)}\rightarrow 0$ as $m\rightarrow\infty$ that
\begin{align*}
   \norm{H^m}{\Leb{s,q,w}((0,T)\times B_m)} &\leq Nm \norm{D H^m}{\Leb{s,q,w}((0,T)\times B_m)} \\
   &\leq N \norm{1_{B_m\setminus B_{m/2}} G}{\Leb{s,q,w}(\mathbb{R}^d_T)}\rightarrow 0
\end{align*}
as $m\rightarrow\infty$.

Define
\[
g^m(t,x):=\chi_m(x)g(t,x)+\chi_m(x)\int_0^t c_m(s)ds
\]
and
\[
G^m(t,x):=\begin{cases}
\chi_m(x) G(t,x) +H^m(t,x)&\quad \text{in } (0,T)\times B_m,\\
0&\quad \text{in } (0,T)\times (\mathbb{R}^d\setminus B_m).\\
\end{cases}
\]
Then it is easy to see that $g^m \in \odSob{1}{s,q,w}(\mathbb{R}^d_T)$ and $\partial_t g^m(t,x)= \Div G^m(t,x)$ in the sense of \eqref{eq:compatibility-nondiv}, and the rest of the result follows from the dominated convergence theorem. This completes the proof of existence part of Theorem \ref{thm:mixed-simple}.

It remains us to show the uniqueness part. We first take the curl operation to the equation in the weak sense. Then $\omega_{kl}=D_l u^k-D_k u^l \in \Leb{s,q,w}(\mathbb{R}^d_T)$ is a very weak solution to the heat equation with simple coefficients, i.e.,
\begin{equation}\label{eq:duality-vorticity}
\int_{\mathbb{R}^d_T} \omega_{kl} (\partial_t \psi + a^{ij}(t)D_{ij}\psi)\myd{xdt}=0
\end{equation}
for all $\psi \in C_0^\infty ([0,T)\times\mathbb{R}^d)$. By a standard density argument, the identity holds for all $\psi \in \Sob{1,2}{s',q',\tilde{w}}(\mathbb{R}^d_T)$ with $\psi(T,x)=0$, where $\tilde{w}=w_1^{-1/(q-1)} w_2^{-1/(s-1)}$. By Theorem \ref{thm:classical-Lp} (i), given $\varphi \in C_0^\infty(\mathbb{R}^d_T)$, there exists a unique $\psi \in \Sob{1,2}{s',q',\tilde{w}}(\mathbb{R}^d_T)$ with $\psi(T,x)=0$ such that
\[
    \partial_t \psi + a^{ij}(t)D_{ij} \psi =\varphi\quad \text{in } \mathbb{R}^d_T.
\]
If we put this $\psi$ in \eqref{eq:duality-vorticity}, then we have
\[
\int_{\mathbb{R}^d_T} \varphi\omega_{kl} \myd{xdt}=0
\]
for all $\varphi \in C_0^\infty (\mathbb{R}^d_T)$. Hence $\omega_{kl}$ is identically zero in $\mathbb{R}^d_T$. Since $u \in \oSob{1,2}{s,q,w}(\mathbb{R}^d_T)^d$ satisfies
\[ \Delta u^l = \sum_{k\neq l} D_k (D_k u^l-D_l u^k)=0\quad \text{in } \mathbb{R}^d_T \]
for all $l=1,\dots,d$,  it follows from the mean value property of harmonic functions, H\"older's inequality, and the $A_q$-condition that
\begin{align*}
     |u(t,x)|&\leq  \fint_{B_{R}(x)} |u(t,y)|\myd{y}\\
     & \leq  \frac{1}{|B_R(x)|} |B_R(x)|^{1-1/q} \norm{u(t,\cdot)}{\Leb{q,w_1}(\mathbb{R}^d)} \left(\fint_{B_R(x)} w_1^{-\frac{1}{q-1}} \myd{y}\right)^{1-{1}/{q}} \\
     & \leq  \frac{[w_1]_{A_q}^{1/q}}{w_1(B_R)^{1/q}} \norm{u(t,\cdot)}{\Leb{q,w_1}(\mathbb{R}^d)}
\end{align*}
for a.e. $t\in (0,T)$, for all $x\in \mathbb{R}^d$, and for all $R>0$. Since $w_1(B_R)\rightarrow \infty$ as $R\rightarrow \infty$ (Proposition \ref{prop:weight-property} (vi)), it follows that  $u=0$ for a.e. {on $\mathbb{R}^d_T$} and hence $\nabla p=0$. This completes the proof of Theorem \ref{thm:mixed-simple}.
\end{proof}

\section{Proof of Theorem \ref{thm:mixed-simple-div}}\label{app:C}
This section is devoted to a proof of Theorem \ref{thm:mixed-simple-div} which concerns the solvability of Stokes equations in divergence form with simple coefficients.

\begin{proof}
We first show the existence of weak solutions. Consider
\begin{equation}\label{eq:Stokes-div-1}
\left\{
\begin{alignedat}{2}
\partial_t u_1-D_i(a^{ij}D_j u_1)+\nabla \pi &=\Div \boldF&&\quad \text{in } \mathbb{R}^d_T,\\
\Div u_1 &=0&&\quad \text{in } \mathbb{R}^d_T,\\
u_1&=0 &&\quad \text{on } \{t=0\}\times\mathbb{R}^d
\end{alignedat}
\right.
\end{equation}
and
\begin{equation}\label{eq:Stokes-div-2}
\left\{
\begin{alignedat}{2}
\partial_t u_2-D_i(a^{ij}D_j u_2)+\nabla \tilde{\pi}  &=0&&\quad \text{in } \mathbb{R}^d_T,\\
\Div u_2 &=g&&\quad \text{in } \mathbb{R}^d_T,\\
u_2&=0 &&\quad \text{on } \{t=0\}\times\mathbb{R}^d.
\end{alignedat}
\right.
\end{equation}

Write $\boldF =(f^1,f^2,\dots,f^d)$, where $f^i$ is a vector field. Then by Theorem \ref{thm:mixed-simple}, there exists a strong solution $(v_k,\pi_k)$  satisfying
\[   v_k \in \oSob{1,2}{s,q,w}(\mathbb{R}^d_T)^d,\quad \nabla \pi_k \in \Leb{s,q,w}(\mathbb{R}^d_T)^d,\]
and
\begin{equation}\label{eq:strong-sol-k}
\partial_t v_k-D_i(a^{ij}(t)D_{j}v_k)  +\nabla \pi_k =f^k,\quad \Div v_k=0
\end{equation}
for $k=1,\dots, d$. Moreover, we have
\[
     \norm{D^2v_k}{\Leb{s,q,w}(\mathbb{R}^d_T)}+\norm{\nabla \pi_k}{\Leb{s,q,w}(\mathbb{R}^d_T)}\leq N_1\norm{f^k}{\Leb{s,q,w}(\mathbb{R}^d_T)}
\]
and
\[              \norm{v_k}{\Sob{1,2}{s,q,w}(\mathbb{R}^d_T)}+\norm{\nabla \pi_k}{\Leb{s,q,w}(\mathbb{R}^d_T)}\leq N_2 \norm{f^k}{\Leb{s,q,w}(\mathbb{R}^d_T)}\]
for some constants $N_1=N_1(d,s,q,K_0,\nu)>0$, $N_2=N_2(d,s,q,K_0,\nu,T)>0$ and for all $k=1,\dots,d$.

Define $u_1=\sum_{k=1}^d D_k v_k$ and $\pi=\sum_{k=1}^d D_k \pi_k$. Then $(u_1,\pi) \in \odSob{1}{s,q,w}(\mathbb{R}^d_T)^d\times \Leb{s,q,w}(\mathbb{R}^d_T)$ is a weak solution of \eqref{eq:Stokes-div-1}. Indeed, since $(v_k,\pi_k)$ is a strong solution of \eqref{eq:strong-sol-k}, we have
\begin{align*}
-\int_{\mathbb{R}^d_T} v_k \cdot \partial_t \phi \myd{x}dt+\int_{\mathbb{R}^d_T} (a^{ij}(t)D_j v_k) \cdot D_i \phi - \pi_k \Div \phi \myd{x}dt=\int_{\mathbb{R}^d_T} f^k \cdot \phi \myd{x}dt
\end{align*}
for all $\phi \in C_0^\infty([0,T)\times\mathbb{R}^d)^d$. For $\psi \in C_0^\infty([0,T)\times\mathbb{R}^d)^d$, put $\phi=D_k \psi$ in the identity. Then we have
\[
\begin{aligned}
&-\int_{\mathbb{R}^d_T} (D_k v_k)\partial_t \psi \myd{x}dt+\int_{\mathbb{R}^d_T} (a^{ij}(t) D_{jk} v_k) \cdot D_i \psi  -(D_k \pi_k) \Div \psi \myd{x}dt\\
&=-\int_{\mathbb{R}^d_T} f^k \cdot D_k \psi \myd{x}dt
\end{aligned}
\]
for $k=1,\dots, d$. Hence by summing it over $k$, we get
\[
-\int_{\mathbb{R}^d_T} u_1\cdot \partial_t\psi \myd{x}dt+\int_{\mathbb{R}^d_T} (a^{ij}D_j v )\cdot D_i  \psi -\pi\Div \psi \myd{x}dt=-\int_{\mathbb{R}^d_T} \boldF : \nabla \psi \myd{x}dt
\]
for all $\psi \in C_0^\infty([0,T)\times\mathbb{R}^d)^d$.
Moreover, it follows from \eqref{eq:simple-Hessian} and \eqref{eq:simple-norm-estimate} that
\begin{align}\label{eq:gradient-v}
\norm{D u_1}{\Leb{s,q,w}(\mathbb{R}^d_T)}+\norm{\pi}{\Leb{s,q,w}(\mathbb{R}^d_T)}&\leq N_1 \sum_{k=1}^d \norm{D^2 v_k}{\Leb{s,q,w}(\mathbb{R}^d_T)} \\
&\leq N_1\sum_{k=1}^d \norm{f_k}{\Leb{s,q,w}(\mathbb{R}^d_T)}\leq N_1 \norm{\boldF}{\Leb{s,q,w}(\mathbb{R}^d_T)}\nonumber
\end{align}
and
\begin{align}
\norm{u_1}{\dSob{1}{s,q,w}(\mathbb{R}^d_T)}+\norm{\pi}{\Leb{s,q,w}(\mathbb{R}^d_T)}&\leq N_2 \sum_{k=1}^d \left(\norm{v_k}{\Sob{1,2}{s,q,w}(\mathbb{R}^d_T)}+\norm{\nabla \pi_k}{\Leb{s,q,w}(\mathbb{R}^d_T)} \right)\label{eq:div-sob-v}\\
&\leq  N_2 \norm{\boldF}{\Leb{s,q,w}(\mathbb{R}^d_T)}\nonumber
\end{align}
for some constants $N_1=N_1(d,s,q,K_0,\nu)>0$ and $N_2=N_2(d,s,q,K_0,\nu,T)>0$.

To find a solution $(w,\tilde{\pi})$ to \eqref{eq:Stokes-div-2}, define
\[    \tilde{\pi}=\sum_{i,j=1}^d \mathcal{R}_i \mathcal{R}_j (G^{ij}-ga^{ij}(t)),\]
where $\mathcal{R}_j$ denotes the $j$th Riesz transform. Then by $\Leb{q,w_1}$-boundedness of Riesz transforms (see e.g. \cite[\S 4.2, Chapter V]{S93}), we have  $\tilde{\pi}\in \Leb{s,q,w}(\mathbb{R}^d_T)$ and
\[     \norm{\tilde{\pi}}{\Leb{s,q,w}(\mathbb{R}^d_T)}\leq N( \norm{\boldG}{\Leb{s,q,w}(\mathbb{R}^d_T)}+\norm{g}{\Leb{s,q,w}(\mathbb{R}^d_T)})
\]
for some constant $N=N(d,s,q,K_0,\nu)>0$.   Since
\[
   -\mathcal{R}_i \mathcal{R}_j (\Delta \psi)=D_{ij} \psi
\]
for all $\psi \in C_0^\infty(\mathbb{R}^d)$ and $\mathcal{R}_i \mathcal{R}_j$ is self-adjoint on $\Leb{2}$, it follows that
\begin{equation}\label{eq:pressure-identity}
\begin{aligned}
-    \int_{\mathbb{R}^d} \tilde{\pi}(t,x)\Delta \psi(x) \myd{x}&=\int_{\mathbb{R}^d} (G^{ij}(t,x)-a^{ij}(t)g(t,x))D_{ij}\psi(x)\myd{x}
\end{aligned}
\end{equation}
for all $\psi \in C_0^\infty(\mathbb{R}^d)$ and for a.e. $t\in (0,T)$.

By \eqref{eq:pressure-identity} and the compatibility condition \eqref{eq:compatibility-divergence},  the identity
\begin{equation}\label{eq:pressure-identity-2}
-\int_{\mathbb{R}^d_T} \tilde{\pi} \Delta \psi \myd{x}dt =-\int_{\mathbb{R}^d_T} g (\psi_t+a^{ij}(t)D_{ij}\psi) \myd{x}dt
\end{equation}	
holds for all $\psi \in C_0^\infty([0,T)\times \mathbb{R}^d)$.

On the other hand, it follows from Theorem \ref{thm:classical-Lp} (i) that there exists a unique $\Phi \in \oSob{1,2}{s,q,w}(\mathbb{R}^d_T)$ satisfying
\begin{equation}\label{eq:duality-problem}
   \partial_t \Phi -a^{ij}(t)D_{ij} \Phi =\tilde{\pi}\quad{\text{in $\mathbb{R}^d_T$}.}
\end{equation}
Moreover, we have
\[    \norm{D^2 \Phi}{\Leb{s,q,w}(\mathbb{R}^d_T)}\leq N_1 \norm{\tilde{\pi}}{\Leb{s,q,w}(\mathbb{R}^d_T)}\]
and
\[    \norm{\Phi}{\Sob{1,2}{s,q,w}(\mathbb{R}^d_T)}\leq N_2 \norm{\tilde{\pi}}{\Leb{s,q,w}(\mathbb{R}^d_T)}\]
for some constants $N_1=N(d,s,q,\nu,K_0)>0$ and $N_2=N(d,s,q,\nu,K_0,T)>0$.

We show that $-\Delta \Phi=g$. Since $\tilde{\pi}$ satisfies \eqref{eq:pressure-identity-2} and $\Phi$ satisfies \eqref{eq:duality-problem}, we have
\begin{align*}
\int_{\mathbb{R}^d_T} (\partial_t \Phi-a^{ij}(t)D_{ij} \Phi)\Delta \psi \myd{x}dt&=\int_{\mathbb{R}^d_T} g(\partial_t \psi+a^{ij}(t)D_{ij}\psi)\myd{x}dt
\end{align*}
for all $\psi \in C_0^\infty([0,T)\times \mathbb{R}^d)$. Integration by part gives
\begin{align*}
  \int_{\mathbb{R}^d_T} (\partial_t \Phi-a^{ij}(t)D_{ij} \Phi)\Delta \psi \myd{x}dt&=\int_{\mathbb{R}^d_T} \Phi \Delta (-\partial_t \psi - a^{ij}(t)D_{ij}\psi) \myd{x}dt\\
  &=-\int_{\mathbb{R}^d_T} \Delta \Phi (\partial_t \psi+a^{ij}(t)D_{ij} \psi)\myd{x}dt
\end{align*}
for all $\psi \in C_0^\infty([0,T)\times\mathbb{R}^d)$. Hence
\begin{equation}\label{eq:duality-identity-3}
   -\int_{\mathbb{R}^d_T} (\Delta\Phi)  (\partial_t \psi +a^{ij}(t)D_{ij} \psi) \myd{x}dt=\int_{\mathbb{R}^d_T} g (\partial_t\psi +a^{ij}(t)D_{ij} \psi)\myd{x}dt
\end{equation}
for all $\psi \in C_0^\infty([0,T)\times \mathbb{R}^d)$.  Then by a standard density argument, we see that the identity holds for all $\psi \in \Sob{1,2}{s',q',\tilde{w}}(\mathbb{R}^d_T)$ with $\psi(T,x)=0$, where $\tilde{w}=w_1^{-1/(q-1)}w_2^{-1/(s-1)}$.

Given $\varphi \in C_0^\infty(\mathbb{R}^d_T)$, it follows from Theorem \ref{thm:classical-Lp} (i) that there exists a unique $\psi \in \Sob{1,2}{s',q',\tilde{w}}(\mathbb{R}^d_T)$ satisfying $\psi(T,x)=0$ and
\[   \partial_t \psi+a^{ij}(t)D_{ij} \psi =\varphi\quad\text{in } \mathbb{R}^d_T.    \]
Hence by \eqref{eq:duality-identity-3}, we have
\[   -\int_{\mathbb{R}^d_T} (\Delta \Phi)\varphi \myd{x}dt=\int_{\mathbb{R}^d_T} g \varphi \myd{x}dt \]
for all $\varphi \in C_0^\infty([0,T)\times\mathbb{R}^d)$, which implies that $-\Delta \Phi=g$ in $\mathbb{R}^d_T$. 	
Put $u_2=-\nabla \Phi$. Then by \eqref{eq:duality-problem}, it is easy to show that $(u_2,\tilde{\pi})$ is a weak solution to \eqref{eq:Stokes-div-2} satisfying
\begin{equation}\label{eq:gradient-w}
\begin{aligned}
\norm{Du_2}{\Leb{s,q,w}(\mathbb{R}^d_T)}+\norm{\tilde{\pi}}{\Leb{s,q,w}(\mathbb{R}^d_T)}    &\leq N_1\norm{\tilde{\pi}}{\Leb{s,q,w}(\mathbb{R}^d_T)}\\
    &\leq N_1\left(\norm{\boldG}{\Leb{s,q,w}(\mathbb{R}^d_T)}+\norm{g}{\Leb{s,q,w}(\mathbb{R}^d_T)}\right)
\end{aligned}
\end{equation}
and
\begin{equation}\label{eq:div-sob-w}
  \norm{u_2}{\dSob{1}{s,q,w}(\mathbb{R}^d_T)}+\norm{\tilde{\pi}}{\Leb{s,q,w}(\mathbb{R}^d_T)}\leq N_2\left(\norm{\boldG}{\Leb{s,q,w}(\mathbb{R}^d_T)}+\norm{g}{\Leb{s,q,w}(\mathbb{R}^d_T)} \right)
\end{equation}
for some constants $N_1=N_1(d,s,q,\nu,K_0)>0$ and $N_2=N_2(d,s,q,\nu,K_0,T)>0$.

{
Since $u_2=-\nabla \Phi$ and $-\Delta \Phi = g$, it follows that $u_2 \in \odSob{1}{s,q,w}(\mathbb{R}^d_T)$ and
\begin{align*}
        \int_{\mathbb{R}^d} \nabla u_2^l \cdot \nabla \phi \myd{x}&=-\int_{\mathbb{R}^d} \nabla (D_l \Phi) \cdot \nabla \phi \myd{x}\\
        &=\int_{\mathbb{R}^d} \nabla \Phi \cdot \nabla (D_l \phi)\myd{x}=  -\int_{\mathbb{R}^d} g(D_l \phi)\myd{x}
\end{align*}
for all $\phi \in C_0^\infty(\mathbb{R}^d)$ and for a.e. $t\in (0,T)$. Hence it follows from Corollary \ref{cor:a-priori-regularity} that
\begin{equation}\label{eq:gradient-w-2}
\norm{Du_2}{\Leb{s,q,w}(\mathbb{R}^d_T)}\leq N_1\norm{g}{\Leb{s,q,w}(\mathbb{R}^d_T)}
\end{equation}
for some constant $N_1=N_1(d,s,q,\nu,K_0)>0$.

Define $u=u_1+u_2$ and $p=\pi+\tilde{\pi}$. Then $(u,p)$ is a weak solution to \eqref{eq:Stokes-div-simple} in $\mathbb{R}^d_T$ satisfying $u\in \odSob{1}{s,q,w}(\mathbb{R}^d_T)^d$ and $p\in \Leb{s,q,w}(\mathbb{R}^d_T)$. By \eqref{eq:gradient-v} and \eqref{eq:gradient-w-2}, we have
\begin{align*}
\norm{Du}{\Leb{s,q,w}(\mathbb{R}^d_T)}&\leq  \norm{Du_1}{\Leb{s,q,w}(\mathbb{R}^d_T)}+\norm{Du_2}{\Leb{s,q,w}(\mathbb{R}^d_T)} \\
&\leq N_1 \left( \norm{\boldF}{\Leb{s,q,w}(\mathbb{R}^d_T)}+\norm{g}{\Leb{s,q,w}(\mathbb{R}^d_T)}\right)
\end{align*}
for some constant $N_1=N_1(d,s,q,\nu,K_0)>0$.}

Similarly, by \eqref{eq:gradient-v} and \eqref{eq:gradient-w}, we have
\begin{align*}
  \norm{p}{\Leb{s,q,w}(\mathbb{R}^d_T)}&\leq \norm{Du_1}{\Leb{s,q,w}(\mathbb{R}^d_T)}+\norm{p_1}{\Leb{s,q,w}(\mathbb{R}^d_T)}\\
       &\relphantom{=}+\norm{Du_2}{\Leb{s,q,w}(\mathbb{R}^d_T)}+\norm{p_2}{\Leb{s,q,w}(\mathbb{R}^d_T)}\\
        &\leq N_1\left(\norm{\boldF}{\Leb{s,q,w}(\mathbb{R}^d_T)}+\norm{\boldG}{\Leb{s,q,w}(\mathbb{R}^d_T)}+\norm{g}{\Leb{s,q,w}(\mathbb{R}^d_T)}\right).
\end{align*}
Moreover, it follows from \eqref{eq:div-sob-v} and \eqref{eq:div-sob-w} that
\[
\norm{u}{\dSob{1}{s,q,w}(\mathbb{R}^d_T)}+\norm{p}{\Leb{s,q,w}(\mathbb{R}^d_T)}\leq N_2 \left(\norm{\boldF}{\Leb{s,q,w}(\mathbb{R}^d_T)}+\norm{\boldG}{\Leb{s,q,w}(\mathbb{R}^d_T)}+\norm{g}{\Leb{s,q,w}(\mathbb{R}^d_T)}\right)
\]
for some constants $N_1=N_1(d,s,q,\nu,K_0)>0$ and $N_2=N_2(d,s,q,\nu,K_0,T)>0$.

{
 It remains to show the uniqueness of weak solutions. Suppose that $(u,p)$ satisfies
 \[     u\in \odSob{1}{s,q,w}(\mathbb{R}^d_T)^d,\quad p \in \Leb{s,q,w}(\mathbb{R}^d_T), \]
 and
\begin{equation}\label{eq:uniqueness-function}
   \int_{\mathbb{R}^d_T} u\cdot(\partial_t \phi + a^{ij}(t)D_{ij} \phi)+p\Div\phi \myd{x}dt=0
\end{equation}
for all $\phi \in C_0^\infty([0,T)\times\mathbb{R}^d)^d$. For $\psi \in C_0^\infty(\mathbb{R}^d_T)$, put $\phi=\nabla \psi$ in \eqref{eq:uniqueness-function}. Since $\Div u=0$ in $\mathbb{R}^d_T$, we get
\[   \int_{\mathbb{R}^d_T} p \Delta \psi\myd{x}dt=0 \]
for all $\psi \in C_0^\infty(\mathbb{R}^d_T)$. This implies that $p$ is  harmonic in $\mathbb{R}^d$ a.e. $t\in (0,T)$. Then following exactly the same argument as in the proof of uniqueness part of Theorem \ref{thm:mixed-simple}, one can show that $p$ is identically zero. By \eqref{eq:uniqueness-function}, $u\in \odSob{1}{s,q,w}(\mathbb{R}^d_T)^d$ is a weak solution to
\[    \partial_t u -D_i(a^{ij}D_j u)=0\quad \text{in } \mathbb{R}^d_T.\]
Therefore, it follows from Theorem \ref{thm:classical-Lp} (ii) that $u$ is identically zero, which completes the proof of Theorem \ref{thm:mixed-simple-div}. }
\end{proof}

\section*{Acknowledgements}
H. Dong was partially supported by Simons Fellows Award 007638 and the NSF under agreement DMS-2055244. H. Kwon was partially supported by the NSF under agreement DMS-2055244.

\bibliographystyle{amsplain}
\providecommand{\bysame}{\leavevmode\hbox to3em{\hrulefill}\thinspace}
\providecommand{\MR}{\relax\ifhmode\unskip\space\fi MR }
\providecommand{\MRhref}[2]{%
  \href{http://www.ams.org/mathscinet-getitem?mr=#1}{#2}
}
\providecommand{\href}[2]{#2}

\end{document}